\documentclass[preprint,12pt,3p]{elsarticle}
\pdfoutput=1

\newif\ifAnimations
\Animationsfalse

\usepackage{graphicx}
\usepackage{epstopdf}
\graphicspath{{./Figures/}}
\usepackage{tipa}

\usepackage[cmex10]{amsmath}
\usepackage{amsthm}
\usepackage{amssymb}
\usepackage{amsfonts}
\usepackage{textcomp}
\usepackage{bm}
\usepackage{cancel}
\usepackage{upgreek}
\usepackage{algorithm}
\usepackage{algorithmicx}
\usepackage{algpseudocode}
\usepackage{mathrsfs}
\usepackage{marvosym}

\usepackage{tabularx}

\usepackage{multicol}
\usepackage{float}
\usepackage{url}
\usepackage[tight,footnotesize]{subfigure}
\usepackage{enumitem}
\usepackage{calc}
\usepackage{tikz}

\usepackage[above,below]{placeins}

\SetLabelAlign{parright}{\parbox[t]{\labelwidth}{\raggedleft#1}}

\usepackage{booktabs}
\newcommand{\head}[1]{{\small\textbf{#1}}}

\usepackage{color}
\definecolor{MyGrey}{RGB}{234,234,234}
\definecolor{MyBlue}{RGB}{0,0,255}
\definecolor{DarkSlateGray}{rgb}{0.1836,  0.3086  ,  0.3086}
\definecolor{Chocolate}{rgb}{0.8203                                                                                                                          ,  0.4102  ,  0.1172}
\definecolor{FireBrick}{rgb}{ 0.6953   ,  0.1328  ,  0.1328}
\definecolor{MediumBlue}{rgb}{ 0       ,  0       ,  0.8008}
\definecolor{MyGreen}{rgb}{ 0          ,  0.5000  ,  0     }
\definecolor{MyRed}{rgb}{ 1        ,  0.0000  ,  0     }
\definecolor{MyPurple}{rgb}{0.5     ,    0  ,0.5}
\definecolor{MyGrey}{rgb}{0.5,0.5,0.5}

\definecolor{myOrange}{rgb}{1,0.64453,0}
\definecolor{myCyan}{rgb}{0,1,1}

\usepackage{nameref}
\makeatletter
\let\orgdescriptionlabel\descriptionlabel
\renewcommand*{\descriptionlabel}[1]{%
  \let\orglabel\label
  \let\label\@gobble
  \phantomsection
  \edef\@currentlabel{#1}%
  \let\label\orglabel
  \orgdescriptionlabel{#1}%
}
\makeatother

\def\Xint#1{\mathchoice
   {\XXint\displaystyle\textstyle{#1}}%
   {\XXint\textstyle\scriptstyle{#1}}%
   {\XXint\scriptstyle\scriptscriptstyle{#1}}%
   {\XXint\scriptscriptstyle\scriptscriptstyle{#1}}%
   \!\int}
\def\XXint#1#2#3{{\setbox0=\hbox{$#1{#2#3}{\int}$}
     \vcenter{\hbox{$#2#3$}}\kern-.625\wd0}}
\def\dashint{\Xint-}
\def\XintS#1{\mathchoice
   {\XXintS\displaystyle\textstyle{#1}}%
   {\XXintS\textstyle\scriptstyle{#1}}%
   {\XXintS\scriptstyle\scriptscriptstyle{#1}}%
   {\XXintS\scriptscriptstyle\scriptscriptstyle{#1}}%
   \!\int}
\def\XXintS#1#2#3{{\setbox0=\hbox{$#1{#2#3}{\int}$}
     \vcenter{\hbox{$#2#3$}}\kern-.5\wd0}}
\def\dashintS{\XintS-} 

\let\oldsqrt\sqrt
\def\sqrt{\mathpalette\DHLhksqrt}
\def\DHLhksqrt#1#2{%
\setbox0=\hbox{$#1\oldsqrt{#2\,}$}\dimen0=\ht0
\advance\dimen0-0.2\ht0
\setbox2=\hbox{\vrule height\ht0 depth -\dimen0}%
{\box0\lower0.4pt\box2}}

\newcommand\solidSrule[1][1cm]{\rule[0.5ex]{#1}{.75pt}}
\newcommand\solidMrule[1][1cm]{\rule[0.5ex]{#1}{1.25pt}}
\newcommand\solidLrule[1][1cm]{\rule[0.5ex]{#1}{2.25pt}}
\newcommand\dashedrule{\mbox{\solidSrule[.75mm]\hspace{0.25mm}\solidSrule[0.75mm]\hspace{0.25mm}\solidSrule[0.75mm]}}

\newcommand{\tikzcircle}[2][red,fill=red]{\tikz[baseline=-0.65ex]{\draw[#1,radius=#2] (0,0) circle ;}}%
\newcommand{\tikzsquare}{\tikz[baseline=-0.3ex]{\draw[fill=red, red] (0,0) rectangle (1mm,1mm);}}%

\newcommand{\lineDotBlue}{\mbox{\textcolor{MyBlue}{\solidMrule[1.5mm]}\rlap{\textcolor{MyBlue}{\solidMrule[3mm]}}\tikzcircle[MyBlue, fill=MyBlue]{.6mm}\hspace{0.5mm}\textcolor{MyBlue}{\solidMrule[.9mm]}}}
\newcommand{\lineDotBlueS}{\mbox{\textcolor{MyBlue}{\solidMrule[1.5mm]}\rlap{\textcolor{MyBlue}{\solidMrule[3mm]}}\tikzcircle[MyBlue, fill=MyBlue]{.6mm}\hspace{0.5mm}\textcolor{MyBlue}{\solidMrule[1mm]}}}
\newcommand{\lineSqrRed}{\mbox{\textcolor{MyRed}{\solidMrule[1.5mm]}\rlap{\textcolor{MyRed}{\solidMrule[2.7mm]}}\tikzsquare\hspace{0.5mm}\textcolor{MyRed}{\solidMrule[1mm]}}}
\newcommand{\lineSqrRedS}{\mbox{\textcolor{MyRed}{\solidMrule[1.5mm]}\rlap{\textcolor{MyRed}{\solidMrule[2.7mm]}}\tikzsquare\hspace{0.5mm}\textcolor{MyRed}{\solidMrule[1mm]}}}

\newcommand{\rtailt}{{\text{\textrtailt}}}

\usepackage{hyperref}
\hypersetup{
	colorlinks,
	citecolor=black,
	filecolor=black,
	linkcolor=black,
	urlcolor=blue
	bookmarks=true,
}

\usepackage{bookmark}

\journal{Applied and Computational Harmonic Analysis}

\begin{document}

\begin{frontmatter}

\title{Theory of the Hilbert Spectrum}

\author[label1]{Steven~Sandoval\corref{cor1}}
\address[label1]{School of Electrical, Computer and Energy Engineering, Arizona State University, Tempe, AZ\\85287 USA}

\cortext[cor1]{Corresponding author}

\ead{spsandov@asu.edu}
\ead[url]{http://www.HilbertSpectrum.com}

\author[label5]{Phillip~L.~De~Leon}
\address[label5]{Klipsch School of Electrical and Computer Engineering, New Mexico State University, Las Cruces NM 88003 USA}
\ead{pdeleon@nmsu.edu}

\hypertarget{AbstractMark}{}
\bookmark[level=part,dest=AbstractMark]{Abstract}
\begin{abstract}
This paper is a contribution to the old problem of representing a signal in the coordinates of time and frequency. As the starting point, we abandon Gabor's complex extension and re-evaluate fundamental principles of time-frequency analysis. We provide a multicomponent model of a signal that enables rigorous definition of instantaneous frequency on a per-component basis. Within our framework, we have shifted all uncertainty of the latent signal to its quadrature. In this approach, uncertainty is not a fundamental limitation of analysis, but rather a manifestation of the limited view of the observer. With the appropriate assumptions made on the signal model, the instantaneous amplitude and instantaneous frequency can be obtained exactly, hence exact representation of a signal in the coordinates of time and frequency can be achieved. However, uncertainty now arises in obtaining the correct assumptions, i.e.~how to correctly choose the quadrature of the components.

\end{abstract}

\begin{keyword}
Hilbert Space \sep Signal Analysis \sep Instantaneous Frequency \sep Hilbert Spectrum \sep Latent Signal Analysis \sep AM--FM Modeling 
\end{keyword}
\end{frontmatter}

\hypertarget{Part1}{}
\begin{center}\textbf{Part I. Latent Signal Analysis and the Analytic Signal}\end{center}
\bookmark[level=part,dest=Part1]{Part I. Latent Signal Analysis and the Analytic Signal}

\begin{quote}
Interest in the proper definition of instantaneous frequency first arose with the advent of frequency modulation for radio transmission in the 1920's ....  [T]he ``analytic signal procedure'' devised by Gabor results in a complex signal that has spectrum identical to that of the real signal for positive frequencies and zero for the negative frequencies .... However, instantaneous frequency is a primitive concept and not a question of mere mathematical definition .... One should keep an open mind regarding the proper definition of the complex signal, that is, the appropriate way to define phase, amplitude, and instantaneous frequency. Probably the last word on the subject has not yet been said.---Cohen \cite{cohen1995time}
\end{quote}

In the first part of this paper, we present the Latent Signal Analysis (LSA) problem as a recasting of the classic complex extension problem. Almost universally, the solution approach has been to use the Hilbert Transform (HT) to construct Gabor's Analytic Signal (AS). This approach relies on Harmonic Correspondence (HC), which may lead to incorrect Instantaneous Amplitude (IA) and Instantaneous Frequency (IF) parameters. We show that by relaxing HC, the resulting complex extension can still be an analytic function and we can arrive at alternate IA/IF parameterizations which may be more accurate at describing the latent signal. Although the existence of other IA/IF parameterizations is not new, Vakman argued that the AS is the only physically-justifiable complex extension. We argue that by modifying the differential equation for simple harmonic motion \cite{shankar2014fundamentals, FundAcoust}, our parameterizations are also physically justified.

\section{Introduction}
\label{sec:introP1}

Many physical phenomenon are characterized by a complex signal
\begin{equation}
z(t) = x(t)+jy(t)=  \rho(t) e^{j \Theta(t)}
\label{eq:ComplexSignal}
\end{equation} 
where $\rho(t)$ is the signal's IA, $\Theta(t)$ the signal's phase, and $\Omega(t)=\dfrac{d}{dt}\Theta(t)$ the signal's IF. We will call $z(t)$ the \emph{latent signal} because only the real part, $x(t)$ is observed and the imaginary part, $y(t)$ is hidden, i.e.~the act of observation corresponds to the real operator
\begin{equation}
x(t) = \Re\{z(t)\}.
\label{eq:RealObservation}
\end{equation}

\begin{figure}[ht]
	\centering
	  \begin{minipage}[b]{0.49\linewidth}
  		\centering
  		\subfigure[]{
  		\includegraphics[width = 0.95\linewidth]{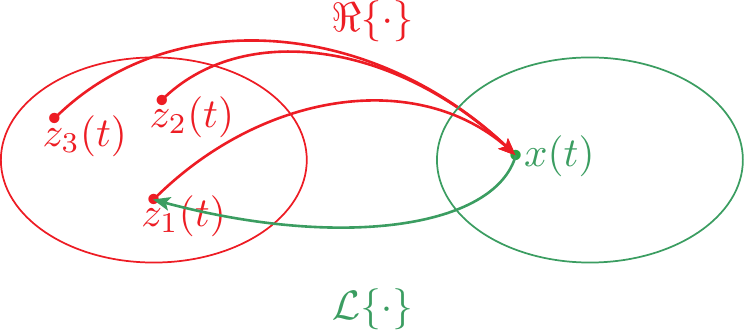}
  		\label{fig:mappinga}
  	}
  	\end{minipage}
  	\begin{minipage}[b]{0.49\linewidth}
  		\centering	
  		\subfigure[]{
  		\includegraphics[width = 0.95\linewidth]{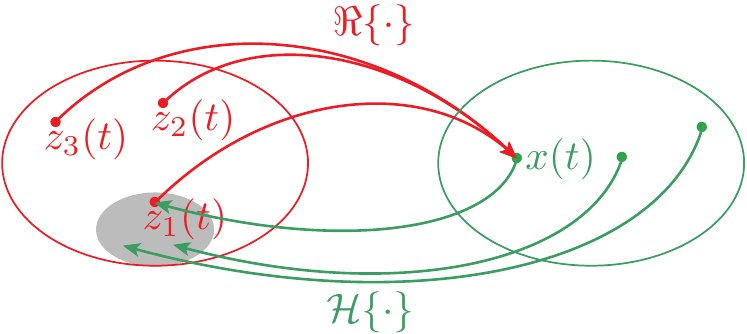}
  		\label{fig:mappingb}
  	}
  	\end{minipage}
    \caption{Set diagrams for the LSA problem. (a) There are an infinite number of latent signals, $z(t)$ that map to the observed signal $x(t)$, according to $x(t) = \Re\{z(t)\}$.  We seek a rule, $\mathcal{L}\{\cdot\}$ to determine an appropriate latent signal $z(t)$. (b) Almost universally the rule chosen is $\mathcal{H}\{\cdot\}$, thus this limits us to only a subset of latent signals as illustrated by the shaded part.}
	\label{fig:mapping}
\end{figure}

It is often desirable to analyze the latent signal---which is known to completely parameterize the physical phenomenon. Thus the complex extension problem becomes that of determining $z(t)$ from the observation $x(t)$, i.e.~recovering the quadrature $y(t)$. In the classical approach, we seek a unique rule $\mathcal{L}\{\cdot\}$ such that the estimate of $y(t)$, $\hat{y}(t) = \mathcal{L}\{x(t)\}$.  These relations and the complex extension problem are illustrated in Fig.~\ref{fig:mappinga}, where we see many latent signals mapping to a single observation under the real operator. In the LSA problem, given $x(t)$ we seek $z(t)$.

In the context of time-frequency signal analysis, we desire the instantaneous parameters, $\rho(t)$ and $\Omega(t)$.
Once we have determined a rule for $\hat{y}(t)$, the instantaneous estimates are given by
\begin{equation}
	\hat{\rho}(t)    = \pm\left|\hat{z}(t)\right|= \pm\sqrt{x^2(t)+\hat{y}^2(t)}
	\label{eq:instAmp}
\end{equation}
and
\begin{equation}
	\hat{\Omega}(t) = \frac{d}{dt}\left[\arctan\left(\frac{\hat{y}(t)}{x(t)}\right)\right].
	\label{eq:instFreq}
\end{equation}
The very definition of IA and IF hinges on $\hat{y}(t)$ and hence $\mathcal{L}\{\cdot\}$.

The complex extension problem is well-known. In 1937, Carson and Fry formally defined the IF based on the phase derivative of a complex FM signal \cite{CarsonFry,VanDerPol1946}. In Gabor's seminal 1946 paper \cite{gaborTOC}, he introduced a Quadrature Method (QM) as a practical approach for obtaining the complex extension of a real signal and showed its equivalence to the HT. In this context, the rule is given by
\begin{equation}
\hat{y}(t) = \mathcal{H}\{x(t)\}
\label{eq:HTRule}
\end{equation}
where $\mathcal{H}\{\cdot\}$ is the HT operator and 
\begin{equation}
\hat{z}(t) = x(t) +j \mathcal{H}\{x(t)\}
\label{eq:ASdef}
\end{equation}
is termed the AS. This is illustrated in Fig.~\ref{fig:mappingb} where we see using the HT to estimate the quadrature leads us to a subset of the latent signals.

However shortly afterward, Shekel pointed out the ambiguity problem in Ville's work defining the instantaneous parameters of a real signal \cite{shekel1953instantaneous,picinbono1997instantaneous,BoashashP1}.  As an example, consider
\begin{equation}
    x(t) = \Re\left\{ a_0(t)e^{j[ \int\limits_{-\infty}^{t} \omega_0(\uptau)d\uptau+\phi_0 ]}\right\}.
\end{equation}
There are an infinite set of pairs $a_0(t)$ and $\omega_0(t)$ for which $x(t)$ may be equivalently described and hence an infinite set of IA/IF parameterizations. Shekel pointed this ambiguity out ``with the hope of banishing it [IF] forever from the dictionary of the communication engineer.''  Others such as Hupert, suggested that despite the ambiguity, the concept of instantaneous parametrization of real signals was still useful and could be applicable \cite{Hupert1953instantaneous,rihaczek1966hilbert}.

In order to constrain the ambiguity problem so that a unique complex extension can be justified for a real signal, Vakman proposed three conditions tied to physical reality \cite{vakman1972, vakman1976, vakman1979}: 1) amplitude continuity, 2) phase independence of scaling and homogeneity, and 3) HC. With these conditions, Vakman showed that the unique complex extension is given by the rule $\hat{y}(t)=\mathcal{H}\{x(t)\}$.

More recently, other authors have proposed other conditions to constrain the ambiguity, such as bounded amplitude and bounded IF variation, leading to alternate IA/IF solutions \cite{loughlin1996amplitude,loughlin1998bounded}. Finally, we point out that a number of other methods have also been introduced, without explicitly stating the conditions. Vakman has shown that all these methods violate one or more of the conditions that he proposed and almost all retain HC \cite{vakman1996analytic,VakmanBook}. Like many authors, we believe Vakman's conditions justifying Gabor's use of the HT are reasonable and sensible for several special cases of $z(t)$. On the contrary, we believe that Vakman's conditions, in particular HC, can be too restrictive for more general $z(t)$ and in the worst case, may lead to incorrect results and interpretations.

Part I is organized as follows. In Section \ref{sec:theHT}, we review the HT and its motivations and the AS. In Section \ref{sec:relaxing}, we give a proof showing we can still maintain analyticity without HC and also give a proof that no universal rule for the quadrature exists. In Section \ref{sec:examplesP1}, we provide an example of a LSA problem and several solutions. Finally, in Section \ref{sec:sumP1} we summarize Part I.

\section{The Hilbert Transform and Analytic Signal}
\label{sec:theHT}
Although there exist several methods for estimating the instantaneous parameters, use of the HT dominates science and engineering. The HT of $x(t)$ is given by
\begin{equation}
	\label{eq:HT}
   \mathcal{H}\lbrace x(t)\rbrace \equiv -\frac{1}{\pi} \dashint\limits_{-\infty}^{\infty}\frac{x(\uptau)}{\uptau-t}d\uptau
\end{equation}
where $\dashintS$ indicates the Cauchy principle value integral \cite{titchmarsh, LeeWienerLegacy}. The three main motivations for use of the HT are: 1) Vakman's physical conditions, 2) analyticity of the resulting complex extension, and 3) ease of computation via Gabor's QM. In this section, we review the motivations.  

\subsection{Vakman's Physical Conditions}
Vakman proposed conditions in order to constrain the ambiguity in choosing the complex extension \cite{vakman1972,vakman1976,vakman1979,vakmanVainshtein1977,VakmanBook}.  

\emph{Condition 1 Amplitude Continuity:}
Simply stated, amplitude continuity requires that the IA, $\rho(t)$ is a continuous function.
This implies that the rule, $\hat{y}(t) = \mathcal{L}\{x(t)\}$ must be continuous, i.e.
\begin{equation}
\mathcal{L}\{x(t)+\epsilon w(t)\} \rightarrow \mathcal{L}\{x(t)\}~\mathrm{for}~||\epsilon w(t)||\rightarrow 0.
\end{equation}

\emph{Condition 2 Phase Independence of Scaling and Homogeneity:}
Let $x(t)$ have a complex extension, $z(t) = \rho(t) \exp [j\Theta(t)]$. Then for a real constant $c>0$, $c x(t)$ has associated complex extension $z_1(t) = [c \rho(t)] \exp[j\Theta(t)]$, i.e.~only the IA of the complex representation is affected and $\Theta(t)$ and $\Omega(t)$ remain unchanged. This implies that the rule for performing the complex extension is scalable
\begin{equation}
\mathcal{L}\{c x(t)\} = c \mathcal{L}\{x(t)\}.
\end{equation}
Conditions 1 and 2 force the operator to be linear \cite{VakmanBook}.

\emph{Condition 3 Harmonic Correspondence:}
Let $x(t)=a_0\cos(\omega_0 t + \phi_0)$, then HC forces the complex extension, 
\begin{equation}
\hat{z}(t) = a_0e^{j(\omega_0 t+\phi_0)},
\label{eq:shc}
\end{equation}
where we note the IA and IF must be constant. This implies that 
\begin{equation}
\mathcal{L}\{a_0\cos(\omega_0 t+\phi_0)\} = a_0\sin(\omega_0 t+\phi_0)
\label{eq:HarmonicCoor}
\end{equation}
and $\hat{z}(t)$ is a Simple Harmonic Component (SHC) with positive IF.

As Vakman showed \cite{vakman1972}, the HT is the only operator that satisfies the above conditions and as a result, the HT (or Gabor's practical QM implementation) is viewed as the correct way to complex extend a signal.  With his work, Vakman was able to refute most objections as to whether use of the HT is physically justified.

\emph{Condition 4 Phase Continuity:}
In addition to real signals having an ambiguity in instantaneous parameterization, complex signals also have an ambiguity: although we construct $\hat{z}(t)$ using a rule for $y(t)$, we want the IA/IF pair $\rho(t),~\Omega(t)$ for signal analysis, and there can be ambiguity in this coordinate transformation \cite{cohenLoughlinVakman1999}. There exist at least two choices for resolving this ambiguity when it arises:
\begin{itemize}
	\item \textbf{Condition 4a} Positive IA: $\rho(t)\geq0~\forall~t$.
	\item \textbf{Condition 4b} Phase continuity: the phase $\Theta(t)$ is a continuous function.
\end{itemize}
Although condition 4a is the traditional choice, we advocate condition 4b to ensure the IF is well defined. This is apparent in (\ref{eq:instAmp}) where the IA may be negative.

\subsection{Gabor's Quadrature Method}
The HT is an ideal operator that in practice is not physically realizable.  Gabor's QM provides a frequency domain approach that under certain conditions is equivalent to the HT. The steps in Gabor's QM are:
\begin{enumerate}
	\item decompose $x(t)$ into SHCs, i.e.~compute the Fourier spectrum
	\item double the magnitude of the non-negative frequency components
	\item negate the negative frequency components.
\end{enumerate}
Gabor's QM results in a complex signal formulated in terms of non-negative spectral frequencies,
\begin{equation}
\hat{z}(t)=  \sum\limits_{k=0}^{K-1}  a_k \exp\left\lbrace j \left[\omega_kt +\phi_k\right] \right\rbrace
\label{eq:ComplexSignalwSHCs}
\end{equation}
i.e.~we have decomposed the signal into SHCs each having constant IA, $a_k$ and (non-negative) constant IF, $\omega_k$. By comparing the expressions for $\hat{z}(t)$ in (\ref{eq:ComplexSignalwSHCs}) and (\ref{eq:ComplexSignal}), we can effectively collapse $K$ SHCs into a single AM-FM component and then obtain an IA/IF pair. Implementing Gabor's QM using a FT is very convenient and is one reason for the method's popularity.

\subsection{Analyticity of the Analytic Signal}
Unfortunately, the word ``analytic'' has two distinct meanings when used in signal processing and in addition, another meaning in mathematics. A signal is said to be \emph{analytic} if it consists only of non-negative frequency components \cite{cohen1995time,boashash2003time}. \emph{The analytic signal} refers to the complex extension of a real signal using the HT, i.e.~Gabor's method \cite{ville, bedrosian1962analytic, cohen1995time, xia1999analytic, boashash2003time}. If $z(\rtailt)$, where $\rtailt\equiv t + j \tau$, is an \emph{analytic function} then the real and imaginary parts of the complex function,
\begin{equation}
z(\rtailt)=u(t,\tau)+j v(t,\tau)
\label{eq:ComplexZ}
\end{equation}
satisfy the Cauchy-Riemann (CR) conditions
\begin{equation}
\dfrac{\partial}{\partial t}u\left(t,\tau \right) = \dfrac{\partial }{\partial \tau}v\left(t,\tau \right)~~~\text{and}~~~\dfrac{\partial}{\partial \tau}u\left(t,\tau\right) = -\dfrac{\partial}{\partial t}v\left(t,\tau \right).
\label{eq:CR}
\end{equation}
For the AS $\hat{z}(t) = x(t)+j\mathcal{H}\{x(t)\}$, such that $t\gets \rtailt\equiv t + j \tau$, the complex function $\hat{z}(\rtailt)=\hat{u}(t,\tau)+j \hat{v}(t,\tau)$ is an analytic function \cite{brown2009complex}. Hence the reason for calling a HT-extended signal the AS \cite{ville, cohen1995time}.

\section{Relaxing the Condition of Harmonic Correspondence}
\label{sec:relaxing}

\subsection{On Harmonic Correspondence}

We can understand Vakman's motivation for tying HC to physical reality by considering the differential equation
\begin{equation}
\dfrac{d^2}{dt^2}z(t) + \omega_0^2 z(t)=0
\label{eq:simpleharmonicmotion}
\end{equation}
which describes many ideal systems, e.g.~a LC circuit or mass/spring model. The solution to (\ref{eq:simpleharmonicmotion}) is the SHC in (\ref{eq:shc}).

Any deviation from this ideal case, requires modification of the differential equation.  For example, when the differential equation describes a circuit and resistance is included or describes motion and damping is included, the equation becomes
\begin{equation}
\dfrac{d^2}{dt^2}z(t) + c \dfrac{d}{dt}z(t)+\omega_0^2 z(t)=0
\end{equation}
where $c$ is a constant. In this case, the solution includes an AM term and has the form
\begin{equation}
z(t) = a_0 e^{-\nu t} e^{j(\omega_d t+\phi_0)}
\end{equation}
which is not a SHC \cite{FundAcoust}. As another example, when the differential equation is further modified to include time-varying coefficients, non-linearities, or partial derivatives with respect to $\tau$ or spacial variables, the resulting solution may include both AM and FM terms, which is also not a SHC \cite{ van1934nonlinear, VakmanBook, farlow2012partial}.  

While most authors believe HC is a reasonable condition to describe physical phenomena, we advocate that this condition is overly constraining. For many analysis problems, this can lead to incorrect interpretations because real physical systems always deviate from the ideal case. By not assuming HC, we gain a degree of freedom in our analysis that allows us to construct other complex extensions that may be better suited to describing the underlying physical phenomena associated with the signal. We believe that any attempt to find a unique rule to infer $z(t)$ of the form in (\ref{eq:ComplexSignal}) from $x(t)$ in the form of (\ref{eq:RealObservation}) \emph{is fundamentally flawed and that no such universal rule can exist}.

\subsection{On Analyticity of the Complex Extension}

We believe that the use of the term ``analytic'' has in most cases become too restrictive in signal processing. To wit, one may falsely believe that the HT is the only way to complex-extend a \emph{real signal} in order to result in an \emph{analytic function}, when time is considered complex. Other complex extensions of real signals can be constructed that result in analytic functions. Choosing $\hat{y}(t)=\mathcal{H}\{x(t)\}$ to obtain the AS $\hat{z}(t)$ ensures $z(\rtailt)$ is an analytic function. However, there are other choices for $\hat{y}(t)\neq\mathcal{H}\{x(t)\}$ such that $z(\rtailt)$ is an analytic function. 

\textbf{Theorem} If we do not assume HC, there exists at least one choice for the quadrature, $y(t)\neq\mathcal{H}\{x(t)\}$ that results in $z(\rtailt)$ being an analytic function.

\begin{proof}Let $x(t)= a_0\cos(\omega_0 t)$ and choose the quadrature such that
	\begin{equation}
	z(t) = a_0 \cos(\omega_0 t) + j \alpha a_0 \sin( \beta \omega_0 t)
	\label{eq:NotHarmonicCoor}
	\end{equation}
	with real $\alpha$ and $\beta$. The complex function is given by
	\begin{equation}
	z(\rtailt) = a_0 \cos(\omega_0 [t+j\tau]) + j \alpha a_0 \sin(\beta\omega_0 [t+j\tau])
	\label{eq:zrtailt}
	\end{equation}
	where
	\begin{equation}
	u(t,\tau)= a_0 \cos(\omega_0 t)\cosh(\omega_0 \tau)- \alpha a_0 \cos(\beta \omega_0 t)\sinh(\beta \omega_0 \tau)
	\label{eq:uttau}
	\end{equation}
	and
	\begin{equation}
	v(t,\tau)= \alpha a_0 \sin(\beta \omega_0 t)\cosh(\beta\omega_0 \tau)- a_0 \sin(\omega_0 t)\sinh(\omega_0 \tau) .
	\label{eq:vttau}
	\end{equation}
	It can easily be shown that this choice leads to $z(\rtailt)$ satisfying the CR conditions and hence is an analytic function. Any choice of $\alpha\neq1$ or $\beta\neq 1$ does not imply HC.
\end{proof}

Although Gabor's QM provides a simple rule to obtain an IA/IF pair that parameterizes $x(t)$, it may not address the problem of obtaining the latent signal from the observation. The heart of the problem is that no single rule can determine the latent signal from the observation for every possible latent signal because more than one complex signal map to the same real signal under the real operator. Although a unique rule can be constructed to work for a particular $z(t)$, the same rule cannot generally be used. 

\textbf{Corollary} No unique rule for the quadrature, $\hat{y}(t)=\mathcal{L}\{x(t)\}$ exists to obtain the latent signal, $z(t)$ from the observation, $x(t)=\Re\{z(t)\}$ for all $z(t)$.

\begin{proof}
Consider the latent signal, $z(t)$ of the form in (\ref{eq:NotHarmonicCoor}). Assume a unique rule, $\hat{y}(t)=\mathcal{L}\{x(t)\}\equiv\alpha_0 a_0 \sin(\beta_0\omega_0 t)$ exists to obtain $z(t)$. This implies a unique $\hat{z}(t)= a_0 \cos(\omega_0 t) + j \alpha_0 a_0 \sin(\beta_0\omega_0 t)$. If $\alpha_0\neq\alpha$ or $\beta_0\neq\beta$ then $\hat{z}(t)\neq z(t)$ and $\mathcal{L}\{\cdot\}$ is not unique.
\end{proof}

\subsection{On Harmonic Conjugate Functions}

If $z(\rtailt)=u(t,\tau)+j v(t,\tau)$ is an analytic function, then $u(t, \tau)$ and $v(t, \tau)$ are unique and are \emph{harmonic conjugates} \cite{brown2009complex}. Even though $x(t)=u(t,0)$ and $y(t)=v(t,0)$ in (1) this does not imply that $u(t,\tau)= x(\rtailt)$ and $v(t,\tau)= y(\rtailt)$, but rather $u(t,\tau)$ and $v(t,\tau)$ each contain terms from both $x(\rtailt)$ and $y(\rtailt)$:
\begin{eqnarray}
z(\rtailt) &=& x(\rtailt)+j y(\rtailt)\nonumber\\
           &=& \left[x_R(t,\tau)+jx_I(t,\tau)\right]+j\left[y_R(t,\tau)+jy_I(t,\tau)\right]\nonumber\\
           &=& \left[x_R(t,\tau)-y_I(t,\tau)\right]+j\left[x_I(t,\tau)+y_R(t,\tau)\right]\nonumber\\
           &=& u(t,\tau) +j v(t, \tau)
\end{eqnarray}
where $R,~I$ denote the real and imaginary parts, respectively. This is easily seen in (\ref{eq:uttau}), where $\alpha$ and $\beta$ are present in $u(t,\tau)$ despite $\alpha$ and $\beta$ appearing only in $y(t)$ in (\ref{eq:NotHarmonicCoor}). Thus the problem of finding $y(t)$ cannot be solved by finding a harmonic conjugate of $u(t,\tau)$ because $y(t)$ must be known to obtain $u(t,\tau)$.

\subsection{On AM--FM Demodulation}

The HT is widely used as a demodulation algorithm for AM--FM signals. This is typically justified as a valid approach because of Vakman and also Bedrosian's theorem.\footnote{If the product $l(t)h(t)$ consists of a low-frequency factor, $l(t)$ and a high-frequency factor, $h(t)$ of non-overlapping spectra, then $\mathcal{H}\{l(t)h(t)\} = l(t)\mathcal{H}\{h(t)\}$ \cite{Bedrosian1963,Nuttall1966,VakmanBook}.} The HT can be used to determine the IA/IF for a small subset of latent signals, i.e. those with HC. The HT can also be used to closely approximate the IA/IF for latent signals whenever $\mathcal{H}\{x(t)\}\approx y(t)$ \cite{picinbono1997instantaneous}. However, the HT cannot be used in general to obtain the IA/IF of \emph{all} latent signals. 

\section{Example}
\label{sec:examplesP1}

Recall the LSA problem: given the observation $x(t)=\Re\{z(t)\}$, find $z(t)$ or more strictly the instantaneous parameters, $\rho(t)$ and $\Omega(t)$. In this section, we will give three solutions to a LSA problem where two of these solutions have $\hat{y}(t)\neq\mathcal{H}\{x(t)\}$, as illustrated in Fig.~\ref{fig:mappingb}. 

Suppose we have three ideal systems:  a frequency modulator, an amplitude modulator, and a Linear Time-Invariant (LTI) system in the steady state. We observe a triangle waveform $x(t)$ that could have come from any of the three ideal systems. What is the corresponding latent signal $z(t)$ or more strictly the instantaneous parameters, $\rho(t)$ and $\Omega(t)$?

\begin{figure}[ht]
	\centering
  	\begin{minipage}[b]{0.5\linewidth}
  		\centering
  		\includegraphics[width = \linewidth]{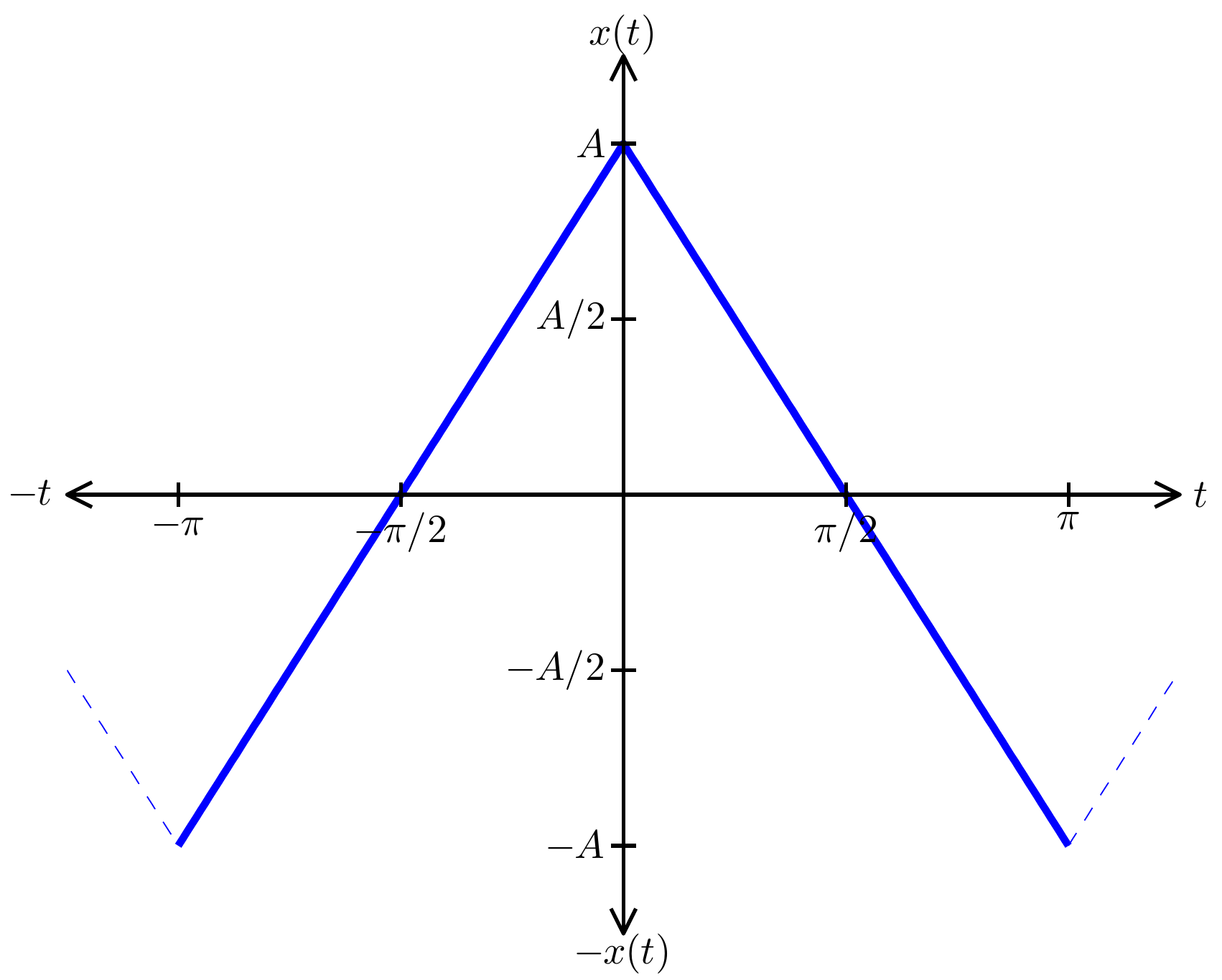}
  		\label{fig:FStriangle}
  	\end{minipage}	
  	\caption{One period of the triangle waveform, $x(t)$ in (\ref{eq:trianglewaveP1}) with amplitude $A$ and $\omega_0=1$.}
  	\label{fig:TRIwave}
\end{figure}

Let the observation be the periodic, even-symmetric triangle waveform where one period is given by
\begin{equation}
x(t)= \left\{      
\begin{array}{ll}
-2A\omega_0 t/\pi+A,& 0\leq t \leq T/2\\
 2A\omega_0 t/\pi+A,&-T/2\leq t \leq 0
\end{array}
\right.
\label{eq:trianglewaveP1}
\end{equation}
with amplitude $A$, period $T$, and fundamental frequency $\omega_0$, as illustrated in Fig.~\ref{fig:TRIwave}.

\subsection{Solution Assuming Harmonic Correspondence}
If we assume HC, i.e.~the ideal system is LTI in the steady state then
\begin{subequations}
	\begin{align}
	\label{eq:trianglesingleHCP1}
	\hat{z}_1(t)   &=  x(t)  +j \mathcal{H}\lbrace x(t) \rbrace  \\[-0.1em]
	&= \sum\limits_{k=0}^{\infty} \frac{8A}{\pi^2(2k+1)^2} \cos\left[ (2k+1) \omega_0 t\right]  +j \sum\limits_{k=0}^{\infty} \frac{8A}{\pi^2(2k+1)^2} \sin\left[ (2k+1) \omega_0 t\right]    
	\end{align}
\end{subequations}
where the IA is
\begin{equation}
\hat{\rho}(t) = \frac{8A}{\pi^2} \tilde{a}_0(t)
\label{eq:trianglesingleamfmiaP1}
\end{equation}
and the IF is
\begin{equation}
\hat{\Omega}(t) = \omega_0+\frac{d}{dt}M_0(t)
\label{eq:trianglesingleamfmifP1}
\end{equation}
where $\tilde{a}_0(t)$ and $M_0(t)$ are related through
\begin{equation}
\tilde{a}_0(t)e^{jM_0(t)} = \sum\limits_{k=0}^{\infty} \frac{1}{(2k+1)^2}e^{j 2k\omega_0 t}.
\end{equation}

\subsection{Solution Assuming Constant IF}
If we assume constant IF $\hat{\Omega}(t)=\omega_0$, i.e.~the ideal system is an amplitude modulator then
\begin{equation}
\hat{z}_2(t)=   \hat{\rho}(t) \cos(\omega_0 t) +j \hat{\rho}(t)\sin(\omega_0 t) 
\label{eq:trianglesingleamP1}
\end{equation}
with the IA given by
\begin{equation}
\hat{\rho}(t) =  x(t)/\cos(\omega_0 t).
\label{eq:trianglesingleamiaP1}
\end{equation}
In this solution, $\hat{y}(t) = \hat{\rho}(t)\sin(\omega_0 t)\neq\mathcal{H}\{\hat{\rho}(t) \cos(\omega_0 t)\}$, because Bedrosian's theorem cannot be applied. 

\subsection{Solution Assuming Constant IA}
If we assume constant IA $\hat{\rho}(t)=A$, i.e.~the ideal system is a frequency modulator then
\begin{equation}
    \hat{z}_3(t) = A \cos \left[\omega_0 t+M_0(t)\right]+j A \sin \left[\omega_0 t+M_0(t)\right] 
    \label{eq:trianglesinglefmP1}
\end{equation}
with IF given by
\begin{equation}
    \hat{\Omega}(t) = \omega_0+\frac{d}{dt}M_0(t)
\label{eq:trianglesinglefmifP1}
\end{equation}\vspace{-2mm}
where
\begin{eqnarray}
    M_0(t)                    & = &  \arccos\left[ x(t)/A \right]-\omega_0 t.
\end{eqnarray}
As in the solution for constant IF, $\hat{y}(t) = A \sin \left[\omega_0 t+M_0(t)\right]  \neq \mathcal{H}\{A \cos \left[\omega_0 t+M_0(t)\right]\}$.

\section{Summary}
\label{sec:sumP1}
We have presented the LSA problem and reviewed the motivation for the use of the HT and AS in signal analysis. We proved the HC condition is not necessary for analyticity and that no unique rule for the complex extension exists to obtain the latent signal. It was argued that by relaxing the HC condition, which forces the use of the HT, we gain alternate choices for the latent signal. Although the HT is widely used for demodulation of AM--FM signals,  it cannot be used to obtain the IA/IF of all latent signals. In a strict sense, Gabor's QM can only be used when the latent signal is a superposition of SHCs with non-negative IF and can approximate the latent signal whenever $\mathcal{H}\{x(t)\}\approx y(t)$, e.g.~communication signals where Bedrosian's theorem is approximately satisfied.

\clearpage
\hypertarget{Part2}{}
\begin{center}\textbf{Part II. Hilbert Spectral Analysis Using Latent AM--FM Components}\end{center}
\setcounter{section}{0}
\bookmark[level=part,dest=Part2]{Part II. Hilbert Spectral Analysis Using Latent AM--FM Components}

\begin{quote}
[W]ell-known methods ... can no longer be applied when the frequency itself is made a function of the time. Even the concept of ``instantaneous-frequency''... is of a somewhat arbitrary but nevertheless highly useful nature, is often misunderstood and even misinterpreted .... The only way at present available to solve these and similar problems is to go back to the very first and fundamental principles. This implies a theoretical treatment beginning with the differential equations of the problem concerned. Unfortunately these equations can seldom be solved in terms of well-known functions, such as real or complex exponentials [simple harmonic components]. I think it is precisely to this fact that the main difficulties can be traced.---Van der Pol \cite{VanDerPol1946}
\end{quote}

In the second part of the paper, the central theme is Hilbert Spectral Analysis (HSA) using a superposition of latent AM--FM components. We present HSA as a generalized LSA problem.  In the general problem, we seek a representation of $z(t)$ consisting of a superposition of latent components, i.e.~a multicomponent model.  Furthermore, rather than seeking a single IA/IF pair for $z(t)$, we seek a set of IA/IF pairs each associated with the components. Although time-frequency analysis has been extensively studied the use of a generalized AM--FM model for this analysis, without the HC condition, has never been proposed. Using this model leads to non-unique signal decompositions due to the theorem proved in Part I. However, by imposing assumptions on the form of the AM--FM component, a unique parameterization in terms of IA and IF can be obtained. This analysis requires abandoning Gabor's complex extension as advocated in Part I and instead allowing the assumptions to imply the complex extension. This model enables us to analyze signals with very few restrictions resulting in alternate and possibly more useful decompositions, especially for nonstationary signals. In this part, AM--FM modeling and HSA theory is presented. Examples using HSA are given and a visualization of the Hilbert spectrum is proposed.

\section{Introduction}
\label{sec:introP2}
It has long been known that simple harmonic analysis of nonstationary signals may lead to incorrect interpretations of an underlying signal model \cite{BoashashP1}. In Gabor's seminal paper, ``Theory of Communications'' he writes: 
\begin{quote}
    The greatest part of the theory of communication has been built up on the basis of Fourier's reciprocal integral relations .... Though mathematically this theorem is beyond reproach ... even experts could not at times conceal an uneasy feeling when it came to the physical interpretation of results obtained by the Fourier method. After having for the first time obtained the spectrum of a frequency-modulated sine wave, Carson wrote: `The foregoing solutions, though unquestionably mathematically correct, are somewhat difficult to reconcile with our physical intuitions ....' \cite{gaborTOC}
\end{quote}

As Gabor noted, Carson in 1922 was the first to rigorously study a FM signal and realize that the frequency components obtained from such a nonstationary signal, do not describe the physical system in a meaningful way \cite{Carson1922Notes}. A similar argument regarding an AM signal was made by Priestly, who wrote,  ``...a nonstationary process in general cannot be represented in a meaningful way by the simple Fourier expansion'' \cite{priestley1988non}.  His example,
\begin{equation}
    x(t) = Ae^{-t^2/o^2} \cos(\omega_0 t + \phi_0 )
\end{equation}
and its Fourier Transform (FT) consists of two Gaussian functions centered at frequencies $\pm\omega_0$. Thus, the FT contains an infinite number of SHCs. This interpretation of the underlying signal model may be incorrect. For example, an alternative signal model consists of just two components at constant frequencies $\pm\omega_0$, with each component having a \emph{time-varying amplitude}, $(A/2)e^{-t^2/o^2}$. Mathematically, these two representations are equally valid and correspond to different families of basic functions used for representation \cite{BoashashP1}. 

For nonstationary signals, SHCs lose effectiveness and thus the idea of time-varying components\footnote{Some authors, such as Ville \cite{ville}, refer to time-varying components as ``instantaneous spectra'' or more generally as functions of time that give the structure of a signal at a given instant.}, i.e.~components with time-varying IA and IF has arisen in order to account for the non-stationarity \cite{priestley1988non,BoashashP1}. However, the concept of IF is not without its own controversy primarily due to four reasons:
\begin{enumerate}
    \item There is an apparent paradox in associating the words ``instantaneous'' and ``frequency'' because frequency \emph{usually} defines the number of cycles undergone during one unit of time \cite{BoashashP1}.

    \item Without assumptions, instantaneous parameterizations of a signal are not unique \cite{shekel1953instantaneous, BoashashP1, cohen1995time, loughlin1998bounded}.

    \item The commonly accepted definition of IF, i.e.~phase derivative of Gabor's AS is correct for only a limited class of signals \cite{cohen1995time}.

    \item Although different quantities, harmonic frequency and IF are often confused likely due to the term ``frequency'' attached to both \cite{fink1966relations, mandel1974interpretation}. Harmonic frequency is a special case of IF and is \emph{only} equivalent under the assumption of SHCs.
\end{enumerate}
Several authors over the previous decades have shown that problems and paradoxes exist that are related to the definition of IF \cite{shekel1953instantaneous, fink1966relations,loughlin1998bounded, cohen1995time,mandel1974interpretation,gupta1975definition, vakmanVainshtein1977, priestley1988non, loughlin1997comments,oliveira1998concept, oliveira1999instantaneous, nho1999instantaneous, chui2015signal}. However on the whole, these problems seem to have been forgotten or ignored. In 1937, Carson and Fry formally defined the IF based on the phase derivative of a complex FM signal assuming the FM message is slowly-varying \cite{CarsonFry}. This assumption, while perfectly valid in communication theory, implies a narrowband component when used in signal analysis.

In order to exploit the mathematical convenience of using complex exponentials in signal analysis, Gabor introduced a QM for obtaining a complex extension of a real signal \cite{gaborTOC}. Gabor's QM assumed positive IF SHCs and was shown to be equivalent to the HT \cite{gaborTOC}. Although Gabor's QM provides a unique method for complex extension, as Cohen pointed out without strict adherence to the assumption, this results in many counter-intuitive consequences \cite{cohen1995time}. Ville defined the IF of a real signal by using Gabor's complex extension and then Carson's definition of IF \cite{BoashashP1}. By defining the IF as the derivative of the phase of Gabor's AS, Ville was able to show the average harmonic frequency is equal to the time average of the IF \cite{nho1999instantaneous}. He then formulated the Wigner-Ville Distribution (WVD) and showed that the first moment of the WVD with respect to frequency yields the IF. Using Gabor's AS to extend a real signal results in a number of very convenient relations. As a result, Gabor's QM is almost universally viewed as the correct way to define the complex signal and subsequently the correct way to define IA, IF, and phase for real signals despite the inherent assumption of HC \cite{cohen1995time}. 

In our research, we have found that the dogmatic use of Gabor's complex extension does not provide the necessary flexibility for modeling nonstationary signals. For many problems, as we will demonstrate, it can be advantageous to allow the assumptions of an underlying signal model to \emph{imply} the complex extension. As a result, assumptions must be made on a per problem basis in order to force a unique decomposition and hence, in order to properly estimate the latent signal and its components. This leads to a decomposition into complex AM--FM components rather than SHCs. As a result of using AM--FM components, we revert back to Carson's original definition of IF and by further relaxing HC and embracing the resulting non-uniqueness, the \emph{associated problems and paradoxes related to IF can all be resolved}. 

Part II of this paper is organized as follows.  In Section~\ref{sec:AMFMmodel}, we propose the AM--FM model. In Section~\ref{sec:HSA}, we propose a general Hilbert spectrum based on the AM--FM model and examine several familiar specializations. In Section~\ref{sec:FDVoLSA}, a frequency domain view of LSA is presented. In Section~\ref{sec:Subtleties}, we discuss subtleties of complex extension of real signals assuming AM-FM components with relaxed HC. In Section~\ref{sec:Solutions}, we provide examples of HSA and in Section \ref{sec:visual} propose a visualization of the Hilbert spectrum to aid in the interpretation. In Section~\ref{sec:SummaryP2}, we summarize Part II.

\section{The AM--FM Model}\label{sec:AMFMmodel}

\subsection{Definitions and Assumptions}
\begin{figure}[h]
	\centering
	\includegraphics[width = 0.75\linewidth]{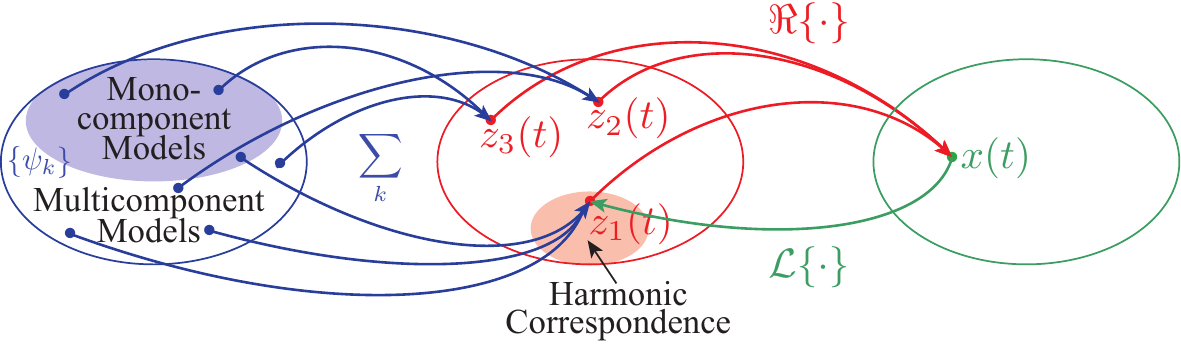}
	\caption{A set diagram for the HSA problem. In the LSA problem, many latent signals $z(t)$ map to the same observed signal $x(t)$ under the real operator.  In the HSA problem, many component sets $\{\psi_k(t)\},~k = 0,1,\ldots,K$ map to the same latent signal $z(t)$ under superposition.  The goal in HSA is to properly choose $\{\psi_k(t)\}$ given $x(t)$.}
	\label{fig:NewSet}
\end{figure}

Our goal in HSA is to decompose a signal into complex AM--FM components, each of which yield IA and IF estimates that are matched to some criteria. That criteria could be some physical or perceptual observation or simply some intuitive notion of an underlying signal model. The decomposition problem is illustrated in Fig.~\ref{fig:NewSet} where we have extended the LSA problem shown in Fig.~\ref{fig:mapping} to show  models of the latent signal that in general, are composed of a set of latent AM-FM components. 

We propose to define the AM--FM model for $z(t)$ as a superposition of $K$ (possibly infinite) AM--FM components
\begin{equation}
	z(t) \equiv  \sum\limits_{k=0}^{K-1}\psi_k\left(t;a_k(t),{\omega}_k(t),\phi_k\right) 
	\label{eq:AMFMmodel}
\end{equation}
where the AM--FM component is defined as 
\begin{subequations}
	\label{eq:AMFMcomp}
	\begin{align}
	  	\label{eq:AMFMcompA}
	    \psi_k(t;a_k(t),{\omega}_k(t),\phi_k) &\equiv a_k(t) \exp\left\lbrace j \left[     \int\limits_{-\infty}^{t}\omega_k(\uptau)d\uptau +\phi_k\right] \right\rbrace\\
	  	\label{eq:AMFMcompB}
	  	&= a_k(t) e^{j \theta_k(t)}\\
	  	\label{eq:AMFMcompC}
	  	&= s_k(t)+j\sigma_k(t)
	\end{align}
\end{subequations}
parameterized by the IA $a_k(t)$, IF $\omega_k(t)$, and phase reference $\phi_k$.  We assume the observed real signal $x(t)$ is related to the latent signal $z(t)$ by (\ref{eq:RealObservation}).

\begin{figure}[ht]
	\centering
	  	\begin{minipage}[b]{0.5\linewidth}
  		\centering
  		\subfigure[]{
  		\includegraphics[width = 0.95\linewidth]{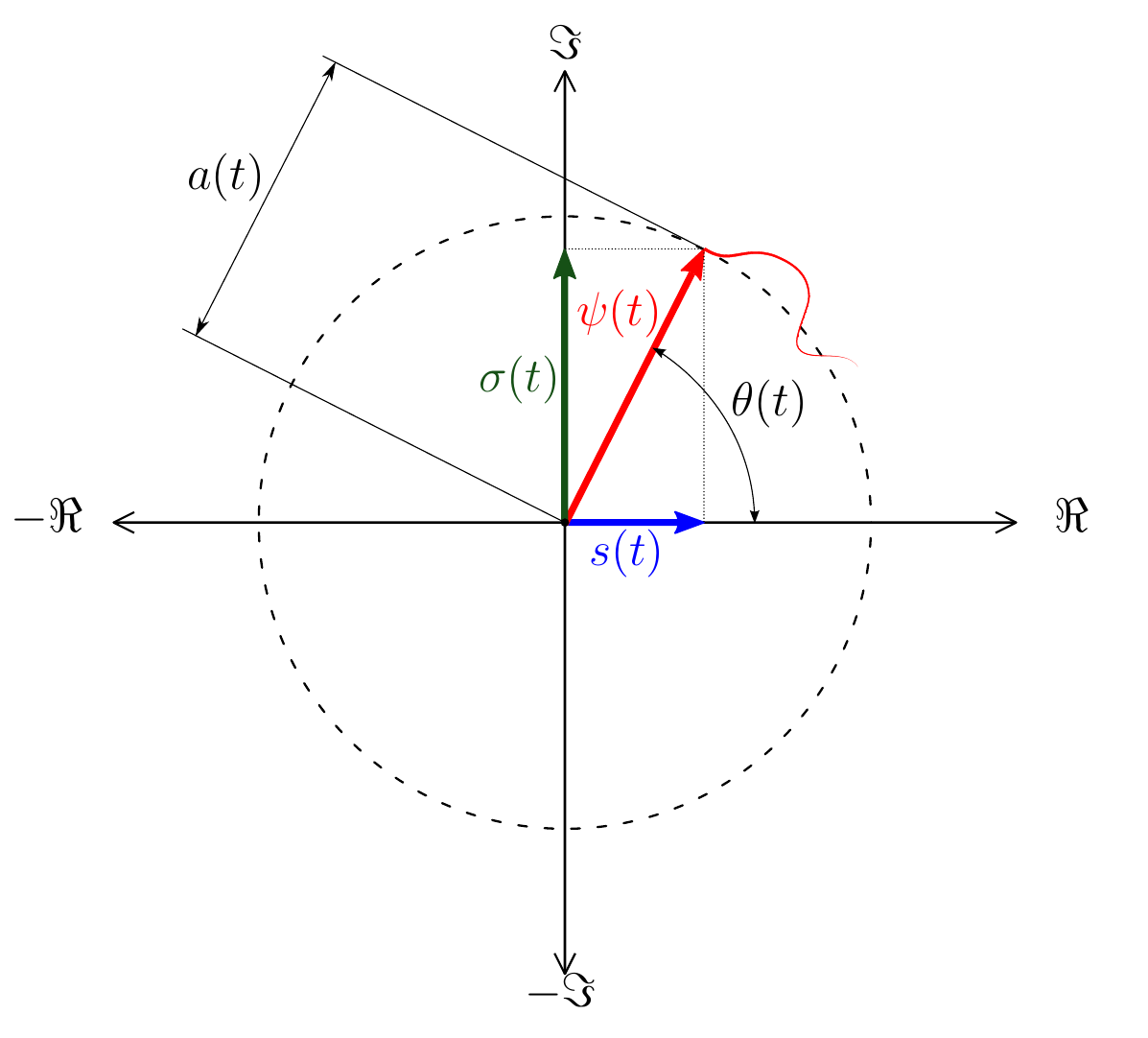}
  		\label{fig:aa}
  	}
  	\end{minipage}\hfill
  	\begin{minipage}[b]{0.5\linewidth}
  		\centering	
  		\subfigure[]{
  		\includegraphics[width = 0.95\linewidth]{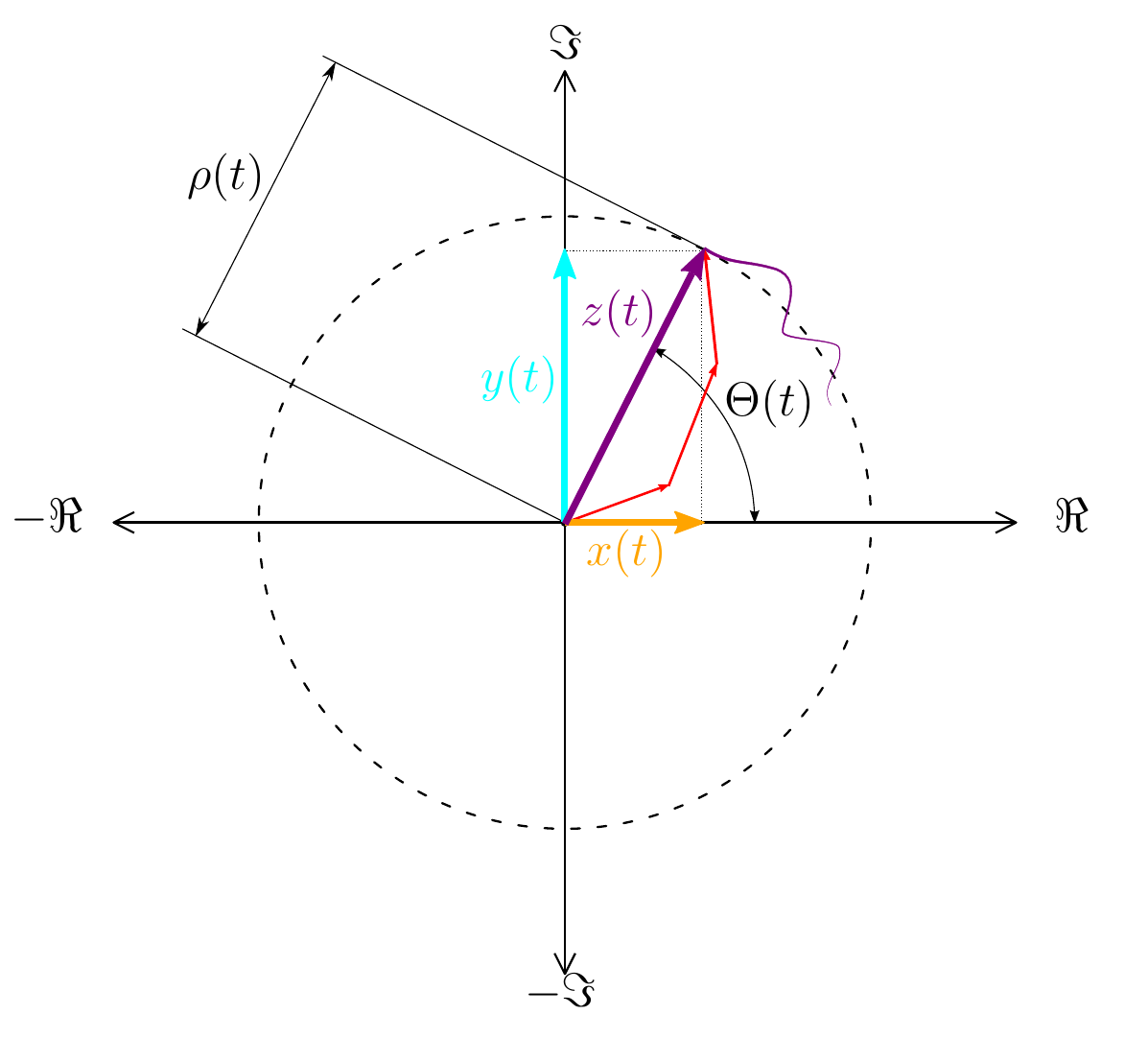}
  		\label{fig:bb}
  	}
  	\end{minipage}
	\caption{(a) Argand diagram of an AM--FM component in (\ref{eq:AMFMcomp}) at some time instant. Each component, $\psi(t)$ (\textcolor{MyRed}{\solidMrule[3.5mm]}) is interpreted as a vector: the IA $a(t)$ is interpreted as the component's length, the phase $\theta(t)$ is interpreted as a component's angular position. Although not shown, the IF $\omega(t)$ is interpreted as a component's angular velocity and phase reference $\phi$ is interpreted as an initial condition. The real part of the component $s(t)$ (\textcolor{MyBlue}{\solidMrule[3.5mm]}) and the component's quadrature $\sigma(t)$ (\textcolor{MyGreen}{\solidMrule[3.5mm]}) are interpreted as orthogonal projections of $\psi(t)$. We have included an example path (\textcolor{MyRed}{\solidSrule[3.5mm]}) taken by $\psi(t)$. (b) Argand diagram of the latent signal $z(t)$ in (\ref{eq:ComplexSignal}) at some time instant, composed of a superposition of components shown in red. The latent signal (\textcolor{MyPurple}{\solidMrule[3.5mm]}) is interpreted as a vector: the IA $\rho(t)$ is interpreted as the latent signal's length, the phase $\Theta(t)$ is interpreted as a latent signal's angular position, and although not shown the IF $\Omega(t)$ is interpreted as a latent signal's angular velocity. The real part of the signal $x(t)$ (\textcolor{myOrange}{\solidMrule[3.5mm]}) and the signal's quadrature $y(t)$ (\textcolor{myCyan}{\solidMrule[3.5mm]}) are interpreted as orthogonal projections of $z(t)$. We have included an example path (\textcolor{MyPurple}{\solidSrule[3.5mm]}) taken by $z(t)$.}
	\label{fig:ArgandComponent}
\end{figure}

When convenient, we will drop the $k$ denoting the parameters of the $k$th component for notational simplification and unspecified references to IA, IF, and phase refer to $a_k(t)$, $\omega_k(t)$, and $\theta_k(t)$ rather than $\rho(t)$, $\Omega(t)$, and $\Theta(t)$, respectively. We note that most of the paradoxes related to IF are a result of not correctly interpreting between these variables in the multicomponent case, because in the monocomponent case they are equivalent. The geometric interpretation of the AM--FM component in (\ref{eq:AMFMcomp}) is illustrated with the Argand diagram in Fig.~\ref{fig:ArgandComponent}. The AM--FM component can be visually interpreted as a single rotating vector in the complex plane with time-varying length and time-varying angular velocity.

Carson expressed the phase $\theta(t)$ in terms of a carrier frequency $\omega_{c}$ and FM message $m(t)$ as \cite{CarsonFry}
\begin{equation}
    \theta(t) = \omega_ct+\lambda\int\limits_{0}^{t}m(\uptau)d\uptau
    \label{eq:CarsonFryIP}
\end{equation}
assuming the modulation index, $\lambda\leq\omega_c$ and $|m(t)|\leq1$. We choose to parameterize the phase as
\begin{eqnarray}
    \theta(t) =  \omega_c t + \int\limits_{-\infty}^{t} m(\uptau)d\uptau + \phi
    \label{eq:SandovalIP}
\end{eqnarray}
which is more general than (\ref{eq:CarsonFryIP}), avoids an arbitrary normalization on $m(t)$, and does not impose an upper bound on the IF. The phase can also be expressed in terms of $\omega_{c}$ and Phase Modulation (PM) message $M(t)$ as
\begin{eqnarray}
    \theta(t) &=&  \omega_c t + M(t)  + \phi
\end{eqnarray}
where the relationship between the messages is given by
\begin{equation}
    M(t) = \int\limits_{-\infty}^{t}{m}(\uptau)d\uptau.
\end{equation}
In the AM--FM model, we assume non-negative IF
\begin{equation}
    {\omega}(t) = \frac{d}{dt}\theta(t)
    = \omega_c+  m(t) \geq0\ \forall \ t.
    \label{eq:NonNegInstFreq}
\end{equation}
In Section \ref{sec:Subtleties}, we discuss situations in which this assumption may be relaxed.  Finally, we assume without loss of generality, that the phase and carrier references are taken at $t=0$, i.e.
\begin{equation}
    \int\limits_{-\infty}^{0}{m}(t)dt=M(0)=0~~~\implies~~~ \phi = \theta(0)~~\text{and}~~\omega_c = \omega(0)-m(0).
\end{equation}

\subsection{Monocomponents and Narrowband Components} \label{ssec:Related}
The concept of a signal being composed of one (mono-) or more (multi-) components and whether those components are narrowband or wideband has caused much confusion in existing literature. There is no clear agreed upon definition for a monocomponent \cite{huang1998empirical}, however, a few definitions have been proposed \cite{boashash1992time,cohen1995time,boashash2003time}.  Some authors describe a monocomponent as defined by a single `ridge' in time and harmonic frequency, corresponding to an elongated region of energy concentration \cite{boashash1992time, cohen1992multicomponent, cohen1995time}. Cohen agrees that if a component is well localized, the crest of the ridge corresponds to the IF and further states that the width of the ridge depends on the energy spread, or instantaneous bandwidth of the component \cite{cohen1992multicomponent, cohen1995time}. A multicomponent signal may then be defined as any signal which is the sum of two or more monocomponents which can only occur if the separation between the ridges is large in comparison to bandwidth of the components \cite{boashash1992time,cohen1992multicomponent,cohen1995time, wei1998instantaneous}. It has also been noted that a signal may be monocomponent at some time instances and multicomponent at other time instances \cite{cohen1992multicomponent}.

With these descriptions of a monocomponent, we point out that the definitions and concepts of IA and IF of a complex AM--FM component as defined by (\ref{eq:AMFMcomp}) are by no means justified only for narrowband components. Also, as pointed out by Cohen \cite{cohen1995time}, a multicomponent signal is not defined by the ability to express a signal as a sum of parts, because there are an infinite number of ways to write a signal as a sum of parts. Rather, it is the nature of the parts in relation to themselves and the signal which determines whether the decomposition is of interest \cite{cohen1988instantaneous}. We strongly believe that the previous description of a monocomponent is unnecessarily restrictive, limiting the monocomponent to be narrowband. Clearly, the IA and IF in (\ref{eq:AMFMcomp}) are well-defined for wideband components.

Simply put, a perfectly-valid component can be defined with (\ref{eq:AMFMcomp}) while foregoing the restrictive narrowband definitions of Carson \cite{CarsonFry}, Ville \cite{ville}, Boashash \cite{BoashashP1,boashash1992time}, and Cohen \cite{cohen1995time}. Such wideband components may result in multiple ridges which may be broken down into narrowband components. However, for many problems a wideband component, such as the AM--FM component in (\ref{eq:AMFMcomp}), can be much more appropriate than multiple, narrowband components. Further, the appearance of structure in the Fourier spectrum can be viewed as an indication that a wideband component is present in the signal.  For narrowband components the Fourier and Hilbert spectra can be quite similar, but for wideband components these can be quite different.

\section{Hilbert Spectral Analysis} \label{sec:HSA}
Huang's original definition of the Hilbert spectrum uses Empirical Mode Decomposition (EMD) to determine a set of Intrinsic Mode Functions (IMFs) which are individually demodulated with the HT to obtain $\{a_k(t)\}$ and $\{\omega_k(t)\}$ \cite{huang1998empirical}.  This analysis is also known as the Hilbert-Huang Transform (HHT) \cite{huang1998empirical}. This definition can be generalized, by recognizing that IMFs are a class of AM--FM components and that other decomposition and demodulation methods exist for obtaining the instantaneous parameters, $\{a_k(t)\}$ and $\{\omega_k(t)\}$.  This generalization can lead to a more powerful and useful signal analysis technique as we will demonstrate.

Towards a more general definition, we define the Hilbert spectrum as a representation or characterization of a signal by an \emph{instantaneous spectrum} which is parameterized by AM--FM components as in (\ref{eq:AMFMcomp}). With this definition the Hilbert spectrum is not unique, therefore two problems must be addressed in order to obtain a unique solution:
\begin{itemize}[labelwidth = \widthof{~~~~~P1:},leftmargin=!]
    \item[\textbf{P1:}] Determine the appropriate complex extension to obtain an AM--FM decomposition with meaningful interpretation.
    \item[\textbf{P2:}] Estimate $\{a_k(t)\}$ and $\{\omega_k(t)\}$, $0 \leq k \leq K-1$.
\end{itemize}

The instantaneous parameters are defined on a complex signal, therefore to perform analysis on a real signal it must be extended to the complex space. There are cases where the extension is straightforward, i.e.~simple harmonic motion, or in communications where assumptions can be matched at the transmitter and receiver, or other cases where the assumptions are implicit, i.e. Fourier analysis. However, in the general case, signal decomposition and component demodulation are ambiguous and due to this, assumptions must always be made in order to arrive at a unique solution. We will address subtleties of complex extension of real signals assuming AM-FM components with relaxed HC in Section~\ref{sec:Subtleties}.

\begin{figure}[ht]
	\centering
	\includegraphics[scale=0.75]{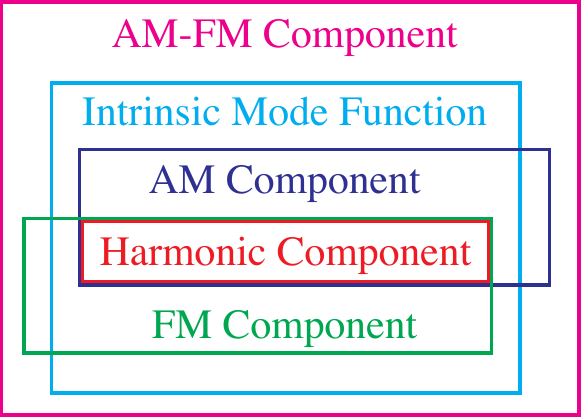} 
	\caption{The AM--FM component in (\ref{eq:AMFMcomp}) may be considered a generalization of other well-known components.  The familiar harmonic component may be considered special cases of the AM component and FM component. Huang's IMF is a special case of the AM--FM component.}
	\label{fig:Classes}
\end{figure}

The class of the solution obtained is dependent on the nature of the assumptions. Fig.~\ref{fig:Classes} illustrates the various classes of components under particular assumptions including the SHC, AM, FM, and AM--FM components as well as the IMF (discussed further in Part III). The AM--FM component in (\ref{eq:AMFMcomp}) may be considered a generalization of all of these components. In the next subsections, we illustrate several forms of the AM--FM component and AM--FM model under various assumptions.

\subsection{Two Conventions for Relating the Real Observation to the Latent Signal}

As noted by Gabor in \cite{gaborTOC}, there are two conventional ways to relate a real signal $x(t)$ to a complex signal $z(t)$. The first convention is
\begin{subequations}
	\label{eq:realofpsi}
	\begin{align}
	  	\label{eq:realofpsiA}
	  	x(t) &= \Re\lbrace z(t)\rbrace \\[-0.1em]
	  	\label{eq:realofpsiB}
         &= \Re\lbrace x(t)+jy(t)\rbrace
	\end{align}
\end{subequations}
which can be thought of as a single vector as in Fig.~\ref{fig:ArgandComponent} and the second convention is
\begin{subequations}
	\label{eq:origplusconj}
	\begin{align}
	  	\label{eq:origplusconjA}
	  	x(t) &=\left[ {z}(t)+{z}^*(t)\right]/2 \\[-0.1em]
	  	\label{eq:origplusconjB}
         &= [ x(t)+jy(t) +  x(t)-jy(t)]/2
	\end{align}
\end{subequations}
which can be thought of as two vectors. If $z(t)$ is a superposition of AM--FM components, the first convention can be viewed as a \emph{single} rotating vector per component, while the second convention can be viewed as a \emph{pair} of vectors per component rotating in opposite directions.

Equation (\ref{eq:origplusconj}) is implicit in standard Fourier analysis, with the interpretation that the latent signal is real-valued and its spectrum\footnote{In this work, use of the unspecified word ``spectrum'' strictly means Fourier spectrum and not ``Hilbert spectrum.''} consists of Hermitian symmetric, positive and negative frequency components. The convention in (\ref{eq:realofpsi}) is adopted for the AM--FM model in (\ref{eq:AMFMmodel}), with the interpretation that only the real part of the signal, $x(t)$ is observed.  Inherent in both these conventions is the ambiguity associated with the quadrature $y(t)$, i.e.~an infinite number of choices exist for $y(t)$ that lead to the same $x(t)$. Thus, we view the quadrature $y(t)$ as a \emph{free parameter} that can be \emph{chosen} to achieve alternate interpretations.

We use the term quadrature in a general sense. For mathematical and physical reasons, the quadrature signal is most often chosen offset in phase by one-quarter cycle ($\pi/2$ radians) relative to the real signal. More specifically, if we assume a single SHC, $z(t)=\psi_0(t;a_0,\omega_0,\phi_0)$ then the AM--FM model is given by~(\ref{eq:shc}). Although choosing the quadrature signal given the real signal is trivial for the SHC, the same is not true for the component in (\ref{eq:AMFMcomp}) where the assumptions $a_0(t)=a_0$, and $\omega_0(t)=\omega_0$ are removed, i.e.~the concept of a quadrature signal for the AM--FM component is ambiguous---as will be demonstrated in the following sections.

\subsection{Simple Harmonic Component}
The simplest form of the AM--FM model has $K=1$, $a_0(t) = a_0$, and  $\omega_0(t)= \omega_0$ in which case the AM--FM component in (\ref{eq:AMFMcomp}) becomes
\begin{equation}
  		\psi_0(t;a_0,\omega_0,\phi_0) = a_0 e^{j( \omega_0 t + \phi_0)}.
	\label{eq:simpleHamonicComponent}
\end{equation}
We recognize the significance of this form as representing simple harmonic motion as we discussed in Part I as Vakman's Condition 3.

\subsection{Superposition of Simple Harmonic Components}
Perhaps the most familiar form of the AM--FM model has $K=\infty$, $a_k(t) = a_k$, and $\omega_k(t) = k\omega_0$
\begin{equation}
	z(t) = \sum\limits_{k=0}^{\infty} a_k e^{j (k \omega_0 t+\phi_k)}
	\label{eq:FourierSeriesModel}
\end{equation}
and has real observation given by (\ref{eq:RealObservation}) repeated here as
\begin{subequations}
	\label{eq:FourierSeries}
	\begin{align}
	  	\label{eq:FourierSeriesA}
	  	x(t) &= \Re\left\lbrace z(t) \right\rbrace\\[-0.1em]
	  	\label{eq:FourierSeriesB}
	  	&= \sum\limits_{k=0}^{\infty} a_k  \cos(k \omega_0 t+\phi_k)\\[-0.5em]
	  	\label{eq:FourierSeriesC}
	  	&= A_0 +\sum\limits_{k=1}^{\infty} [A_k \cos(k \omega_0 t) + B_k \sin(k \omega_0 t)].
	\end{align}
\end{subequations}
Equation (\ref{eq:FourierSeries}) is recognized as a Fourier Series (FS) with a particular convention used for the complex extension \cite{fourier1807memoire, Fourier1822analytical}. Note that $\phi_k$ is required for proper \emph{synthesis} since without it, $A_k=a_k$ and $B_k=0$ and the resulting $x(t)$ is always even, which may not be true. The FS can be considered the simplest of AM--FM superposition models corresponding to the assumption of SHCs.

The FT is the limiting form of the FS as the fundamental, $\omega_0\rightarrow 0$ \cite{OppenheimWillsky}. When viewing the FT as a special case of the AM--FM model, subtle and important observations can be made:
\begin{itemize}
    \item In the decomposition problem \textbf{P1}, the AM--FM model corresponding to a separable solution with $X(\omega)$ (not time-varying) and $e^{j\omega t}$ (constant frequency), \emph{is} the FT, i.e.~the inverse FT is a superposition (inner product) of $X(\omega)$ and $e^{j\omega t}$. 

    \item The FT yields a solution that can be described as a superposition of non-time-varying components described by constant state parameters, $a_k(t)=a_k$, $\omega_k(t)=k\omega_0$, and $\phi_k$ which corresponds to a constant model state. However, this does not in any practical way constrain the real superposition $x(t)$. Although analysis of a time-varying system can be accomplished using a constant model, like the FT, in most cases this forces $K=\infty$ which may not accurately describe a physical system.

    \item HSA can be considered a generalization of Fourier analysis \cite{huang1998empirical}. However, rather than a generalization of Fourier analysis obtained by generalizing a kernel function (Gabor transform, Wigner-Ville distribution, wavelet transform, etc.) \cite{Allen1977ShortTime,  cohen1989time, jones1990instantaneous, cohen1993anything,lovell1993relationship, cohen1995time, QianJoint}, an AM--FM component can be viewed as allowing for \emph{additional} degrees of freedom in the bases of Fourier analysis. Although it may be convenient to think of (\ref{eq:AMFMcomp}) as a basis for HSA because they span the entire Hilbert space, they do not meet the linear independence property of a formal basis \cite{strang09}.
\end{itemize}

\subsection{The AM Component}
One way to generalize the SHC in (\ref{eq:simpleHamonicComponent}) is to relax the constraint of constant amplitude which leads to an AM component
\begin{eqnarray}
	\psi_0(t;a_0(t),\omega_0,\phi_0) &=&  a_0(t)e^{j(\omega_0 t + \phi_0)}.
	\label{eq:AMcomp}
\end{eqnarray}
The simple AM component was originally developed in the context of communication theory \cite{armstrong1915some, voelcker1966towardI, voelcker1966towardII}.

\subsection{Superposition of AM Components}
Another form of the AM--FM model is a superposition of AM components that have frequencies as integer multiples of a fundamental
\begin{equation}
	z(t) = \sum\limits_{k=0}^{\infty} a_k(t) e^{j (k \omega_0 t+\phi_k)} 
    \label{eq:AMSuperposition}
\end{equation}
with real observation given by (\ref{eq:RealObservation}).

The AM superposition was first conceived by Gabor and lead to the development of the Short-Time Fourier Transform (STFT) \cite{gaborTOC} 
\begin{subequations}
	\label{eq:AMsuperposition}
	\begin{align}
	  	\label{eq:AMsuperpositionA}
	  	X_t(\omega) &= \int x(\uptau)h(\uptau-t) e^{-j\omega \uptau}d\uptau\\[-0.75em]
  \notag&\updownarrow \\[-0.75em]
	  	\label{eq:AMsuperpositionB}
        x_\omega(t) &= \frac{1}{2\pi} \int X(\text{\textcloseomega}) H(\omega-\text{\textcloseomega})e^{j\text{\textcloseomega} t} d\text{\textcloseomega}.
	\end{align}
\end{subequations}
where $h(t)\leftrightarrow H(\omega)$ is a window function. There are two possible interpretations of the STFT in context of the AM--FM model:
\begin{enumerate}
    \item The window function $h(t)$ is applied to the signal $x(t)$ in order to use the FT, i.e.~simple harmonic analysis at each $t$. This interpretation couples $H(\omega)$ to $X(\omega)$.
    \item The window function $h(t)$ is applied to the $e^{-j\omega t}$ and is effectively an AM superposition analysis where the window is an assumption on $a_k(t)$ in the AM--FM model. This interpretation couples $H(\omega)$ to $e^{j\omega t}$.
\end{enumerate}
The second interpretation provides an important and alternate view of the STFT. What is typically viewed as imprecision in the ability to compute time-\emph{frequency} parameters in the first interpretation can alternatively be viewed in the second interpretation as a precise ability to compute time-\emph{component} parameters using the modified basis, $h(t-\uptau) e^{-j\omega t}$.

\subsection{FM Component}
Another way to generalize the SHC in (\ref{eq:simpleHamonicComponent}) is to relax the constraint of constant frequency which leads to a simple FM component
\begin{eqnarray}
	\psi_0(t;a_0,\omega_0(t),\phi_0) &=&  a_0e^{j[\int\limits_{-\infty}^{t}\omega_0(\uptau) d\uptau + \phi_0]}.
	\label{eq:FMcomp}
\end{eqnarray}
The FM component was also developed in the context of communication theory \cite{Carson1922Notes,CarsonFry}. The FM component has been used in signal synthesis where Chowning observed that a single, wideband FM component is perceived by the ear as spectrally rich, i.e.~multiple SHCs \cite{FMSynthesis}. This FM-synthesis method has been used in commercial audio synthesizers \cite{bedaux1974micro,chowning1977synthesis}.

\subsection{Superposition of FM Components}
Signal analysis using a superposition of FM components is not usually considered due to the loss of linear independence of the components in the model and the restriction of constant amplitude. For this reason, signal analysis using a superposition of AM--FM components is more useful than a superposition of FM components.

\subsection{Other AM--FM Models}
\label{ssec:Previous}
Alternatives to the proposed AM--FM model in Section \ref{sec:AMFMmodel} have been considered.  However, common to these are restrictive assumptions which limit the utility as we highlight below. Previous AM--FM models for signal analysis/synthesis usually fall into one of three main groups: 1) HT \cite{feldman1994non,rao2000decomposing,gianfelici2007multicomponent,feldman2011hilbert}, 2) peak tracking/sinusoidal modeling \cite{mcaulay1986speech,rao1990estimation,Pantazis2011,Boashash2013}, and 3) Teager energy operator \cite{maragos1993energy,bovik1993fm,fertig1996instantaneous,potamianos1999speech,boudraa2004if,boudraa2011instantaneous}. However, some models exist that do not fall into any of these groups \cite{Quatieriseparation,fertig1996instantaneous}. A historical summary of AM--FM modeling is presented by Gianfelici \cite{gianfelici2007StateOfArt}. A review of algorithms for estimating IF is presented by Boashash \cite{boashash1990algorithms,BoashashP2}.

The HT permits the unambiguous definition of IA, IF, and phase of any real signal (random or deterministic) \cite{vakmanVainshtein1977}. Thus, it provides a direct means of performing AM--FM analysis, however, it requires that the quadrature signal be defined in a consistent manner in all cases. Implicit in the HT is HC, as a result the more likely this assumption is true, the better the Hilbert-transformed signal approximates the true quadrature signal and the more likely the HT-based solution provides an accurate model of the underlying physical synthesis model \cite{boashash1992time}. Direct application of the HT does not however address the signal decomposition problem, as a result, methods for using the HT for decomposition have been proposed \cite{gianfelici2007multicomponent,feldman2008compare,feldman2011hilbert}.

In peak tracking, one accepts the narrowband definition of a component and as a result, each component appears in the time-frequency plane as a single ridge of energy concentration. Thus, a signal can be parameterized by tracking its ridges in location, intensity, and possibly bandwidth. The time-frequency distribution is usually derived from a STFT \cite{mcaulay1986speech} but a generalized time-frequency distribution can also be used \cite{rao1990estimation}. Regardless of the time-frequency distribution used, the narrowband assumption is inherent. Extensions of the sinusoidal model have been proposed, such as the harmonic plus noise, and adaptive quasi-harmonic model \cite{Pantazis2011}.

The Teager Energy Separation Algorithm (ESA) provides a method of estimating the IA and IF of an AM--FM component, assuming HC. The Teager Energy Operator (TEO) is defined as \cite{quatieri02}
\begin{eqnarray}
	\varPsi\left\{x(t)\right\} &\equiv& \left[\dot{x}(t)\right]^{2}-x(t)\ddot{x}(t)
\end{eqnarray}
where $\dot{x}(t)$ and $\ddot{x}(t)$ are the first and second time derivatives of $x(t)$, respectively. The TEO applied to a single narrowband AM--FM component in (\ref{eq:AMFMcomp}), results in \cite{quatieri02}
\begin{eqnarray}
	\varPsi\left\{\psi_0(t)\right\} &\approx& \left[a_0(t)\omega_0(t)\right]^{2}.
\end{eqnarray}
Assuming the modulating signals $a_0(t)$ and $m_0(t)$ are bandlimited, it can be shown that
\begin{eqnarray}
	a_0(t) &\approx& \dfrac{\varPsi\left\{\psi_0(t)\right\}}{\sqrt{\varPsi\left\{\dot{\psi_0}(t)\right\}}}
\end{eqnarray}
and
\begin{eqnarray}
	\omega_0(t) &\approx& \sqrt{\dfrac{\varPsi\left\{\dot{\psi_0}(t)\right\}}{\varPsi\left\{\psi_0(t)\right\}}}
\end{eqnarray}
thus providing a method to estimate narrowband monocomponent IA/IF parameters \cite{quatieri02, potamianos1994comparison, diop2011joint, santhanam2000multicomponent}. 

The AM--FM model in these groups \emph{all} rely on a rigid narrowband component.  Surprisingly, the use of wideband components is a well-known means of synthesis \cite{moorer1976synthesis, chowning1977synthesis, FMSynthesis, dodge1997computer}. However, \emph{analysis} using the wideband components in (\ref{eq:AMFMcomp}) in the general form, has not been considered. It is our belief that the true power of HSA can only be fully recognized with the use of wideband components without the inherent assumption of HC.

\section{Frequency Domain View of Latent Signal Analysis} \label{sec:FDVoLSA}
In Part I, we introduced the LSA problem in the time domain.  The frequency domain view of the LSA problem is to estimate the latent spectrum $Z(\omega)$ from $X(\omega)$. Because $x(t)=\Re\{z(t)\}$, this imposes structure on the spectrum of $x(t)$, namely $X(\omega)=[Z(\omega)+Z^*(-\omega)]/2$, i.e.~Hermitian symmetry. As discussed in Part I, the conventional way to estimate $z(t)$ from $x(t)$ is via the HT as in (\ref{eq:ASdef}). Equivalently in the frequency domain, the conventional way to estimate $Z(\omega)$ from $X(\omega)$ is with Gabor's QM. Using Gabor's QM effectively forces the FT of (\ref{eq:HTRule}) $\hat{Y}(\omega) = -j\mathrm{sgn}(\omega)X(\omega)$, and $\hat{Z}(\omega)$ is always spectrally one-sided \cite{Boashash2013}. In the frequency domain, by relaxing HC we gain the freedom to choose $\hat{Y}(\omega)$ appropriately.

As an example of this frequency domain view, consider the observed spectrum
\begin{equation}
    X(\omega) = a_0 \pi [\delta(\omega+\omega_0) + \delta(\omega-\omega_0)]
\end{equation}
where $\delta(\cdot)$ is the Dirac delta function. If we assume HC,
\begin{equation}
    \hat{Y}(\omega) = -j a_0 \pi [-\delta(\omega+\omega_0) + \delta(\omega-\omega_0) ]
\end{equation} 
and the corresponding latent spectrum is given by
\begin{equation}
    \begin{split}
            \hat{Z}(\omega) &=  X(\omega) +j \hat{Y}(\omega)   \\
                            &=  2 a_0  \pi \delta(\omega-\omega_0).
    \end{split}
    \label{eq:exFT1}
\end{equation} 
If we do not assume HC, there are many choices for $\hat{Y}(\omega)$ leading to many possible latent spectra including,
\begin{equation}
            \hat{Z}(\omega) = 2 a_0  \pi \delta(\omega+\omega_0),
            \label{eq:exFT2}
\end{equation}
\begin{equation}
            \hat{Z}(\omega) = a_0 \pi [(1-\alpha)\delta(\omega+\omega_0) + (1+\alpha)\delta(\omega-\omega_0)],
\label{eq:exFT3}
\end{equation} 
and
\begin{equation}
            \hat{Z}(\omega) = a_0 \pi  [\delta(\omega+\omega_0) + \delta(\omega-\omega_0) + \alpha\delta(\omega-\beta\omega_0) - \alpha\delta(\omega+\beta\omega_0)] 
\label{eq:exFT4}
\end{equation} 
where $\alpha, \beta \in\mathbb{R}$ and (\ref{eq:exFT4}) corresponds to the FT of (\ref{eq:NotHarmonicCoor}). Because of the structure imposed on the spectrum $X(\omega)=[Z(\omega)+Z^*(-\omega)]/2$, the latent spectra in (\ref{eq:exFT1})-(\ref{eq:exFT4}) all yield the same $X(\omega)$.

In Table \ref{tab:3structures} we summarize the spectral structure of a complex signal depending on whether or not HC is assumed. If $z(t)$ takes on only real values, then the spectrum has Hermitian symmetry (first row). If $z(t)$ takes on complex values and has HC, then the Fourier spectrum is single-sided (second row). If $z(t)$ takes on complex values and does not have HC, the spectrum does not have Hermitian symmetry.

\begin{table}[ht]
\caption{Structure of the Fourier spectrum under the assumption of a model composed of Simple Harmonic Components (SHCs).  For a complex signal that takes on only real values, the Fourier spectrum is Hermitian symmetric.  For a complex signal with HC, the Fourier spectrum is one-sided.{\protect\footnotemark} For a complex signal without HC, the Fourier spectrum is not Hermitian symmetric.}
\begin{center}
\begin{tabular}{ll}
\toprule
\head{Class of  $z(t)$} & \head{FT Structure} \\
\midrule
$z(t)\in\mathbb{R}$ & $Z(\omega)=Z^*(-\omega)$\\
$z(t)\in\mathbb{C}$ w/ HC & $Z(\omega)=0,~\omega<0$\\
$z(t)\in\mathbb{C}$ w/o HC & $Z(\omega)\neq Z^*(-\omega)$\\
\bottomrule
\end{tabular}
\end{center}
\label{tab:3structures}
\end{table}%
\footnotetext{For convenience, we omit the case that $Z(\omega)=0,~\omega>0$.}

\begin{figure}[ht]
	\centering
	
		\centering
		\subfigure[]{
			\begin{minipage}[b]{\linewidth}	
				\begin{minipage}[b]{0.33\linewidth}
					\centerline{\includegraphics[width = 0.95\linewidth]{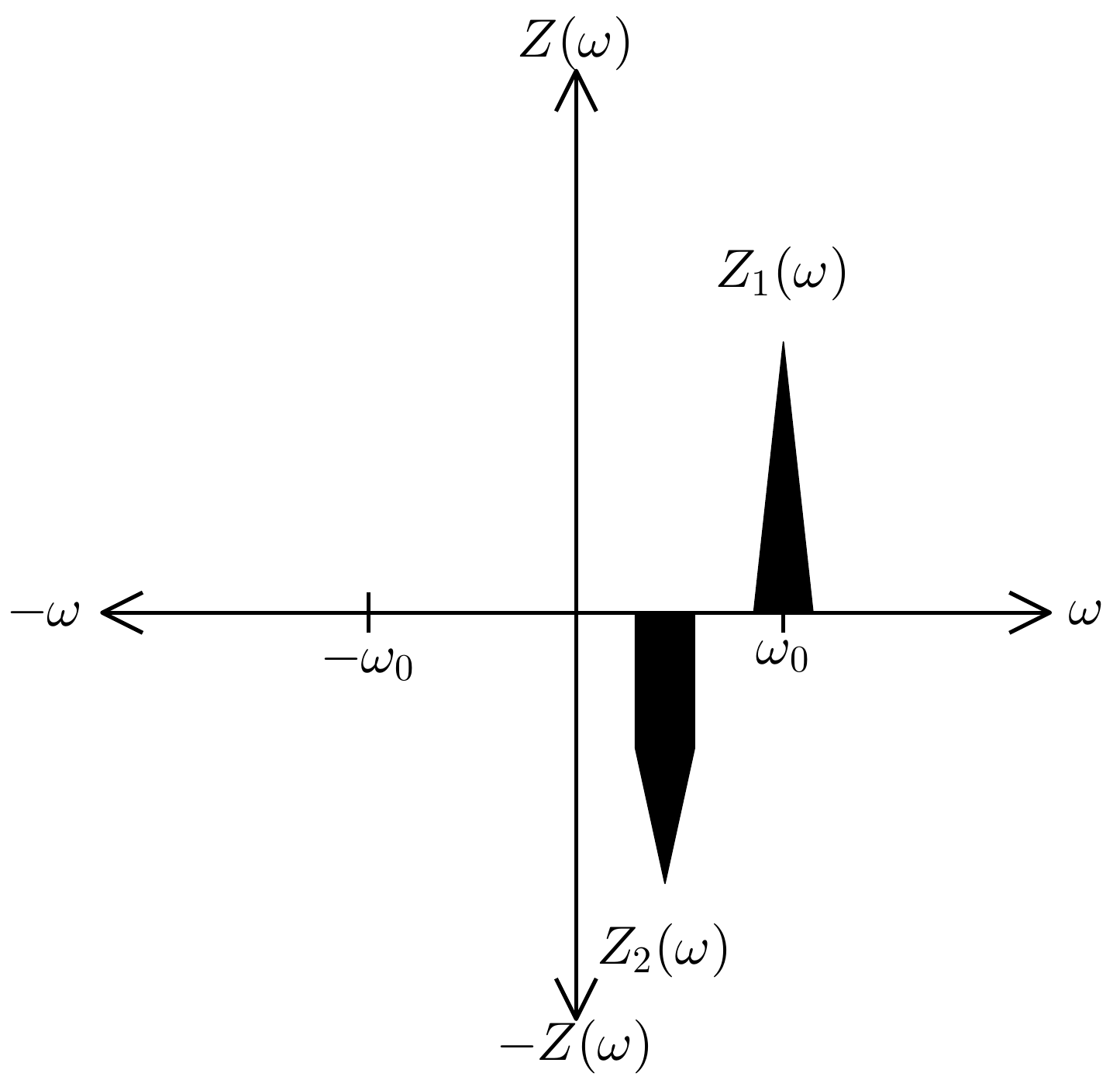}}
					\centerline{\footnotesize Latent Spectrum}			
				\end{minipage}\begin{minipage}[b]{0.33\linewidth}
					\centerline{\includegraphics[width = 0.95\linewidth]{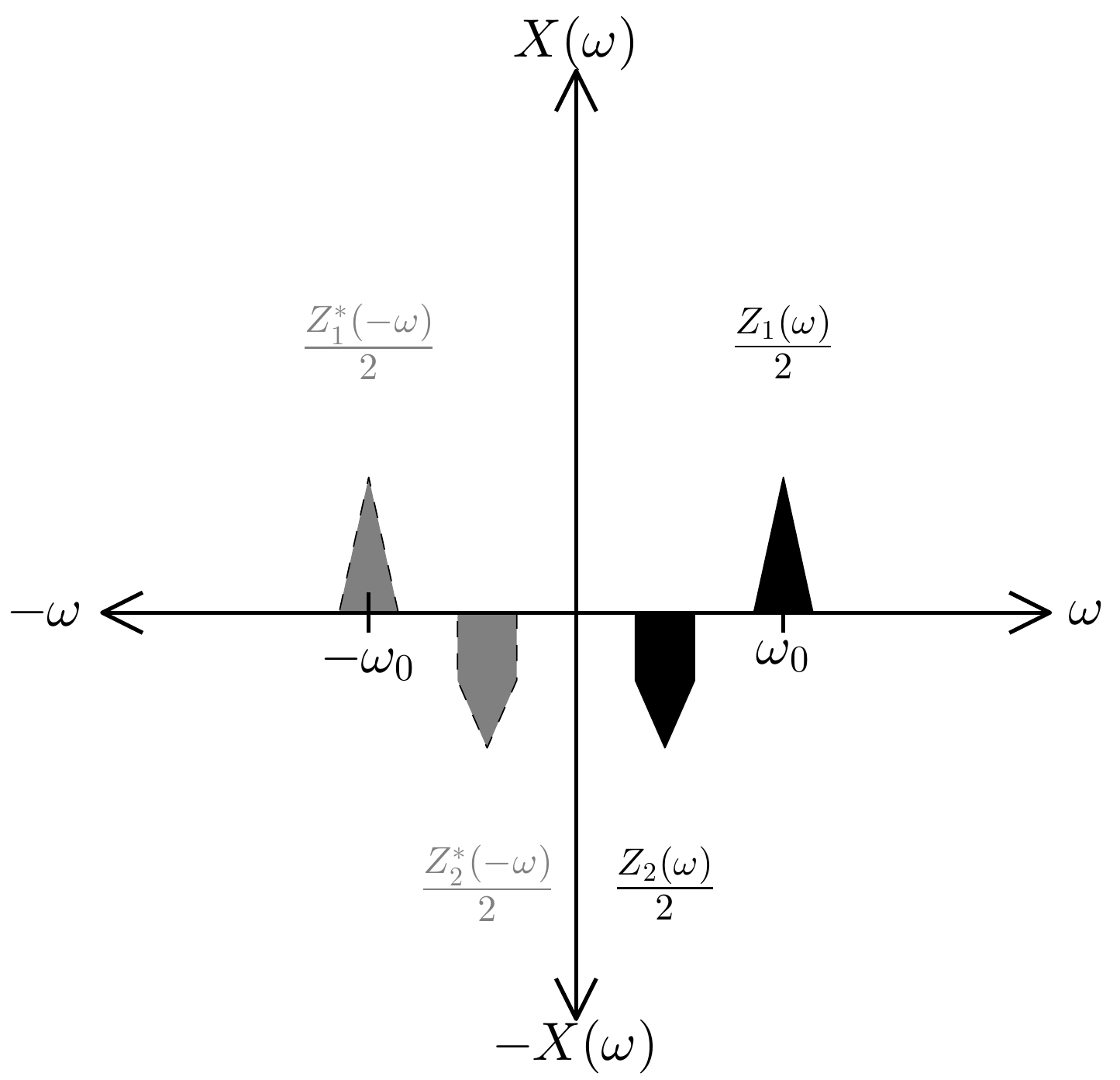}}
					\centerline{\footnotesize Observed Spectrum}		
				\end{minipage}\begin{minipage}[b]{0.33\linewidth}
					\centerline{\includegraphics[width = 0.95\linewidth]{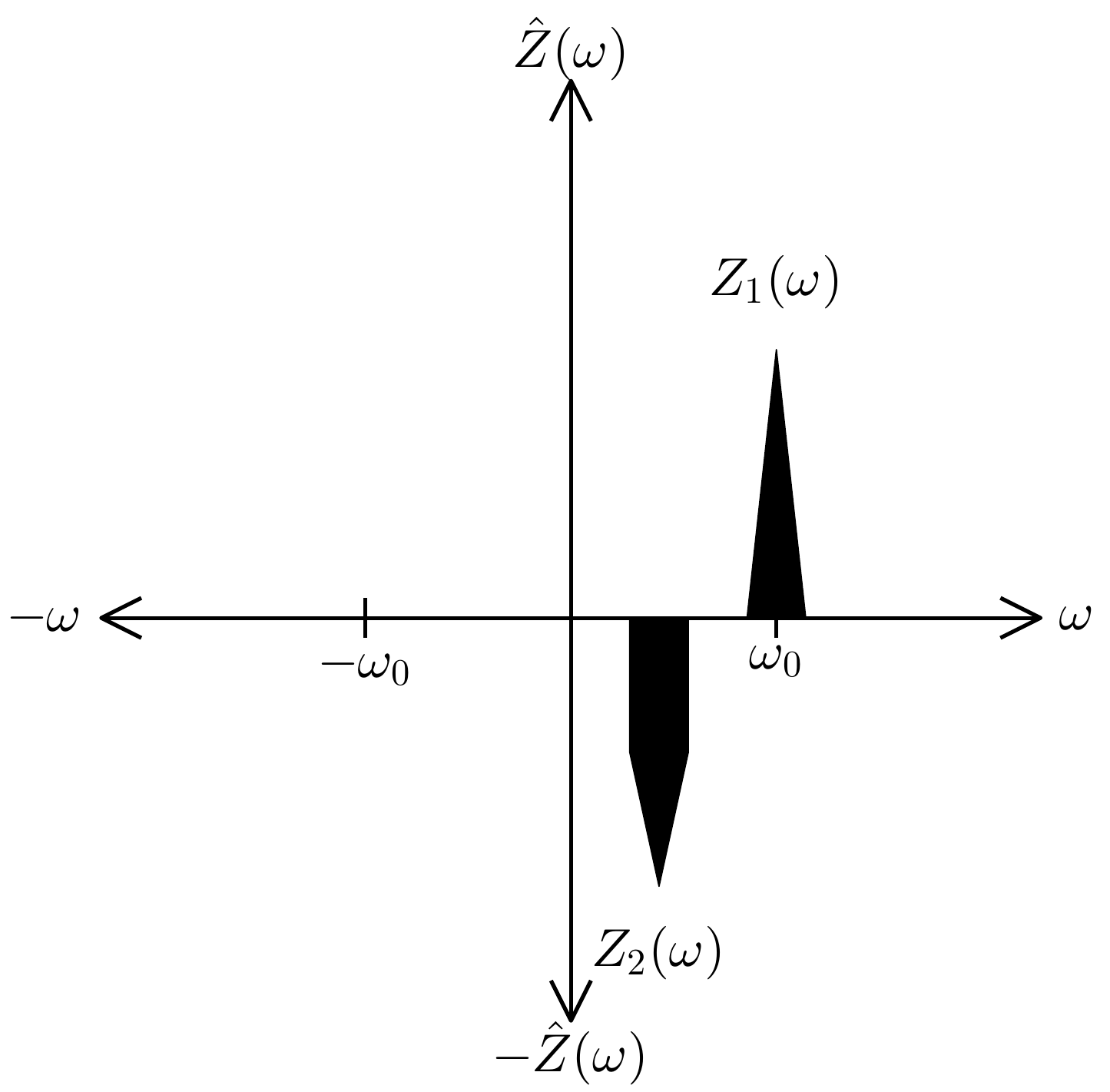}}
					\centerline{\footnotesize Gabor's AS Spectrum}		
				\end{minipage}
			\end{minipage}
			\label{fig:IllustrationA}
		}\\
		\subfigure[]{
					\begin{minipage}[b]{\linewidth}	
						\begin{minipage}[b]{0.33\linewidth}
							\centerline{\includegraphics[width = 0.95\linewidth]{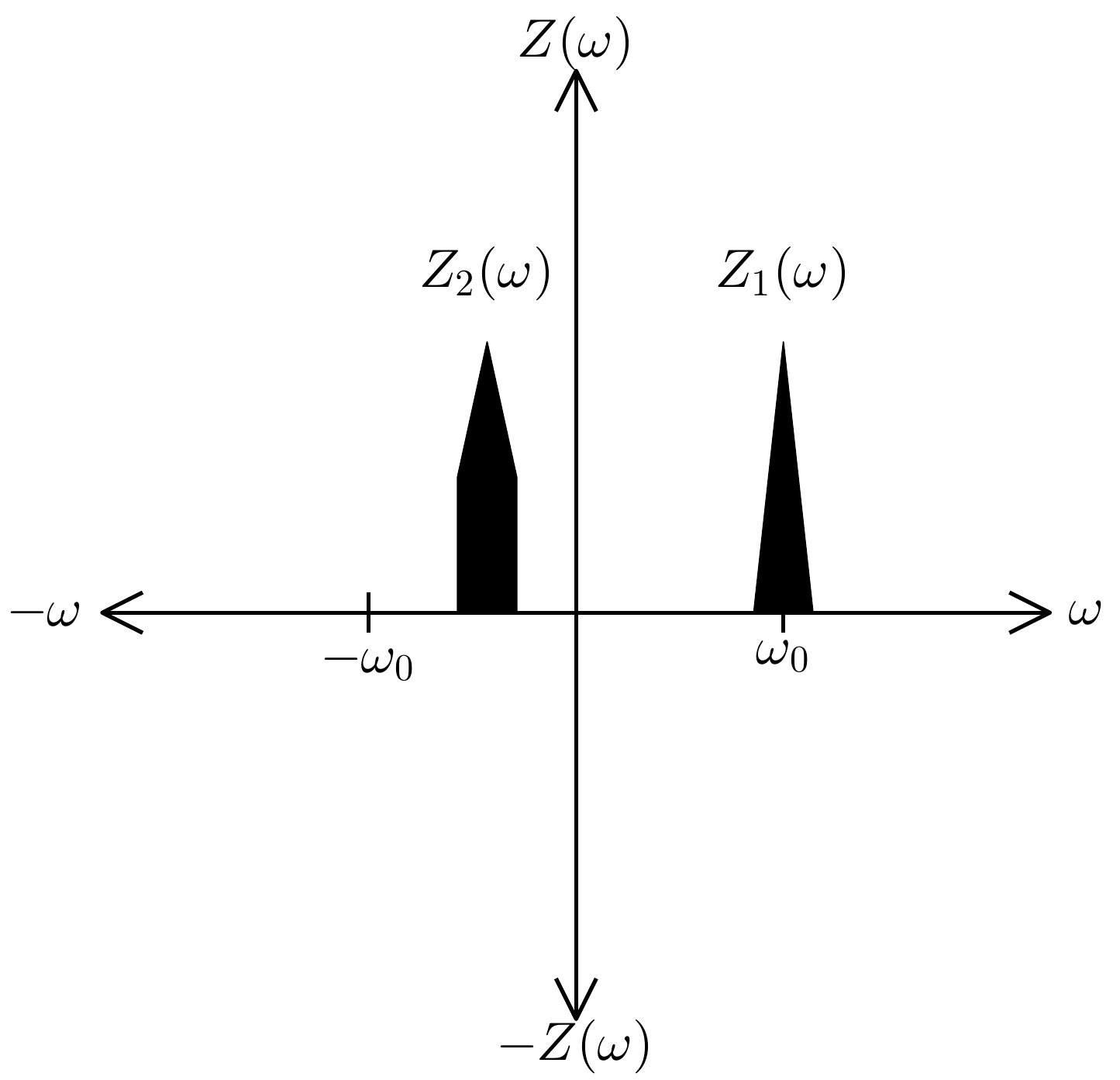}}
							\centerline{\footnotesize Latent Spectrum}			
						\end{minipage}\begin{minipage}[b]{0.33\linewidth}
						\centerline{\includegraphics[width = 0.95\linewidth]{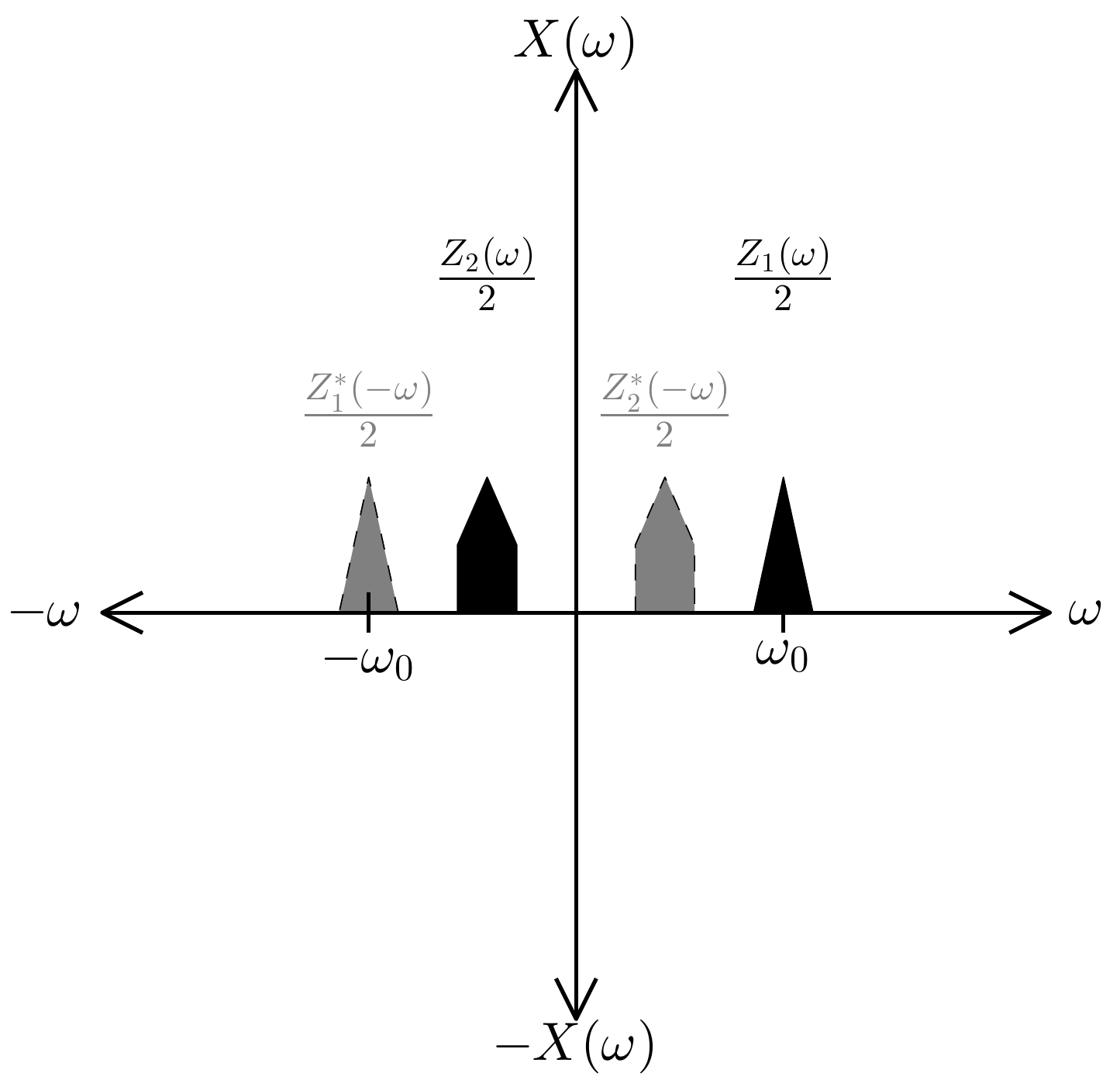}}
						\centerline{\footnotesize Observed Spectrum}		
					\end{minipage}\begin{minipage}[b]{0.33\linewidth}
					\centerline{\includegraphics[width = 0.95\linewidth]{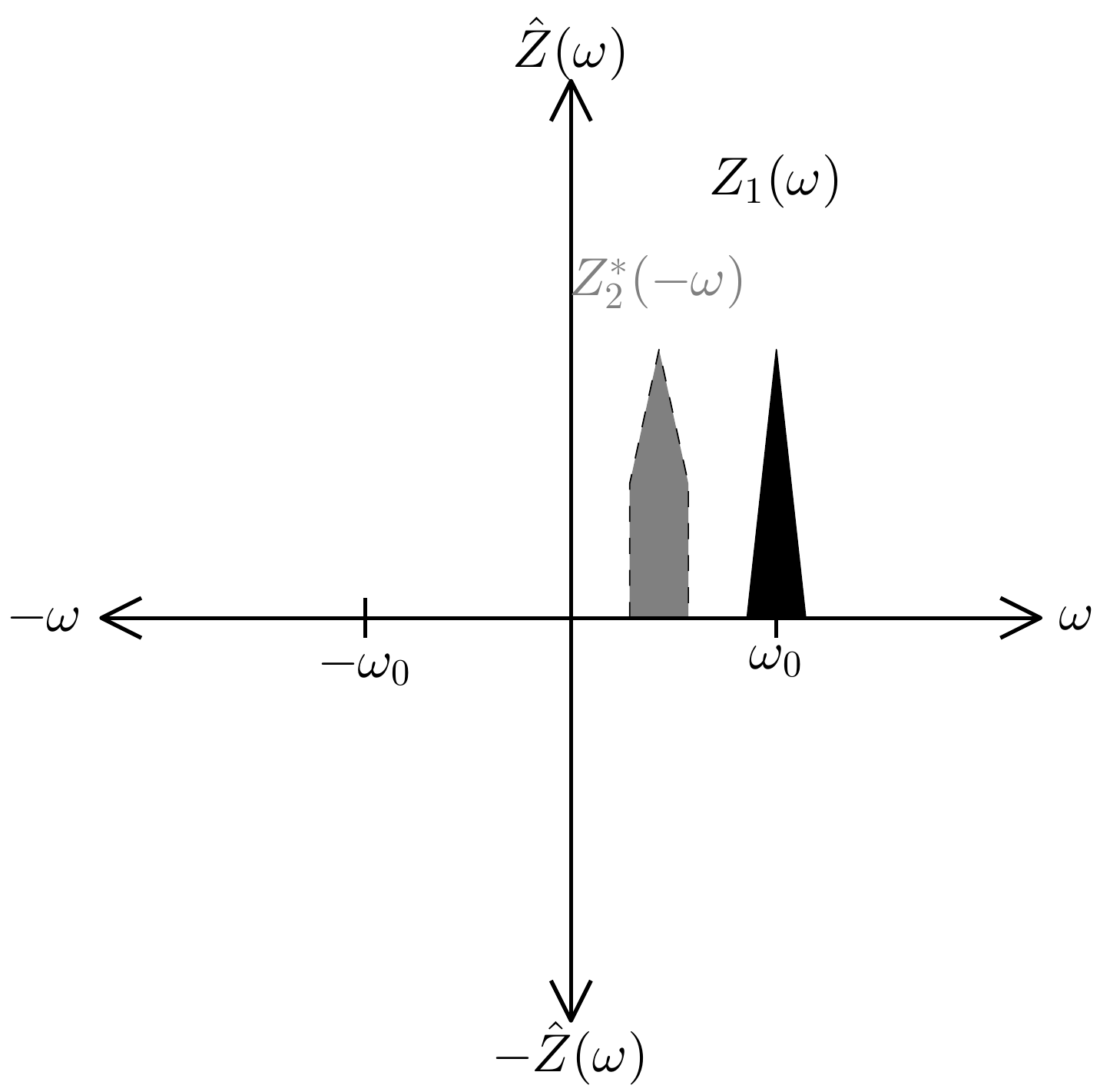}}
					\centerline{\footnotesize Gabor's AS Spectrum}		
				\end{minipage}
			\end{minipage}
			\label{fig:IllustrationB}
		}\\
		\subfigure[]{
						\begin{minipage}[b]{\linewidth}	
							\begin{minipage}[b]{0.33\linewidth}
								\centerline{\includegraphics[width = 0.9\linewidth]{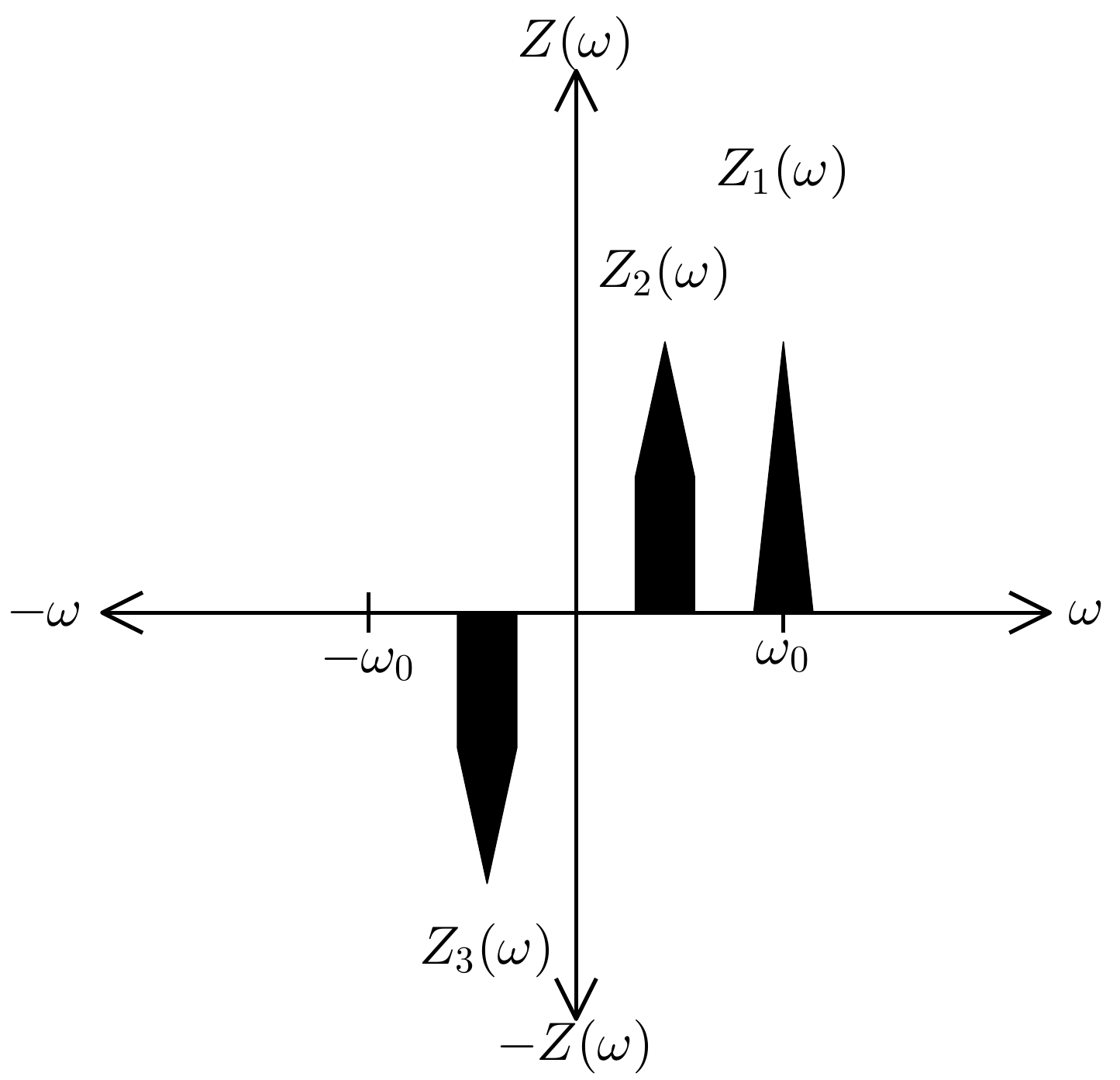}}
								\centerline{\footnotesize Latent Spectrum}			
							\end{minipage}\begin{minipage}[b]{0.33\linewidth}
							\centerline{\includegraphics[width = 0.9\linewidth]{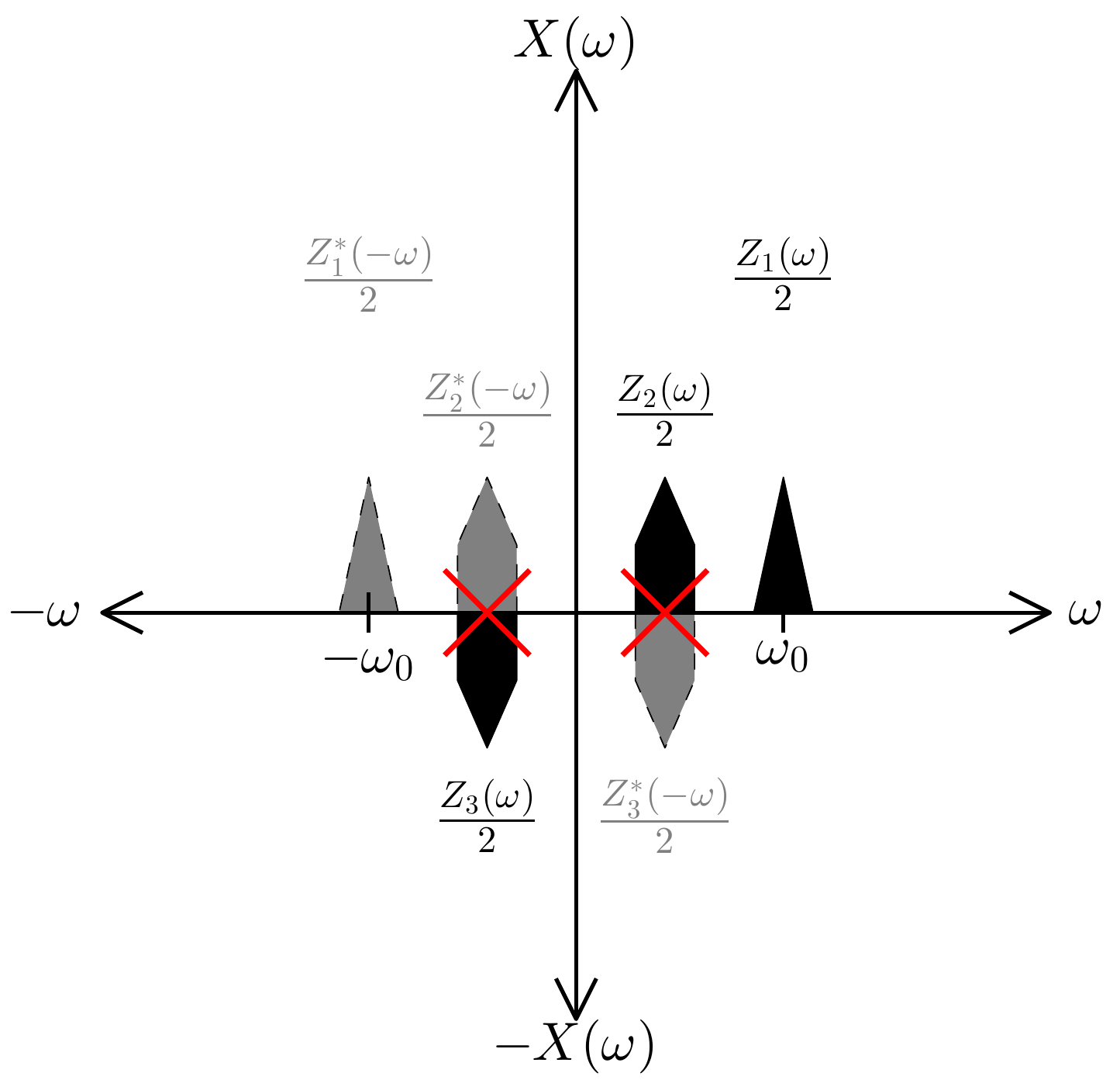}}
							\centerline{\footnotesize Observed Spectrum}		
						\end{minipage}\begin{minipage}[b]{0.33\linewidth}
						\centerline{\includegraphics[width = 0.9\linewidth]{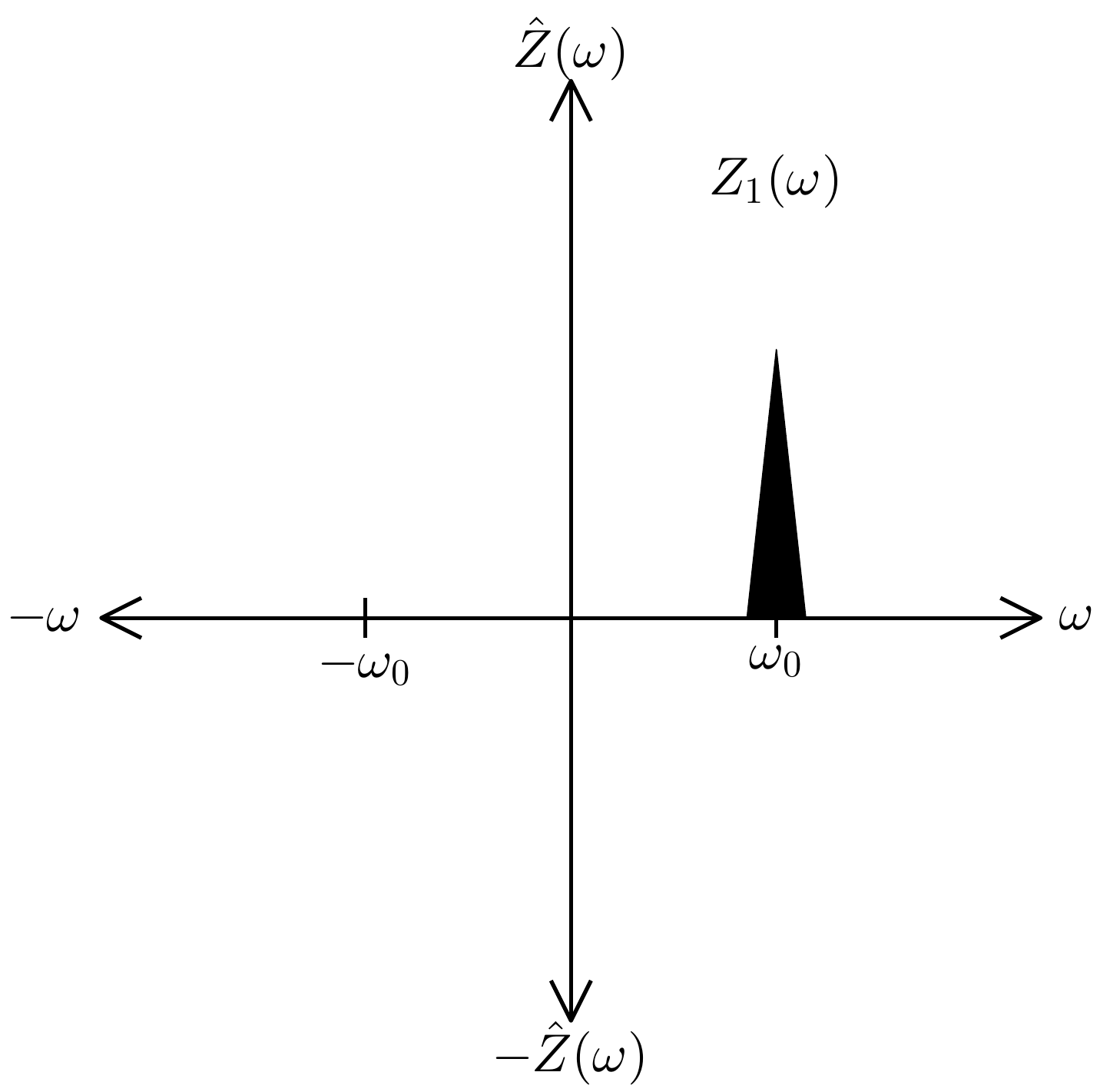}}
						\centerline{\footnotesize Gabor's AS Spectrum}		
					\end{minipage}
				\end{minipage}
				\label{fig:IllustrationC}
			}
	\caption{Illustrations of when Gabor's QM (a) can distinguish $Z(\omega) = Z_1(\omega) + Z_2(\omega)$ from $Z^*(-\omega) = Z_1^*(-\omega) + Z_2^*(-\omega)$ because $Z(\omega)$ has HC, (b) cannot distinguish $Z(\omega)=Z_1(\omega)+Z_2(\omega)$ from $Z^*(-\omega) = Z_1^*(-\omega)+Z_2^*(-\omega)$ because the latent spectrum is two-sided and Hermitian symmetry is imposed by the real observation, and (c) is incorrect because the structure of the latent spectrum along with the Hermitian symmetry imposed by the real observation has \emph{concealed} terms.}
	\label{fig:Illustration}	
\end{figure}

In the frequency domain view of LSA, Gabor's QM can be interpreted as simply a method to distinguish $Z(\omega)$ from $Z^*(-\omega)$ which works only if the harmonic frequencies of $Z(\omega)$ are all non-negative. This is a convenient ``trick'' and a generalization of this concept for the Hilbert spectrum is possible. As an illustration, consider Fig.~\ref{fig:IllustrationA} where the latent spectrum $Z(\omega) = Z_2(\omega)+Z_1(\omega)$ and the real signal under analysis has a spectrum $[Z_1^*(-\omega)+Z_2^*(-\omega)+Z_2(\omega)+Z_1(\omega) ]/2$. Applying Gabor's QM, the complex signal clearly has the spectrum $Z(\omega) = Z_2(\omega)+Z_1(\omega)$. Next consider Fig.~\ref{fig:IllustrationB} where the latent spectrum $Z(\omega) = Z_2(\omega)+Z_1(\omega)$ and the real signal under analysis has a spectrum $[Z_1^*(-\omega)+Z_2(\omega)+Z_2^*(-\omega)+Z_1(\omega) ]/2$. Applying Gabor's QM, the complex signal has the spectrum $Z_2^*(-\omega)+Z_1(\omega)$ which is the incorrect spectral grouping. Finally, consider Fig.~\ref{fig:IllustrationC} where the latent spectrum $Z(\omega) = Z_3(\omega)+Z_2(\omega)+Z_1(\omega)$ and the real signal under analysis has a spectrum $[Z_1^*(-\omega)+ \cancel{Z_2^*(-\omega) + Z_3(\omega)} + \cancel{Z_3^*(-\omega) + Z_2(\omega)} + Z_1(\omega)]/2$; cancellation is due to symmetry of the latent spectrum. Applying Gabor's QM, the complex signal has the spectrum $Z_1(\omega)$ and not only has Gabor's QM failed to yield the latent spectrum but terms are missing entirely.

\section{Subtleties of the Hilbert Spectrum} \label{sec:Subtleties}
In the previous section, we considered the LSA problem where we assumed SHCs but relaxed the HC assumption. By relaxing HC, we also relax the non-negative IF constraint associated with it. In this section, we discuss the implications of assuming non-negative IF and relaxing the assumption of SHCs.

First, consider making the assumption of non-negative IF and a model composed of SHCs. Table \ref{tab:3structures}, second row shows this implies HC on $z(t)$. Next, consider non-negative IF and a model not composed of SHCs. Contrary to Table \ref{tab:3structures} rows two and three where SHCs are assumed, there can exist models with non-negative IF but without HC.  We illustrate the existence of such models with the following example. Consider the complex signal,
\begin{equation}
    z(t) = a_0 \cos(\omega_0 t) + j \alpha a_0 \sin(\omega_0 t)
\end{equation}
and the spectrum given by (\ref{eq:exFT3}) which is not one-sided. However in terms of a single AM--FM component, we can have
\begin{equation}
    z(t) = a_0(t) e^{j\theta_0(t)}
\end{equation}
where
\begin{equation}
    a_0(t) = \sqrt{a_0^2 \cos^2(\omega_0 t) + \alpha^2 a_0^2 \sin^2(\omega_0 t)  }
\end{equation}
and
\begin{equation}
    \theta_0(t) = \arctan\left[ \dfrac{\alpha a_0 \sin(\omega_0 t)}{a_0 \cos(\omega_0 t)}  \right].
\end{equation}
The IF, $\omega_0(t) = \dfrac{d}{dt}\theta_0(t)$ is strictly positive for any positive choice of $\alpha$. Thus the associated \emph{Hilbert} spectrum is one-sided even though the spectrum is two-sided. In other words, saying we assume positive IF is not the same as saying we assume a one-sided spectrum. Although signals can be designed such that $Z(\omega)$, for practical purposes, has a one-sided spectrum, e.g.~communications signals, this does not mean that naturally-occurring signals have, in general, a one-sided spectrum. Therefore, making the one-sided spectrum assumption (and assuming HC) may lead to incorrect model parameters.

Another incorrect assumption that is often made is that if the signal is narrowband, the HT will yield a meaningful complex extension. Unfortunately, this is not the case because narrowband signals exist that do not have Hermitian symmetry and hence use of the HT will yield incorrect results. This is demonstrated by (\ref{eq:exFT3}) which could be considered narrowband, yet use of the HT to extend its real observation would not yield the correct latent signal due to the absence of HC. 

By changing from SHCs to AM--FM components and relaxing the assumption of HC, we have generalized the spectrum to a Hilbert spectrum and effectively changed the definition of IF to be component-dependent. By reverting back to Carson's definition and not using SHCs in the analysis, we move away from Gabor and Ville's specialized definition back towards the generalized definition of IF. 

\begin{figure}[h]
    \centering
    \includegraphics[width = 0.5\linewidth]{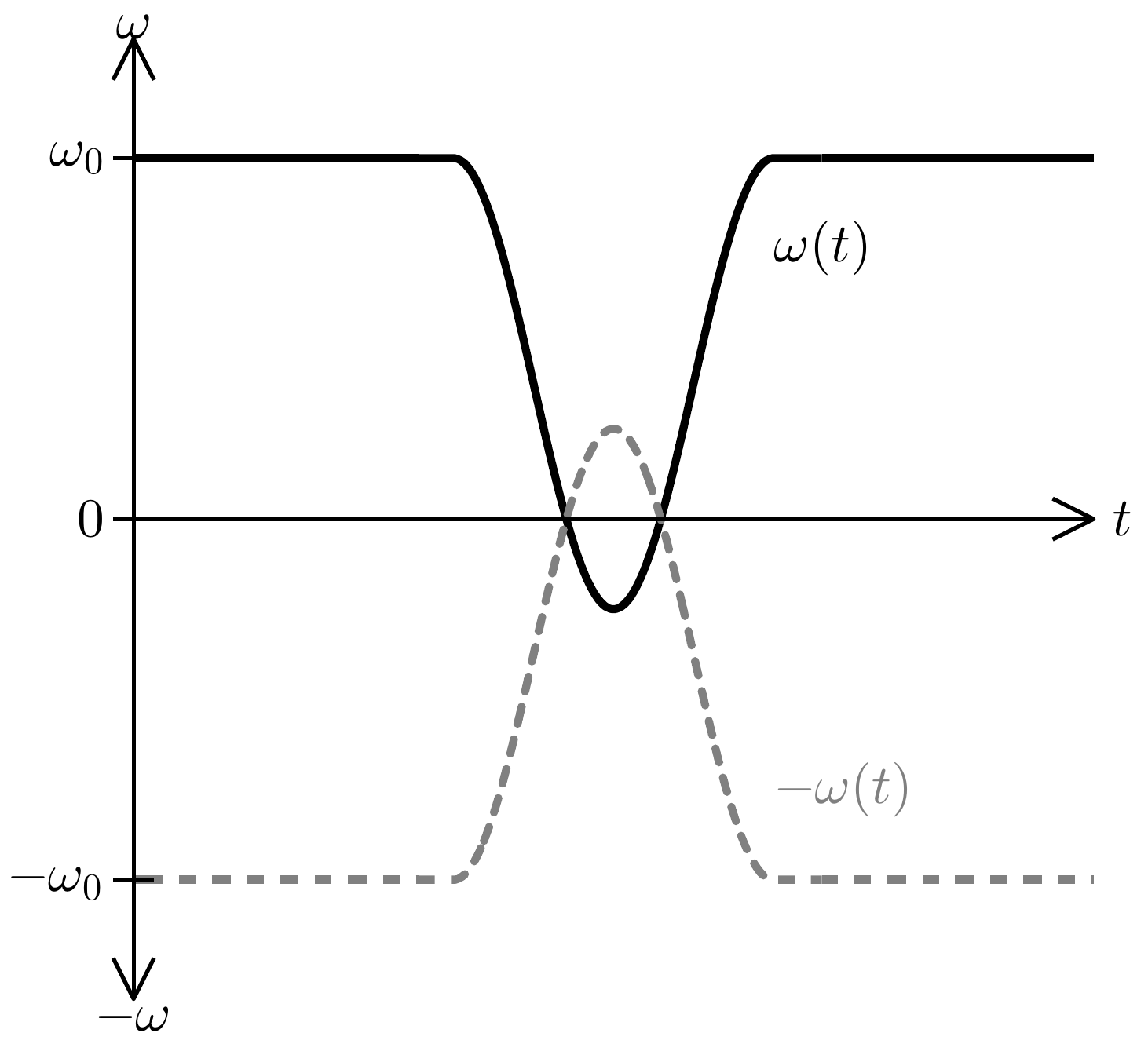}
    \caption{Illustration of when the assumption of non-negative IF cannot distinguish $z(t)$  with associated IF $\omega(t)$ (\solidSrule[3.5mm])  from $z^*(t)$ with associated IF $-\omega(t)$  (\textcolor{MyGrey} {\protect\dashedrule}) because each has both positive and negative values of IF at some time instances.}
    \label{fig:Illustration2}
\end{figure}

While using AM--FM components to define IF, which can change in time, we have to consider the possibility that a component's IF changes sign. Such a sign change is not possible under the assumption of SHCs which have constant IF. In this work, we have arbitrarily assumed that all components \emph{must have non-negative IF for all $t$}, although there may be some signal classes for which this is not true. In these classes, the AM--FM parameters may not match the underlying signal model parameters, although the superposition will yield the correct real signal. This is illustrated in Fig.~\ref{fig:Illustration2} where a component of $z(t)$ with associated IF $\omega(t)$ and the component of $z^*(t)$ with associated IF $-\omega(t)$ cannot be separated by assuming non-negative IF because each has both positive and negative instantaneous frequencies. Therefore in these cases, we need to relax the the assumption of non-negative IF at some time instances in order to properly estimate the latent signal.

\section{Examples using the AM--FM Model}
\label{sec:Solutions}

The example in Part I illustrated solutions to the LSA problem where we arrived at a \emph{single} pair of IA/IF parameters, $\rho(t)$ and $\Omega(t)$, for the latent \emph{signal}. In this section, we illustrate HSA consisting of a \emph{set} of IA/IF pairs, $\{a_k(t), \omega_k(t)\}$ for the signal, one pair for each latent \emph{component}. This is illustrated by the set diagram in Fig.~\ref{fig:NewSet}. In our closed-form solutions, we choose assumptions which lead to various decompositions and interpretations. For convenience, it is assumed that the phase reference $\phi_k = 0$, where possible.

\subsection{Periodic Triangle Waveform Example}
\label{ssec:trianglewaveexample}
As an example of HSA, we consider the triangle waveform from Part I (\ref{eq:trianglewaveP1}) and the sinusoidal FM signal.

\subsubsection{Analysis assuming simple harmonic components}
If we assume the component in (\ref{eq:AMFMcomp}) has the form given in~(\ref{eq:simpleHamonicComponent}), i.e. a SHC, then the AM--FM model can be obtained using Fourier analysis and the triangle waveform may be expressed with an infinite number of components
\begin{equation}
	x (t) = \Re\left\lbrace\sum\limits_{k=0}^{\infty} \frac{8A}{\pi^2}\frac{1}{(2k+1)^2} e^{j (2k+1) \omega_0 t}\right\rbrace
	\label{eq:TRI_FS}
\end{equation}
where the IA is given by
\begin{equation}
	a_k(t) = \frac{8A}{\pi^2}\frac{1}{(2k+1)^2},
	\label{eq:triangleshcia}
\end{equation}
the IF is given by
\begin{equation}
	\omega_k(t) = (2k+1) \omega_0,
	\label{eq:triangleshcif}
\end{equation}
and $\omega_0$ is a constant---conventionally known as the fundamental frequency.

\subsubsection{Analysis assuming a single AM--FM component with harmonic correspondence}
The solution in (\ref{eq:TRI_FS}) may be collapsed into a single, wideband AM--FM component ($K=1$) by rewriting as in Part I, in terms of the AM--FM model as
\begin{equation}
	x (t)	= \Re\left\lbrace a_0(t)e^{j[\omega_0 t+M_0(t)]} \right\rbrace
    \label{eq:trianglesingleamfm}  
\end{equation}
where the IA $a_0(t)$ is given by $\hat{\rho}(t)$ in (\ref{eq:trianglesingleamfmiaP1}) and the IF $\omega(t)$ is given by $\hat{\Omega}(t)$ in (\ref{eq:trianglesingleamfmifP1}).

\subsubsection{Analysis assuming a single AM component}
If we assume a single component ($K=1$) with constant frequency $\omega_0(t)=\omega_0$, the triangle waveform may be expressed as in Part I, in terms of the AM--FM model as
\begin{equation}
	x (t)=  \Re\left\lbrace a_0(t)e^{j\omega_0 t} \right\rbrace
	\label{eq:trianglesingleam}
\end{equation}
where the IA $a_0(t)$ is given by $\hat{\rho}(t)$ in (\ref{eq:trianglesingleamiaP1}).

\subsubsection{Analysis assuming a single FM component}
If we assume a single component ($K=1$) with constant amplitude $a_0(t)=A$, the triangle waveform may be expressed as in Part I, in terms of the AM--FM model as
\begin{equation}
	x ( t) =\Re\left\lbrace A e^{j \left[\omega_0 t+M_0(t)\right]} \right\rbrace
	\label{eq:trianglesinglefm}
\end{equation}
where the IF $\omega_0(t)$ is given by $\hat{\Omega}(t)$ in (\ref{eq:trianglesinglefmifP1}).

\subsection{Sinusoidal FM Example}
A sinusoidal FM signal with carrier frequency $\omega_c$, modulation frequency $\omega_m$, and constant $B$, is expressed as
\begin{equation}
	x (t) = \Re\left\lbrace e^{j\left[\omega_c t + B \sin(\omega_m t)\right]}\right\rbrace.
	\label{eq:sinFM}
\end{equation}

\subsubsection{Analysis assuming a single FM component}
We recognize that (\ref{eq:sinFM}) is already in the form of the AM--FM model, if we assume a single component ($K=1$) with IA 
\begin{equation}
    a_0(t) = 1.
    \label{eq:sinFMfmia}
\end{equation}
This leads to the IF
\begin{equation}
    \omega_0(t) = \omega_{c} + \frac{d}{dt}B\sin(\omega_m t).
    \label{eq:sinFMfmif}
\end{equation}

\subsubsection{Analysis assuming simple harmonic components}
Alternatively, the sinusoidal FM signal can be expressed in terms of the AM--FM model as 
\begin{equation}
	x (t) = \Re\left\lbrace \sum\limits_{k=-\infty}^{\infty}J_k(2\pi B/\omega_m)e^{j[(\omega_c+k\omega_m)t+\phi_k]} \right\rbrace
\end{equation}
where $J_k(\cdot)$ denotes the $k$th-order Bessel function of the first kind \cite{titchmarsh}. This yields an infinite number of components with IA given by
\begin{equation}
	a_k(t)=J_k(2\pi B/\omega_m)
	\label{eq:sinFMfsia}
\end{equation}
and IF given by
\begin{equation}
	\omega_k(t) =\omega_c+k\omega_m.
	\label{eq:sinFMfsif}
\end{equation}
This is an example of a signal class where the assumption of non-negative IF components in (\ref{eq:NonNegInstFreq}) must be relaxed, as a consequence of the summation on $k$ becoming double-sided. For most practical choices of $B$ and $\omega_m$, the $J_k(2\pi B/\omega_m)$ associated with the negative frequencies will be vanishingly small, therefore the assumption of non-negative IF, is for practical purposes, true.

\subsection{Remarks} 
As demonstrated above, we obtain different parameterizations of the same signal by changing assumptions made during analysis. Some parameterizations, such as those obtained using Fourier analysis and the HT, are closely related due to similar assumptions. However, there exist many parameterizations each with different interpretations and without other information, we cannot say that one particular parameterization is any better or worse than any other. This highlights the importance of considering the hidden assumptions when using any particular method (FT, HT, STFT, HHT, etc.) and interpreting meaning from the parameters. If assumptions in analysis do not match the underlying signal model, erroneous interpretations may be made. 

The previous solutions make use of different quadratures each corresponding to different signal models. This demonstrates the non-uniqueness of the quadrature for the AM--FM component, and hence our view of the quadrature as a free parameter. Also, note that the solutions in (\ref{eq:TRI_FS}), (\ref{eq:trianglesingleamfm}), (\ref{eq:trianglesingleam}), etc., each utilize a complex extension that is implied by the particular assumptions made, rather than utilizing an extension based on a single rigorously-defined procedure, like the HT. In fact, for any practically-chosen real signals $x(t)$ and $y(t)$, there exist assumptions leading to instantaneous parameters in which $y(t)$ can be viewed as the quadrature of $x(t)$.

The AM--FM model leads to exact solutions for $a(t)$ and $\omega(t)$, so it might appear to violate the uncertainty principle, i.e. exceed the lower limit of the time-bandwidth product. However, when AM--FM modeling is viewed as a quantum mechanics problem, our casting of the problem is fundamentally different than Gabor's. A comparison of the formal correspondences between quantum mechanics and short-time Fourier analysis and quantum mechanics and Hilbert spectral analysis is given in Table~\ref{tab:Quantum}. In Gabor's casting, \emph{time and frequency} are the complementary variables, i.e.~(position, momentum)$\leftrightarrow$($t$, $\omega$) \cite{cohen1995time, loughlin2004uncertainty} while for the AM--FM model, the \emph{real and quadrature signals} are the complementary variables, i.e.~(position, momentum)$\leftrightarrow$($x(t)$, $y(t)$). Therefore, all uncertainty arises because only the real signal is observed. We believe that our casting is not only more appropriate, but also more useful and powerful.

\begin{table}[h]
\caption{The formal correspondences between quantum mechanics and short-time Fourier analysis (p.~197 in \cite{cohen1995time}) (columns 2 and 1)  and quantum mechanics and Hilbert spectral analysis (columns 2 and 3).}
	\scriptsize
	\centerline{\begin{tabular}{llcllcll}
			\toprule[1.5pt]
			\multicolumn{2}{c}{\head{Fourier Analysis}}& $~$ &\multicolumn{2}{c}{\head{Quantum Mechanics}}&  $~$ & \multicolumn{2}{c}{\head{\textbf{Hilbert Spectral Analysis}}}  \\
			\cmidrule{1-2}                                					\cmidrule{4-5}             	\cmidrule{7-8}                          
			\textbf{Description}   & \textbf{Symbol}    &   $~$ & \textbf{Description} &  \textbf{Symbol}   & $~$ & \textbf{Description}   &   \textbf{Symbol}                                    \\
			\midrule
			Time        &  $t$                  &$\leftrightarrow$ & Position      & $q$                                     & $\leftrightarrow$        & Real Signal   & $x(t)$                                   \\
			Frequency   &  $\omega$             &$\leftrightarrow$ & Momentum      & $p$                                     & $\leftrightarrow$        & Imaginary Signal   & $y(t)$                                   \\
			{\scriptsize\textcolor{MyGrey}{no correspondence}}         &  $~$                  &$\leftrightarrow$ & Time          & $t$                                     & $\leftrightarrow$        & Time          & $t$                                      \\
			Signal      &  $x(t)$               &$\leftrightarrow$ & Wave Function & $\Psi(q,t)$                               & $\leftrightarrow$        & Latent Signal     & $z(t)$                                \\
			Uncertainty &  $BT\geq\frac{1}{2}$  &$\leftrightarrow$ & Uncertainty   & $\sigma_p\sigma_q\geq\frac{1}{2}\hbar$  & $\leftrightarrow$        & Uncertainty   & $x(t) = \Re\{z(t)\}$                  \\
			$~$         & $~$                   &$~$ & $~$           & $~$                                     & $~$        & Instantaneous Frequency  & $\omega_k(t)=\frac{d}{dt}\arg\{\psi_k(t)\}$  \\
			\bottomrule[1.5pt]
		\end{tabular}}
		\label{tab:Quantum}
\end{table}

\section{Proposed Visualization of the Hilbert Spectrum}
\label{sec:visual}
The ability to visualize and interpret model parameters is key to the adoption of any analysis method. Often complex AM--FM signals are plotted as a series of real 2-D plots, i.e.~$s_k(t)$ vs.~$t$, which for AM--FM components, would provide little insight into the underlying signal model. Alternatively, Argand diagrams may be used for AM--FM signal visualization, however, drawbacks include the quadrature signal possibly having no assigned meaning and the revolution rate not intuitively displayed.

A more appropriate visualization plots the Hilbert spectrum in its entirety. Unfortunately, plots of the Hilbert spectrum are often crudely discretized because a clear distinction between instantaneous parameters and spectral parameters is not made \cite{huang1998empirical}. We propose a method for visualizing the Hilbert spectrum, which is both complete, intuitive, and avoids the coarse discretization. The proposed visualization can be considered a (pseudo-) phase-space plot because every degree of freedom or parameter of each AM--FM component is represented. 

By plotting $\omega_k(t)$ vs.~$s_k(t)$ vs.~$t$ as a line in a 3-D space and coloring the line with respect to $|a_k(t)|$ for each component, the simultaneous visualization of multiple parameters for each component is possible. Further, orthographic projections yield common plots: the time-real plane (the real signal waveforms), the time-frequency plane (Hilbert spectrum), and the real-frequency plane (analogous to the Fourier magnitude spectrum). We have found it beneficial to interpret each component as an ``illuminated'' particle moving in 3-D space. Each particle's motion in time and frequency is governed by $\psi_k(t)$. The proposed method allows one to visualize the assumed underlying signal model. 

To illustrate visualization of the Hilbert spectrum, we plot the examples in Section \ref{sec:Solutions}. However, the sophisticated nature of the proposed 3-D visualization of the Hilbert spectrum is not well-accommodated with paper media, and as a result we provide associated {\sc matlab} functions and additional visualizations online \cite{HSAlink}. It is our preference that Hilbert spectra not be visualized with paper media. In the plots, we have utilized a perceptually-motivated colormap in order to improve interpretation over other colormaps \cite{borland2007rainbow, NiccoliPMK}.

\begin{figure}[ht]
\centering
  	\begin{minipage}[b]{0.476\linewidth}
  		\centering
  		\subfigure[]{
  		\includegraphics[width = 0.99\linewidth]{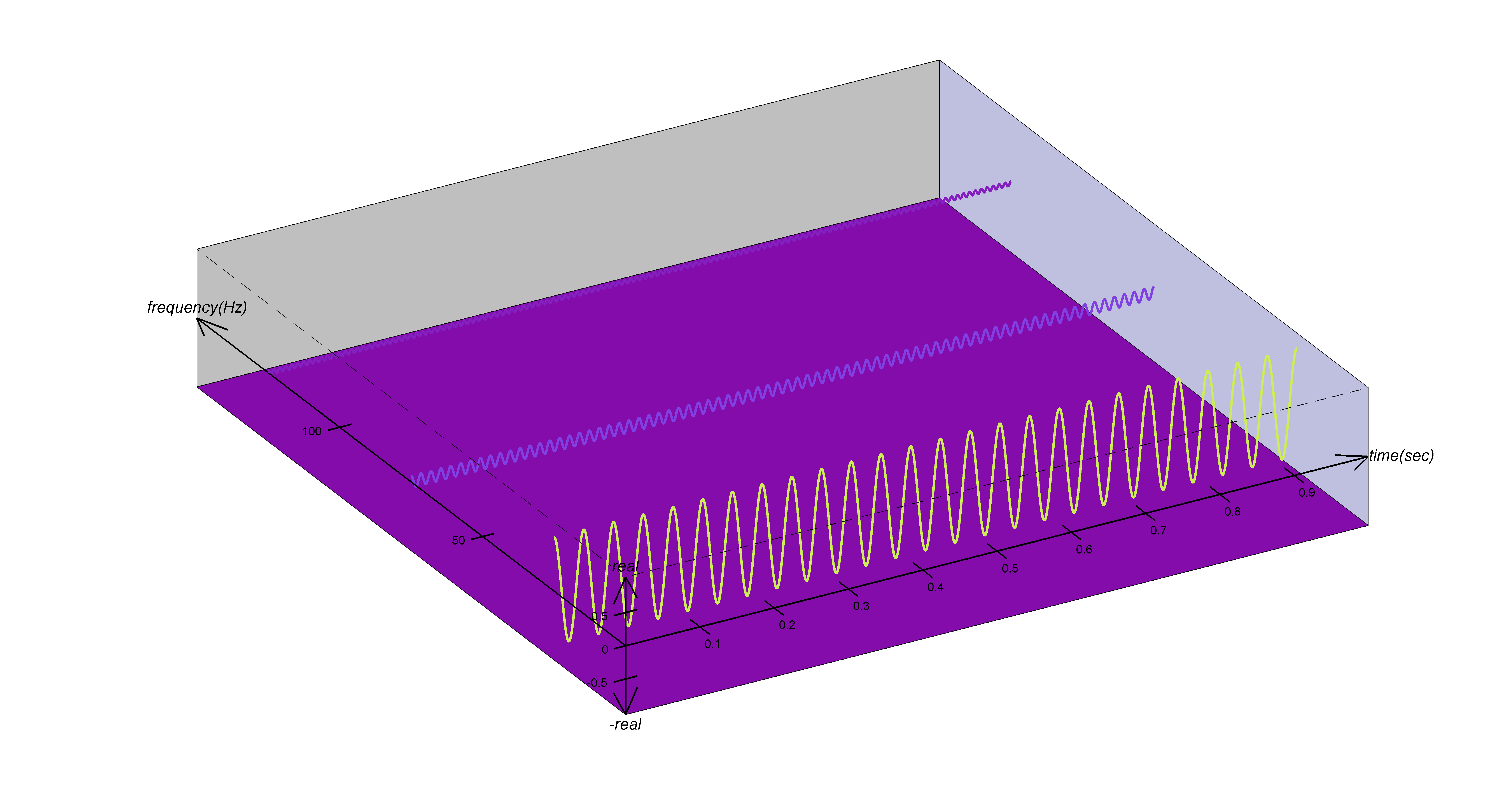}
  		\label{fig:TRI_FS3d}
  	}
  	\end{minipage}
  	\begin{minipage}[b]{0.476\linewidth}
  		\centering
  		\subfigure[]{
  		\includegraphics[width = 0.99\linewidth]{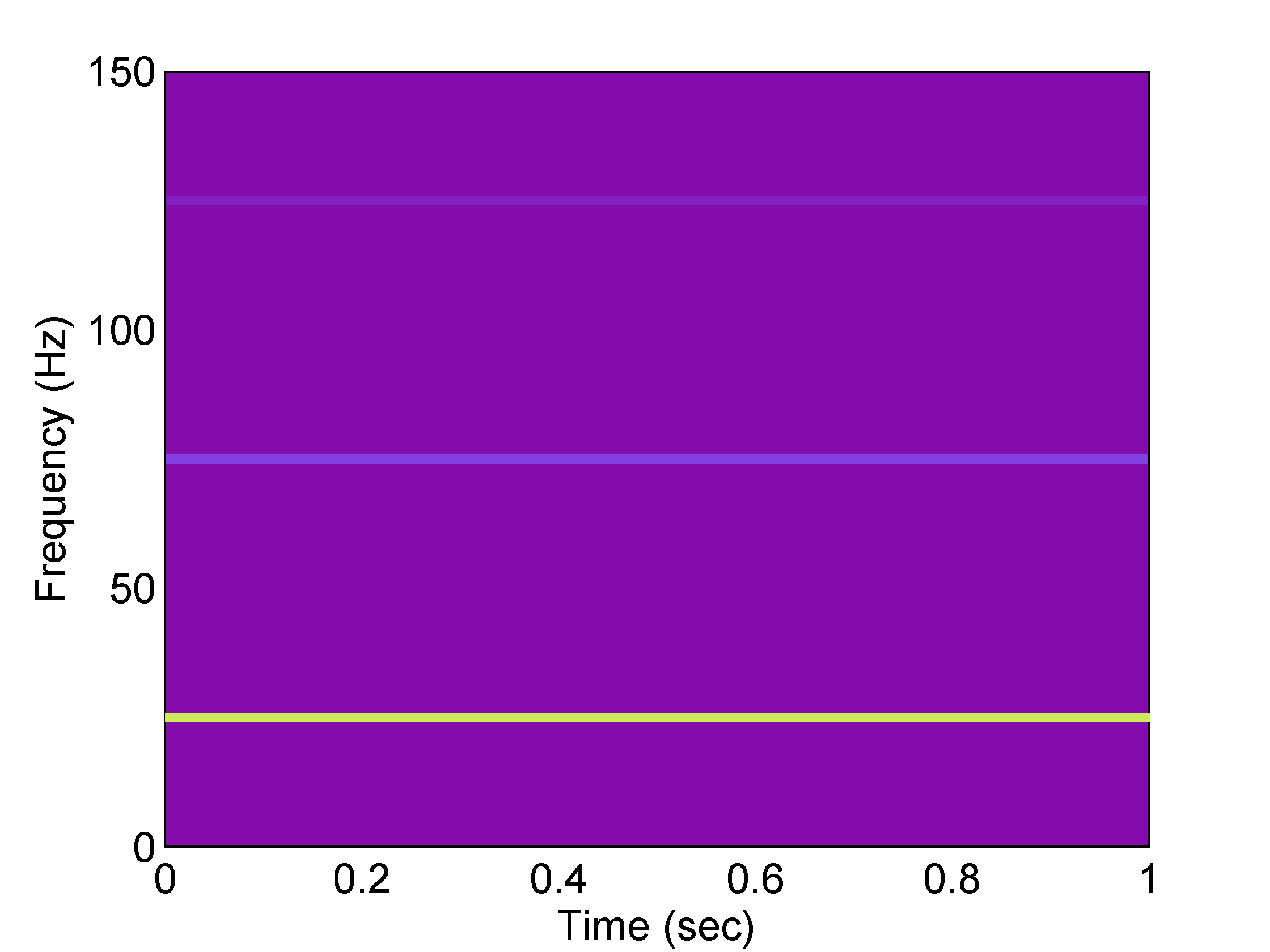}
  		\label{fig:TRI_FS}
  	}
  	\end{minipage}\\
  	\begin{minipage}[b]{0.476\linewidth}
  		\centering	
  		\subfigure[]{
  		\includegraphics[width = 0.99\linewidth]{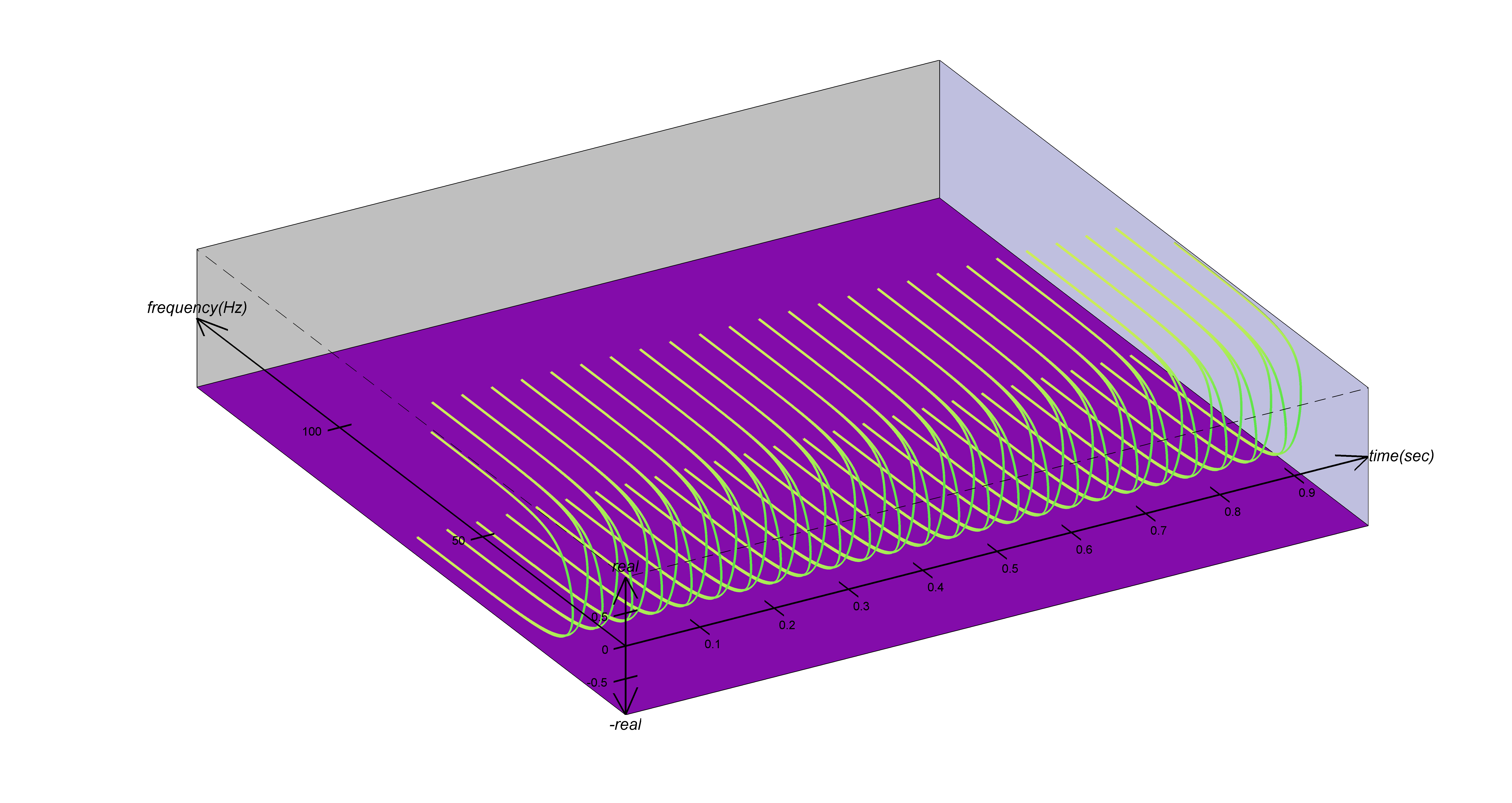}
  		\label{fig:TRI_AMFM3d}
  	}
  	\end{minipage}
  	\begin{minipage}[b]{0.476\linewidth}
  		\centering	
  		\subfigure[]{
  		\includegraphics[width = 0.99\linewidth]{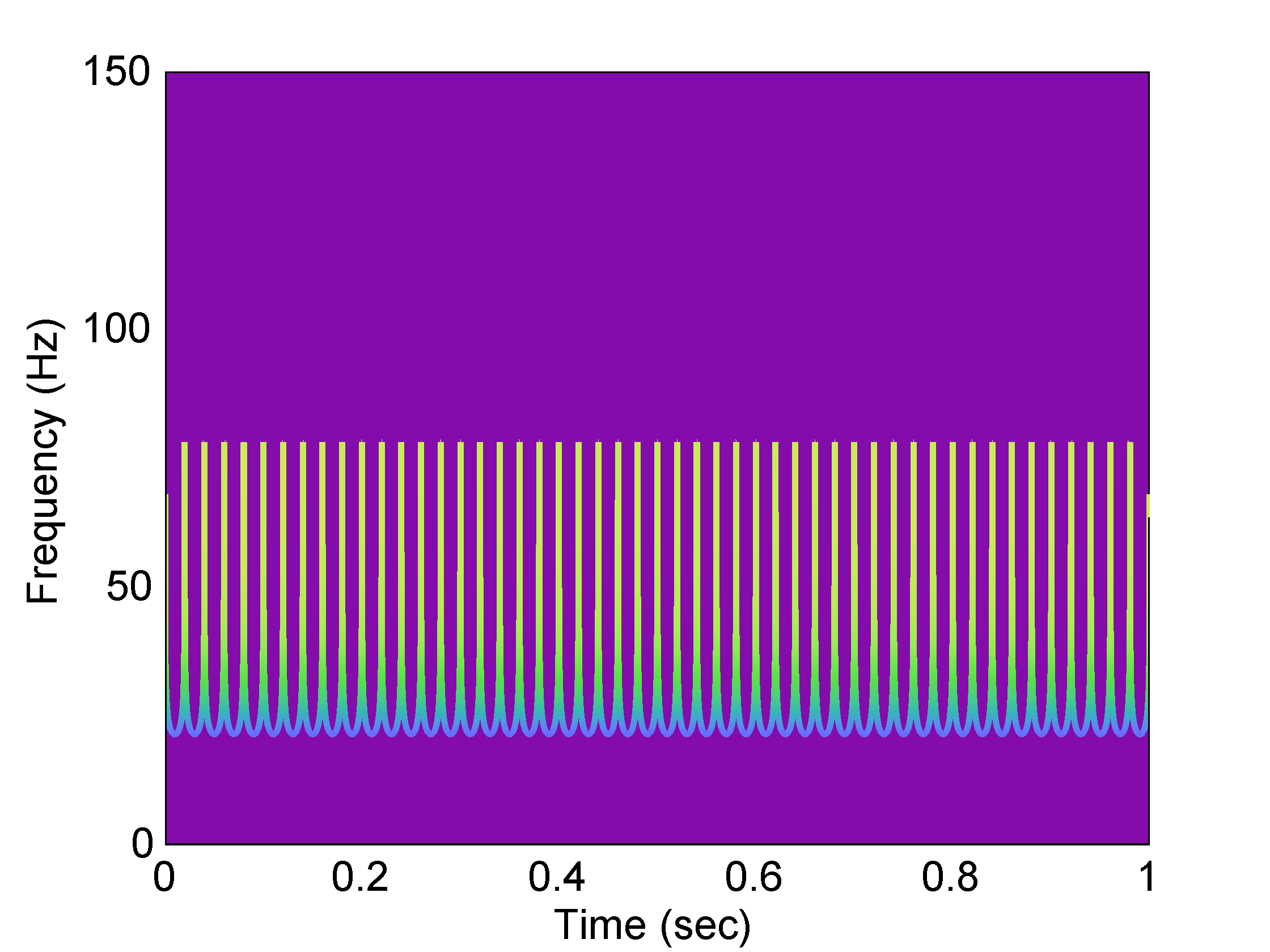}
  		\label{fig:TRI_AMFM}
  	}
  	\end{minipage}\\
  	\begin{minipage}[b]{0.476\linewidth}
  	 		\centering	
  	 		\subfigure[]{
  	 		\includegraphics[width = 0.99\linewidth]{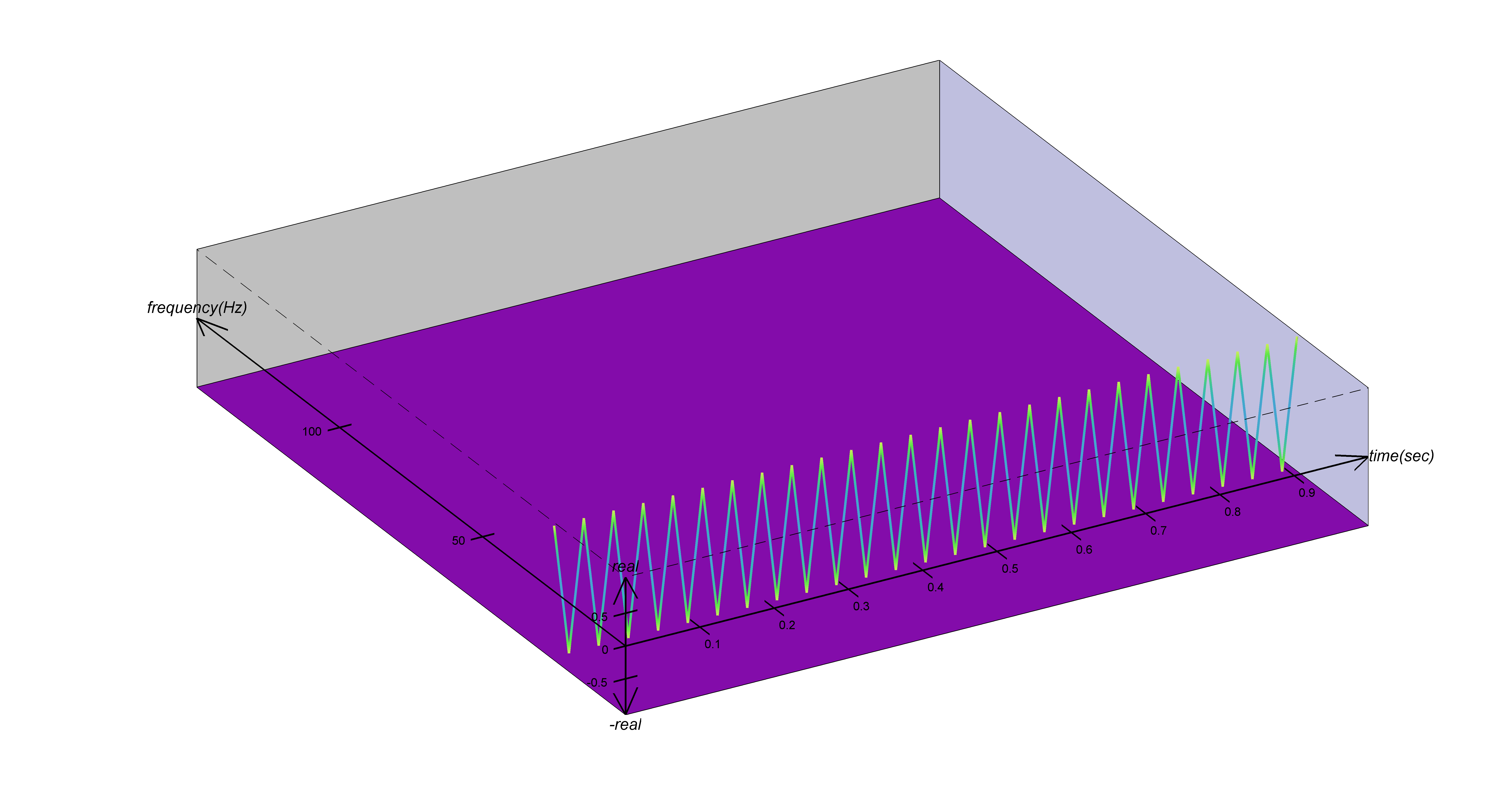}
  	 		\label{fig:TRI_AM3d}
  	 	}
  	\end{minipage}
  	\begin{minipage}[b]{0.476\linewidth}
  	 		\centering	
  	 		\subfigure[]{
  	 		\includegraphics[width = 0.99\linewidth]{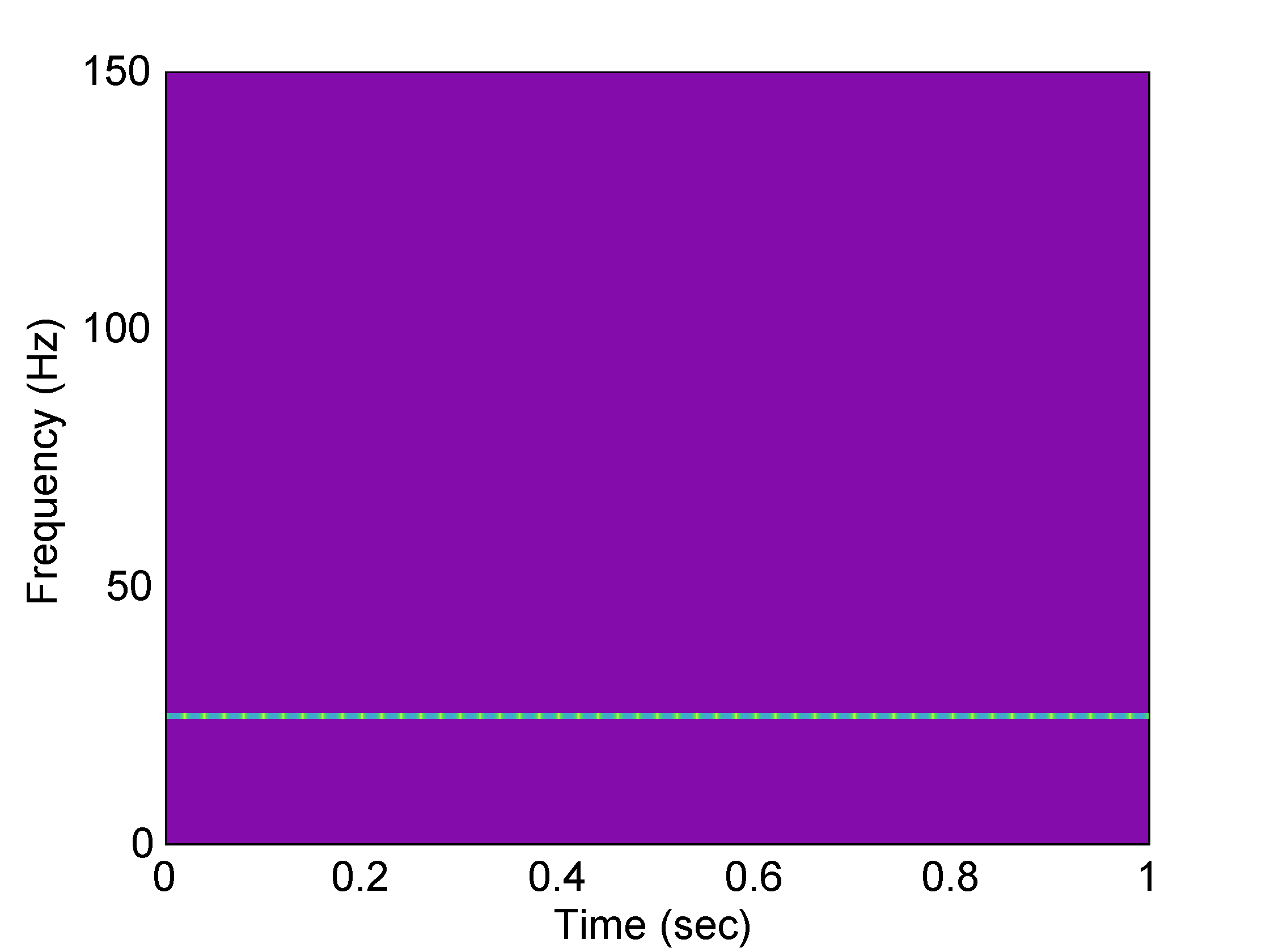}
  	 		\label{fig:TRI_AM}
  	 	}
  	\end{minipage}	
        \caption{Hilbert spectrum for the triangular waveform $x (t)$ with $\omega_0 = 50\pi$ rads/s for the assumptions of (a) SHCs, (c) single AM--FM component with HC, and (e) single AM component. The corresponding time-frequency planes, obtained by projecting out the $s_k(t)$ dimension, are shown in (b),(d),(f). }
  		\label{fig:TRIexample}
\end{figure}

For the triangle waveform example, Figs.~\ref{fig:TRIexample}(a),(c),(e) illustrate the Hilbert spectrum plots for the assumptions of SHC, single AM--FM component with HC, and single AM component, respectively. Fig.~\ref{fig:TRIexample}(a) shows the first three SHCs [constant amplitude and constant frequency in (\ref{eq:triangleshcia})] at integer multiples of a fundamental frequency, $50\pi$ rads/s. Fig.~\ref{fig:TRIexample}(c) shows the single AM--FM component with HC, where we see line color variation indicating a time-varying IA in (\ref{eq:trianglesingleamfmiaP1}) and a clear time-varying IF in (\ref{eq:trianglesingleamfmifP1}). Fig.~\ref{fig:TRIexample}(e) shows the single AM component, where we see a constant IF and color variation indicating the time-varying IA in (\ref{eq:trianglesingleamiaP1}). Figs.~\ref{fig:TRIexample}(b),(d),(f) show the corresponding time-frequency planes of Figs.~\ref{fig:TRIexample}(a),(c),(e) by projecting out the $s_k(t)$ dimension.

\begin{figure}[ht]
\centering
  	\begin{minipage}[b]{0.49\linewidth}
  		\centering
  		\subfigure[]{
  		\includegraphics[width = 0.99\linewidth]{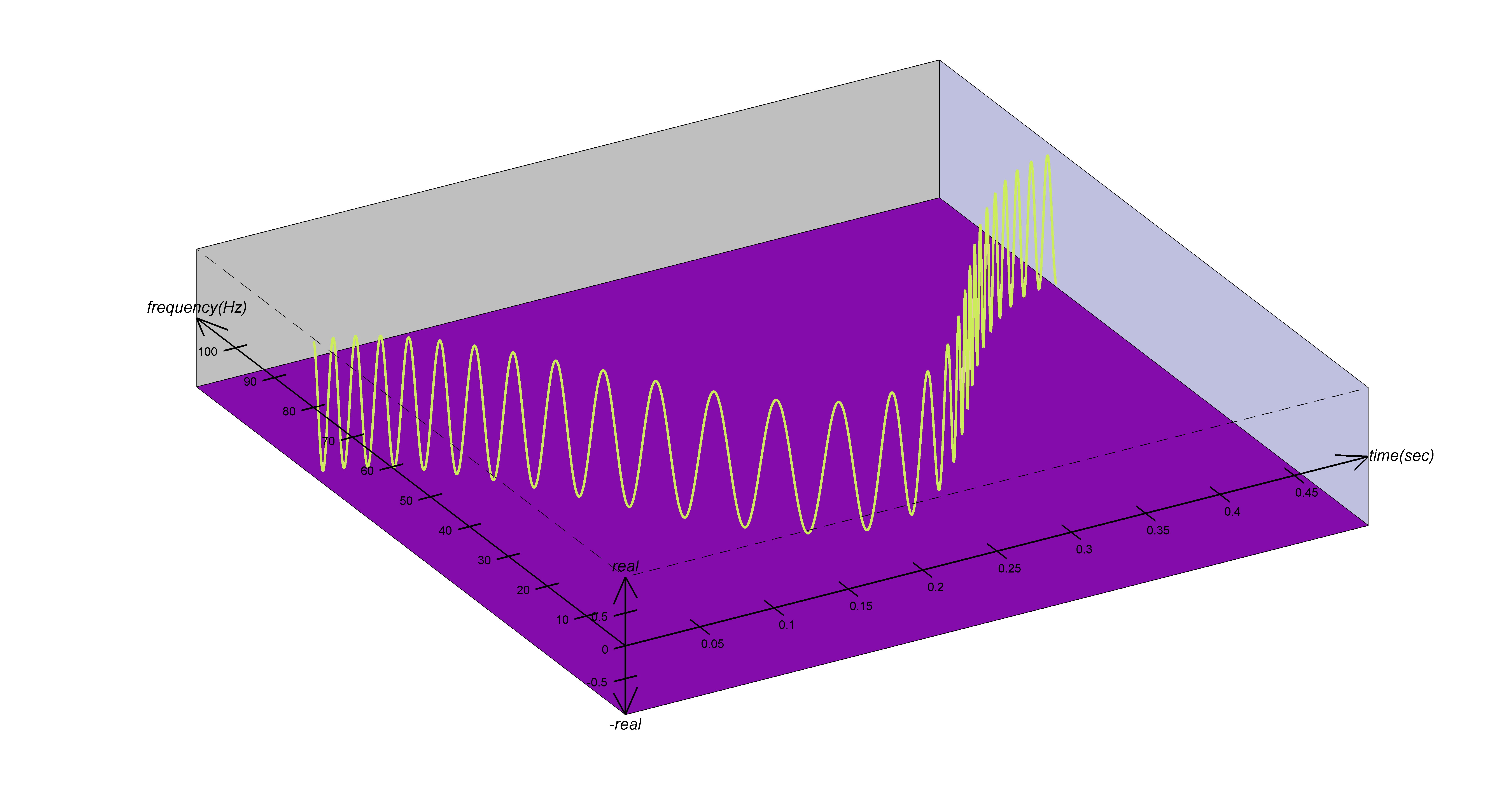}
  		\label{fig:SFM_FM3d}
  	}
  	\end{minipage}
  	  	\begin{minipage}[b]{0.49\linewidth}
  		\centering
  		\subfigure[]{
  		\includegraphics[width = 0.99\linewidth]{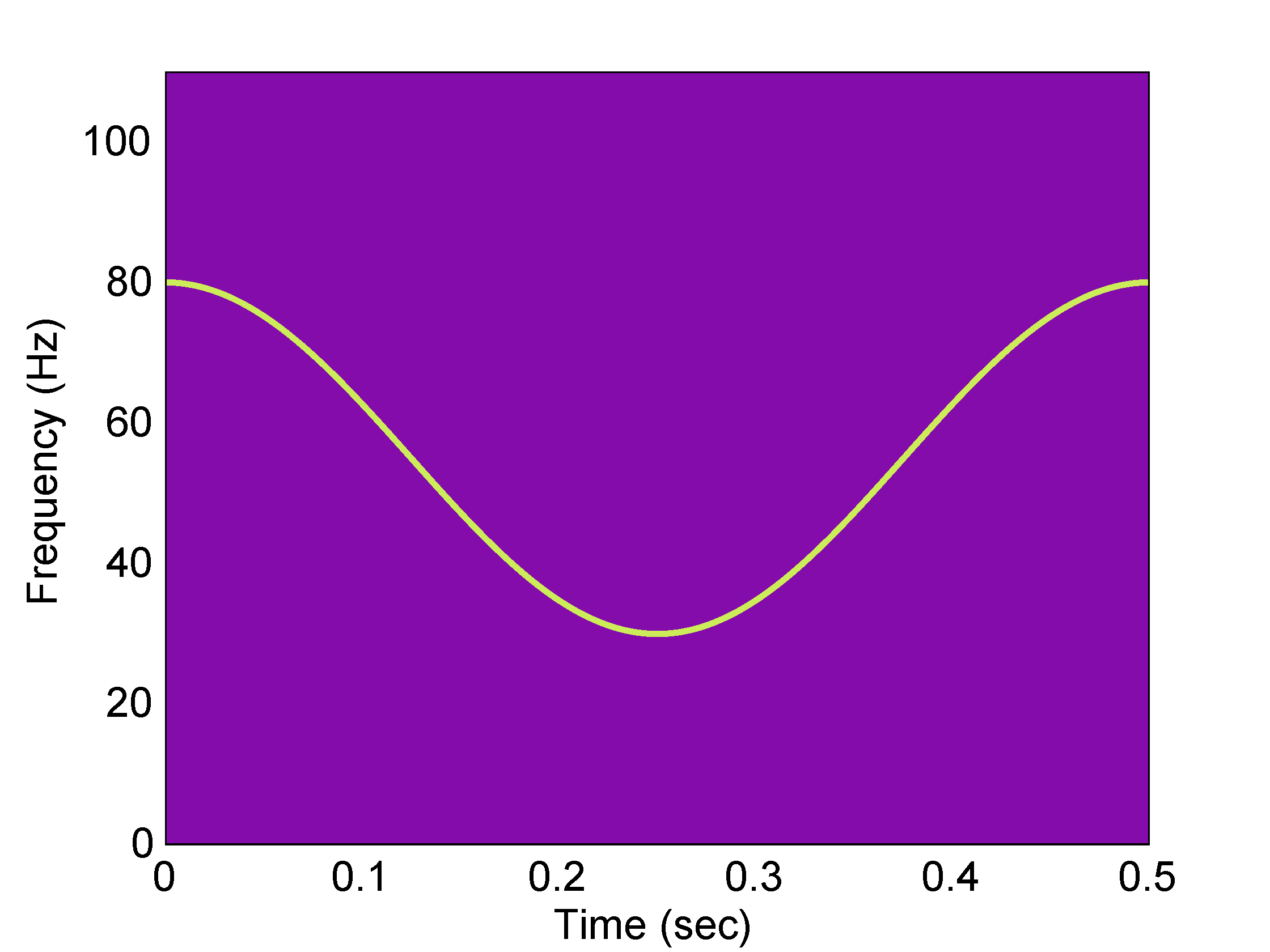}
  		\label{fig:SFM_FM}
  	}
  	\end{minipage}\\
  	\begin{minipage}[b]{0.49\linewidth}
  		\centering	
  		\subfigure[]{
  		\includegraphics[width = 0.99\linewidth]{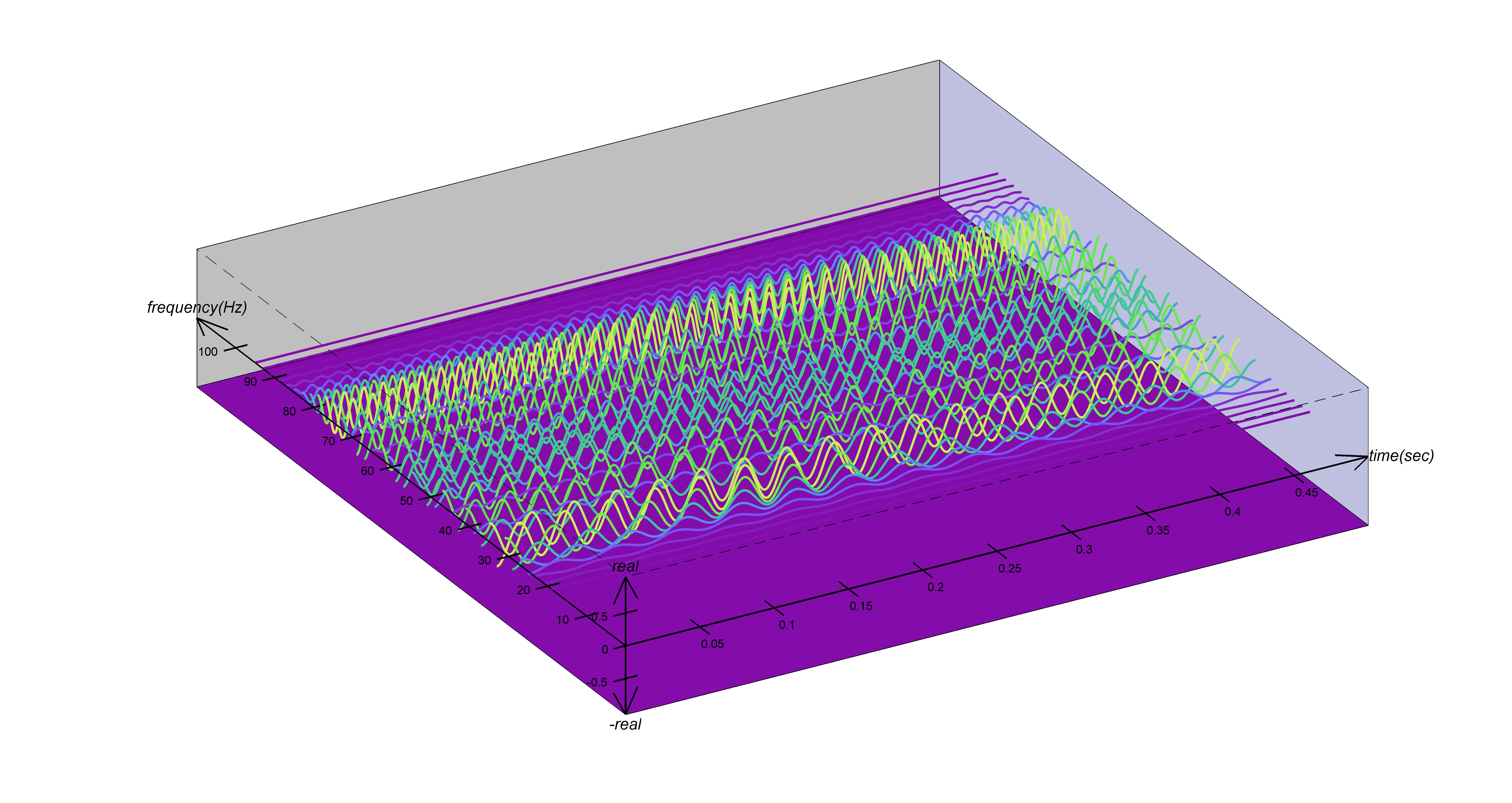}
  		\label{fig:SFM_FT3d}
  	}
  	\end{minipage}
  	\begin{minipage}[b]{0.49\linewidth}
  		\centering	
  		\subfigure[]{
  		\includegraphics[width = 0.99\linewidth]{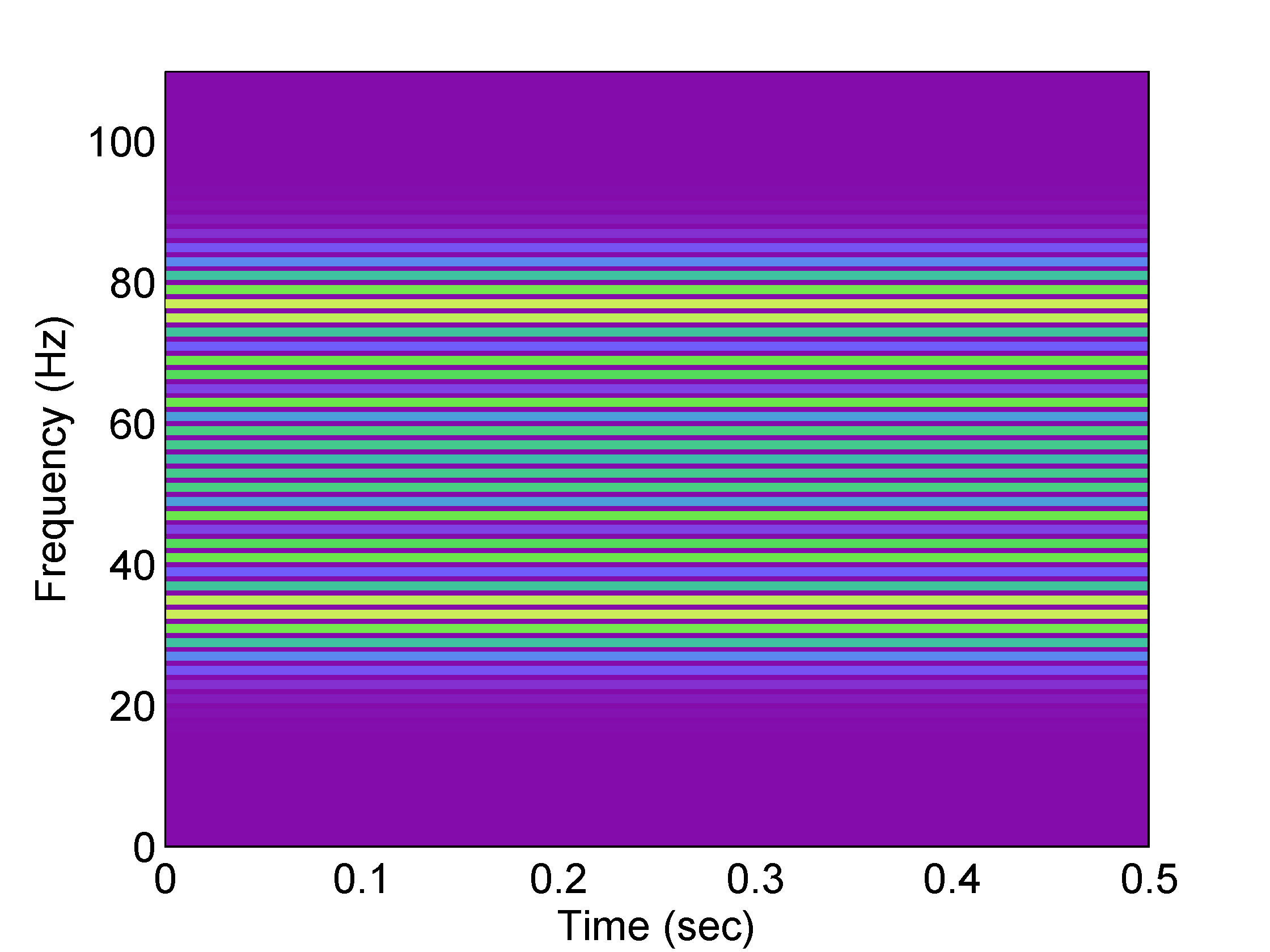}
  		\label{fig:SFM_FT}
  	}
  	\end{minipage}
  	\caption{Hilbert spectrum for the sinusoidal FM waveform $x (t)$ with $\omega_c = 110\pi$ rads/s, $\omega_m = 4\pi$ rads/s, and $B = 25$ for the assumptions of (a) a single FM component and (c) SHCs. The corresponding time-frequency planes, obtained by projecting out the $s_k(t)$ dimension, are shown in (b),(d).}
  	  \label{fig:SFM}	
\end{figure}

For the sinusoidal FM example, Figs.~\ref{fig:SFM}(a) and (c) illustrate the Hilbert spectrum plots for the assumptions of a single FM component and SHCs respectively. Fig.~\ref{fig:SFM}(a) shows the constant IA in (\ref{eq:sinFMfmia}) indicated by a constant line color and time varying IF in (\ref{eq:sinFMfmif}). Fig.~\ref{fig:SFM}(c) shows SHCs at integer multiples of a fundamental frequency $4\pi$ rads/s where each component has constant IA in (\ref{eq:sinFMfsia}). Figs.~\ref{fig:SFM}(b) and (d) show the corresponding time-frequency planes of Figs.~\ref{fig:SFM}(a) and (c) by projecting out the $s_k(t)$ dimension.

\section{Summary} \label{sec:SummaryP2}
In Part II of this paper, we have presented theory for the Hilbert spectrum as a generalized LSA problem. In the general problem, we seek a representation of $z(t)$ consisting of a superposition of latent signals, i.e.~multicomponent model consisting of a superposition of complex AM--FM components. We have used the AM--FM model to define the Hilbert spectrum as parameterized by a set of IA/IF pairs (and phase references) each associated with the components. In the LSA problem, relaxation of the HC condition allowed freedom in the choice of the signal's quadrature and admitted many solutions for the latent signal. In the HSA problem, relaxation of the HC condition allows freedom in the choice of each component's quadrature thereby allowing even greater freedom in the construction of the signal's model. We illustrated how assumptions on the form of the component, i.e.~constant amplitude and constant frequency, time-varying amplitude and constant frequency, constant amplitude and time-varying frequency, etc., lead to familiar specializations of the model. 

We discussed the frequency domain view of LSA including the implications of relaxing HC. From there we discussed the Hilbert spectrum view of LSA, specifically the implications of using a general AM--FM component with relaxed HC and assumed positive IF. Closed form HSA examples were provided in which we assume various forms of the AM-FM component and determined the unique corresponding instantaneous parameterization in terms of IA and IF.  A novel 3D visualization of the Hilbert spectrum was proposed by plotting $\omega(t)$ vs.~$s(t)$ vs.~$t$ and coloring with respect to $|a(t)|$. 

Finally, we discussed the analogy of the Hilbert spectrum to quantum mechanics in which our casting of time-frequency analysis is fundamentally different than Gabor's casting, where the uncertainty in our framework is in the quadrature signal and not the frequency variable. This provides a new and powerful framework for nonstationary signal analysis.

\FloatBarrier

\hypertarget{Part3}{}
\begin{center}\textbf{Part III: Numerical HSA Assuming Intrinsic Mode Functions}\end{center}
\setcounter{section}{0}
\bookmark[level=part,dest=Part3]{Part III: Numerical HSA Assuming Intrinsic Mode Functions}

\begin{quote}
[T]here are some crucial restrictions of the Fourier spectral analysis; the system must be linear; and the data must be strictly periodic or stationary; otherwise, the resulting spectrum will make little physical sense.---Huang \cite{huang1998empirical}
\end{quote}

In Part II, the AM--FM signal model and HSA theory was presented, analysis was performed on example signals, and a visualization of the Hilbert spectrum was proposed. In Part III, numerical algorithms are given for estimating the IAs and IFs of the components in the AM--FM model where now the component is assumed to be an Intrinsic Mode Function (IMF). These algorithms first decompose the signal into IMFs using an improved version of Huang's original Empirical Mode Decomposition (EMD) algorithm and second, demodulate the IMFs to obtain the instantaneous parameters. Unlike previous studies, close attention is paid to the assumptions made in the definition of the IMF which are carried forward to the demodulation step, thereby avoiding any ambiguity associated with obtaining the instantaneous parameters. We begin with a comprehensive review of EMD and several variations, and propose an algorithm for the computation of the Hilbert spectrum assuming IMF components. It is important to note that while IMFs can be considered latent AM--FM components, there are other classes of AM--FM components that are not IMFs as illustrated in Fig.~\ref{fig:Classes}. Examples using the proposed algorithm are provided that highlight alternative decompositions compared to traditional Fourier analysis and demonstrates the advantages of using the HSA framework.

\section{Introduction}
\label{sec:introP3}
In Part II of this paper, we defined the Hilbert spectrum using the AM--FM model in (\ref{eq:AMFMmodel}) where the AM--FM component is defined in (\ref{eq:AMFMcomp}) and parameterized by the IA $a_k(t)$, IF $\omega_k(t)$, and phase reference $\phi_k$.  We assume a real observation, $x(t)$ of a latent signal, $z(t)$.

In Part III of this paper, we turn our attention to numerical computation of instantaneous parameters of the AM--FM model. Many methods for computing the parameters of an AM--FM model already exist, each corresponding to a specific set of assumptions. As discussed in Part II, a FT corresponds to the assumption of SHCs and a STFT corresponds to the assumption of AM components.  Other AM--FM models have been proposed with alternative assumptions, however, common to these are restrictions which limit the utility of the model.

Practical estimation of the instantaneous parameters of the AM--FM model in (\ref{eq:AMFMmodel}) is a two-step process.  First, the signal must be decomposed into a set of latent AM--FM components and second, the instantaneous parameters $\{a_k(t), \omega_k(t)\}$ of each component must be estimated. To date, the most flexible decomposition method for AM--FM modeling is Huang's EMD and its variations \cite{huang1998empirical}. In contrast to most time-frequency analysis methods, EMD makes no assumption of constant amplitude or frequency on the component and has mild bandlimiting assumptions. In \cite{huang1998empirical}, Huang proposed the original EMD algorithm which sequentially determines a set of AM--FM components, i.e.~IMFs via an iterative sifting algorithm. The Ensemble Empirical Mode Decomposition (EEMD) \cite{wu2009ensemble} and tone masking \cite{DeeringMasking} introduced ensemble averaging in order to address the mode mixing problem. The complete EEMD was proposed to address some of the undesirable features of EEMD by averaging at the component-level as each component is estimated rather than averaging at the conclusion of EMD \cite{torres2011complete}. Several improvements to the sifting algorithm have also been proposed including those by Rato \cite{rato2008hht}. These improvements to the original EMD will be more fully described in this part.

In the second step, instantaneous parameters must be estimated through demodulation of the IMFs. Despite the numerous attempts to demodulate IMFs, most of the proposed methods have fallen short because the assumptions made during decomposition have not been maintained during demodulation. In this paper, we identify and pay close attention to these assumptions as we develop a demodulation technique which adheres to the assumptions. We point out that Rato proposed an AM--FM demodulation procedure in which the IA estimation was consistent with the assumptions but the IF estimation was not \cite{rato2008hht}. Also, Huang has examined numerous demodulation methods, including an iterative normalization to obtain an FM signal which is then demodulated to estimate IF using an $\arctan$ approach that unfortunately suffers from numerical instability \cite{huang2009instantaneous}.  The iterative normalization method is consistent with the decomposition assumptions. Utilizing both Rato's AM estimation and Huang's iterative normalization procedure, we propose a mathematically equivalent method to get IF from the FM signal that is more numerically stable. We then incorporate the proposed demodulation and EMD into a single HSA--IMF algorithm which gives very good estimates for the IA and IF parameters of the AM--FM model.

Part III of this paper is organized as follows. In Section \ref{sec:EMD}, we review the EMD algorithm and its variations and improvements by other researchers. By carefully noting the assumptions made in EMD, we propose demodulation of the resulting IMFs and present a complete HSA--IMF algorithm composed of EMD and demodulation.  In Section \ref{sec:Examples}, we provide two sets of examples using the HSA--IMF algorithm.  The first set of examples uses synthetic signals which illustrate how EMD can estimate the parameters of the underlying, assumed AM--FM signal model. The second set of examples uses real-world audio signals including speech which further illustrates AM--FM model parameterizations with alternate interpretations.  These examples clearly demonstrate the alternate decomposition and interpretation offered by HSA. In Section \ref{sec:discussion}, we discuss other issues related to HSA and in Section \ref{sec:further}, we discuss future research. In Section \ref{sec:SummaryP3} we summarize Part III.

\section{Empirical Mode Decomposition} \label{sec:EMD}
EMD consists of an iterative procedure for decomposing a signal into a set of IMFs, $\{\varphi_k(t)\}$ \cite{huang1998empirical}. In \cite{huang1998empirical}, the definition of an IMF is any signal that satisfies two conditions:
\begin{itemize}[labelwidth = \widthof{~~~~~C1:},leftmargin=!]
    \item[\textbf{C1:}] In the whole signal segment, the number of extrema and the number of zero crossings must be either equal or differ at most by one.
    \item[\textbf{C2:}] At any point the mean value of the envelope, defined by the local maxima and the envelope defined by the local minima, is zero.
\end{itemize}
In the context of LSA and AM--FM modeling, an IMF can be thought of as a latent \emph{component},
\begin{equation}
    \varphi_k(t) = \Re \{ \psi_k(t) \}.
    \label{eq:realIMF}
\end{equation}
Thus we view an IMF as an inherently complex-valued component, which is not the conventional interpretation. Furthermore, we argue that the definition of an IMF forces a unique quadrature that in general does not equal $\mathcal{H}\{\varphi_k(t)\}$ as will be discussed in Section \ref{ssec:demod}. 

Unlike traditional methods in time-frequency analysis, EMD is defined by an algorithm rather than transform theory, which has both advantages and disadvantages. One advantage is that EMD utilizes a component that is far less restrictive than other time-frequency methods \cite{cohen1995time}.  One disadvantage is that EMD is understood primarily through experimentation \cite{DeeringMasking}. Empirical experiments using white noise have shown EMD to act as a dyadic filter bank \cite{wu2004study,flandrin2004empirical,FlandrinHuang2005, wu2009ensemble,mandic2011filter}. Using fractional Gaussian noise as a model for broadband noise, it has been shown that the built-in adaptivity of EMD makes it behave spontaneously as a `wavelet-like' filter, i.e.~the result of sifting is a `detail' and a `trend' \cite{flandrin2004empirical}.

Efforts have also been made to replace the sifting algorithm with alternate formulations which are more mathematically grounded, such as techniques based on optimization \cite{ meignen2007new, kopsinis2008investigation, hou2013data}, machine learning \cite{looney2008machine}, PDEs \cite{delechelle2005empirical, sharpley2006analysis, vatchev2008decomposition, diop2009pde, el2010analysis, el2013pde}, and Fourier analysis \cite{niang2010spectral}. We also point out other research on the EMD algorithm including the multivariate EMD \cite{damerval2005fast,rilling2007bivariate,wu2009multi,linderhed2009image,mandic2011filter, nunes2009empirical,  ThePowerAdaptive} and the complex EMD \cite{tanaka2007complex}. As a final note, not all variations of EMD utilize the IMF component or the proposed AM--FM model and thus cannot be considered as a form of HSA. Examples include the Hilbert variational decomposition \cite{feldman2006hvd, feldman2008compare, feldman2011hilbert},  the time-dependent intrinsic correlation \cite{chen2010time}, and synchrosqueezed wavelet transforms \cite{daubechies2011synchrosqueezed, wu2011one, chui2015signal}.

\subsection{The Original EMD Algorithm} \label{ssec:emd}
The original EMD and sifting algorithms proposed by Huang \cite{huang1998empirical} are listed in Algorithms \ref{alg:EMD} and \ref{alg:SIFT}. The purpose of the sifting algorithm is to iteratively identify and remove the trend from the signal, effectively acting as a high pass filter. Step \ref{step:Removal} of Algorithm \ref{alg:EMD} removes the high frequency component $\varphi_{k}(t)$, estimated during the sifting process. The process is then repeated to remove additional IMFs from the signal if they exist. The resulting decomposition is complete and sparse \cite{huang1998empirical, hou2011adaptive, hou2013data}. EMD is formulated with continuous-time signals, however, in practice, EMD is applied to discrete-time signals which may result in errant decompositions. The effects of sampling in the context of EMD have been considered by Rilling and it is generally recommended to oversample but not resample before application of EMD, so that EMD effectively behaves like a continuous operator \cite{rilling2006sampling}.

\begin{algorithm}[H]
\caption{Empirical Mode Decomposition}
  \begin{algorithmic}[1]
    \Procedure{$\lbrace \varphi_k(t)\rbrace={\mathrm{EMD}}$} { $x(t)$ }
      \State initialize: $k=0$ and $x_{-1}(t) = x(t)$
      \While{$x_{k-1}(t)\neq0$ and $x_{k-1}(t)$ is not monotonic}
        \State $\varphi_{k}(t)= {\mathrm{SIFT}(~x_{k-1}(t)~)}$ \label{step:Sift}
        \State $x_{k}(t) = x_{k-1}(t)-\varphi_{k}(t)$ \label{step:Removal}
        \State $k \gets k+1$
      \EndWhile
      \State $\varphi_{k}(t)=x_{k-1}(t)$
    \EndProcedure
  \end{algorithmic}
  \label{alg:EMD}
\end{algorithm}

\begin{algorithm}[H]
\caption{Sifting Algorithm}
  \begin{algorithmic}[1]
    \Procedure{$\varphi(t)={\mathrm{SIFT}}$} { $r(t)$ }  
	\State initialize: $e(t)\neq0$
      \While{$e(t)\neq0$} \label{step:Stop}
        \State find all local maxima: $u_p=r(t_p),\ p = 1,2,\ldots$ \label{step:findMaxima}
        \State find all local minima: $l_q=r(t_q),\ q = 1,2,\ldots$ \label{step:findMinima}
        \State interpolate: $u(t)=\text{CublicSpline}(\{t_p,u_p\})$  \label{step:interpMaxima}
        \State interpolate: $l(t)=\text{CublicSpline}(\{t_q,l_q\})$  \label{step:interpMinima}
        \State $e(t) = [u(t)+l(t)]/2 $\label{step:AvgEnv}.
        \State $r(t) \gets r(t)-e(t)$.\label{step:Alpha}
      \EndWhile
      \State $\varphi(t) = r(t)$
     \EndProcedure 
  \end{algorithmic}
  \label{alg:SIFT}
\end{algorithm}

\subsection{Improving the Sifting Algorithm} \label{ssec:sifting}
If EMD is viewed as an AM--FM decomposition technique, then the sifting algorithm is an iterative way of removing the asymmetry between the upper and lower envelopes in order to transform the input $r(t)$ into an IMF \cite{rato2008hht}. By doing so, low frequency content is discarded at every sifting iteration, effectively making the sifting algorithm behave as a high frequency filter or high frequency component tracker. Due to the doubly iterative nature of EMD and termination conditions, numerical imprecision and differing implementations can lead to very different IMFs. To achieve consistency when using EMD, Rato proposes the following constraints \cite{rato2008hht}:
\begin{itemize}[labelwidth = \widthof{~~~~~\textbf{Time reversal:}},leftmargin=!] 
	\item[\textbf{Scale:}] IMFs should scale with the signal.
	\item[\textbf{Bias:}] Any signal bias, should only be reflected in the trend.
	\item[\textbf{Identity:}] The EMD of an IMF should be the IMF itself.\footnote{As Rato points out, this constraint may be relaxed in some cases \cite{rato2008hht}.}
	\item[\textbf{Time-reversal:}] Time-reversal of the signal should time-reverse the IMFs.
\end{itemize}

Several improvements have been made to the sifting algorithm (Algorithm \ref{alg:SIFT}) to improve the decomposition accuracy \cite{rato2005Modified, rato2008hht}:
\begin{enumerate}
	\item Improving stop criterion robustness in Step \ref{step:Stop} 
	\item Providing a consistent method for identification of extrema in Steps \ref{step:findMaxima} and \ref{step:findMinima} 
	\item Addressing interpolation end effects in Steps \ref{step:interpMaxima} and \ref{step:interpMinima} 
	\item Scaling the mean envelope removal in Step \ref{step:Alpha} 
\end{enumerate}
Although several stopping criteria have been proposed \cite{huang2003confidence,rilling2003algorithms,rato2008hht}, Rato suggested the use of a resolution factor which is the ratio between the energy of the signal at the beginning of sifting $r(t)$ and the energy of the average of the envelopes $e(t)$. If this ratio increases above a predetermined threshold, then the IMF computation terminates. We have found that Rato's use of a parabolic interpolator \cite{rato2008hht} to identify the extrema in Steps \ref{step:findMaxima} and \ref{step:findMinima} is convenient because it uses as few as three samples and it avoids the classificatory function to determine if a particular sample is a maxima, minima, or neither.  End effects appear due to the fact that a given interpolator may not be a good extrapolator \cite{rato2008hht}. In order to deal with this problem, Rato suggested to insert artificial minima and maxima in order to control the behavior of the interpolator. In addition, the mean envelope removal is scaled in Step \ref{step:Alpha},
\begin{equation}
	r(t) \gets r(t) - \alpha e(t)
	\label{eq:introduceAlpha}
\end{equation}
where the step-size $0 < \alpha \leq 1$ increases the number of sifting iterations but improves stability and robustness of the resulting IMFs \cite{rato2008hht}. 

\begin{figure}[ht]
\centering
 		\centering
  		\subfigure[]{
  		\includegraphics[width = 0.56\linewidth]{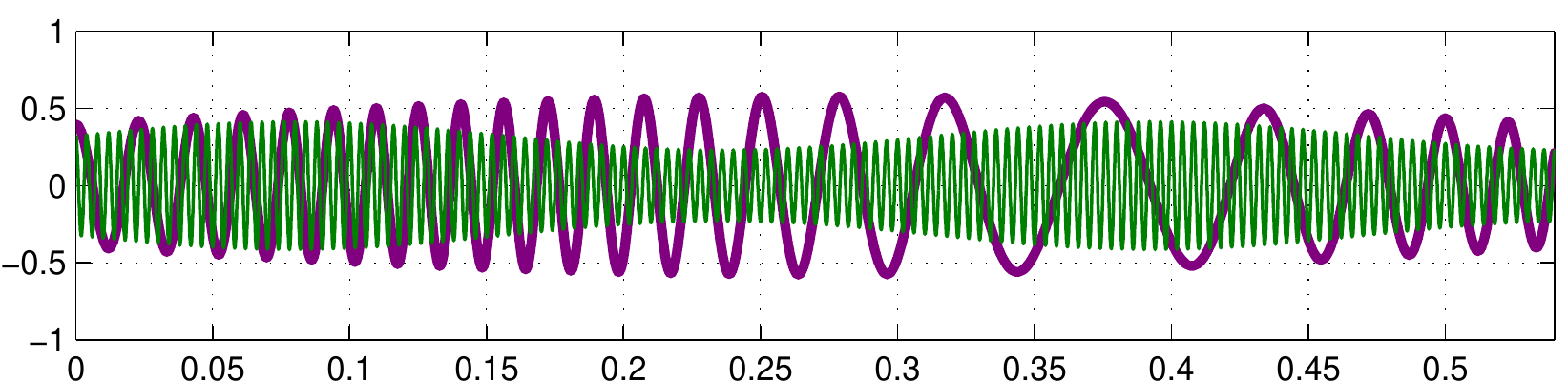}
  		\label{fig:SiftingA}
  		}
  		\subfigure[]{
  		\includegraphics[width = 0.56\linewidth]{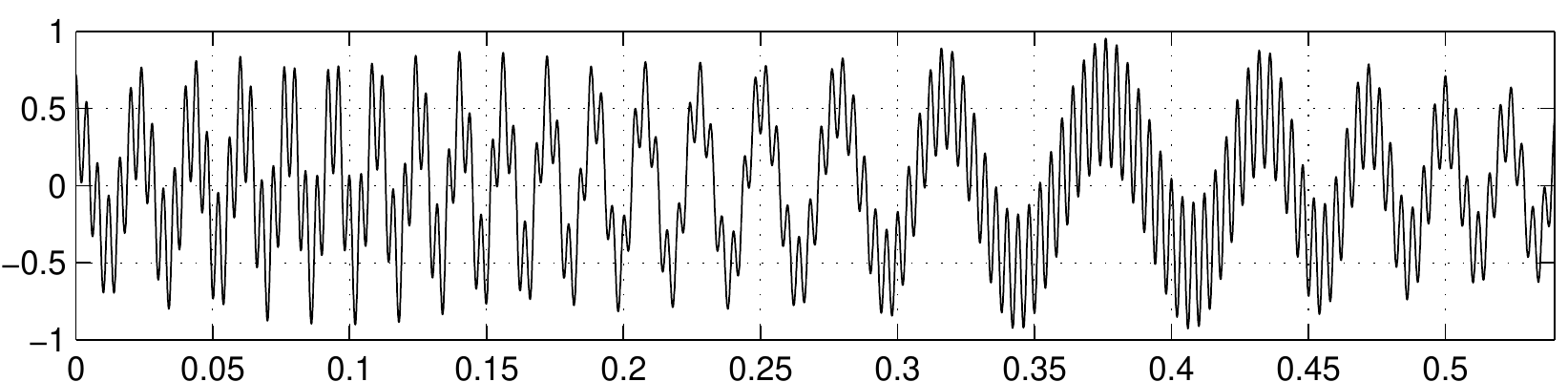}
  		\label{fig:SiftingB}
  		}
  		\subfigure[]{
  		\includegraphics[width = 0.56\linewidth]{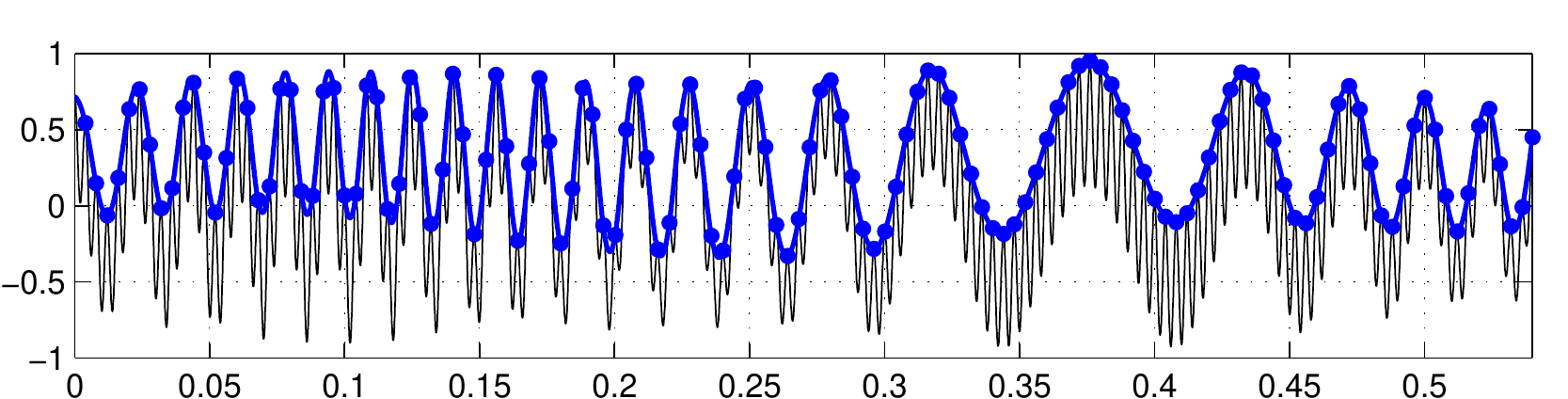}
  		\label{fig:SiftingC}
  		}
  		\subfigure[]{
  		\includegraphics[width = 0.56\linewidth]{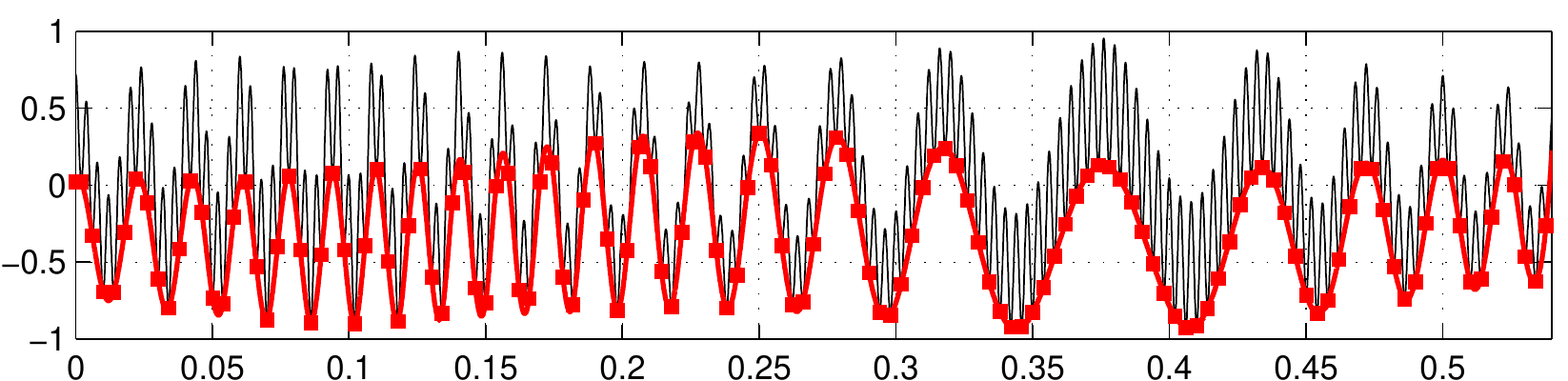}
  		\label{fig:SiftingD}
  		}
  		\subfigure[]{
  		\includegraphics[width = 0.56\linewidth]{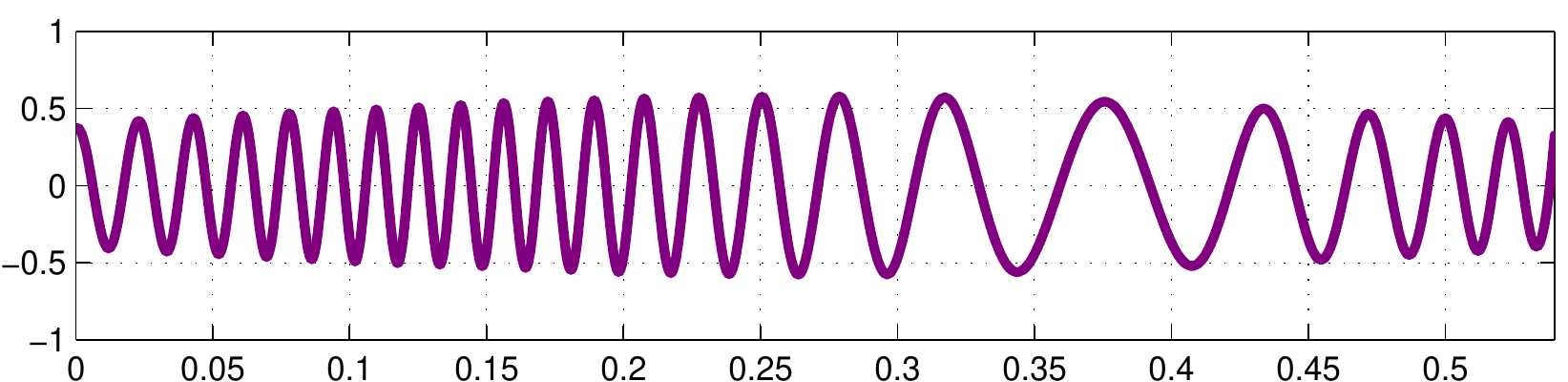}
  		\label{fig:SiftingE}
  		}
  		\subfigure[]{
  		\includegraphics[width = 0.56\linewidth]{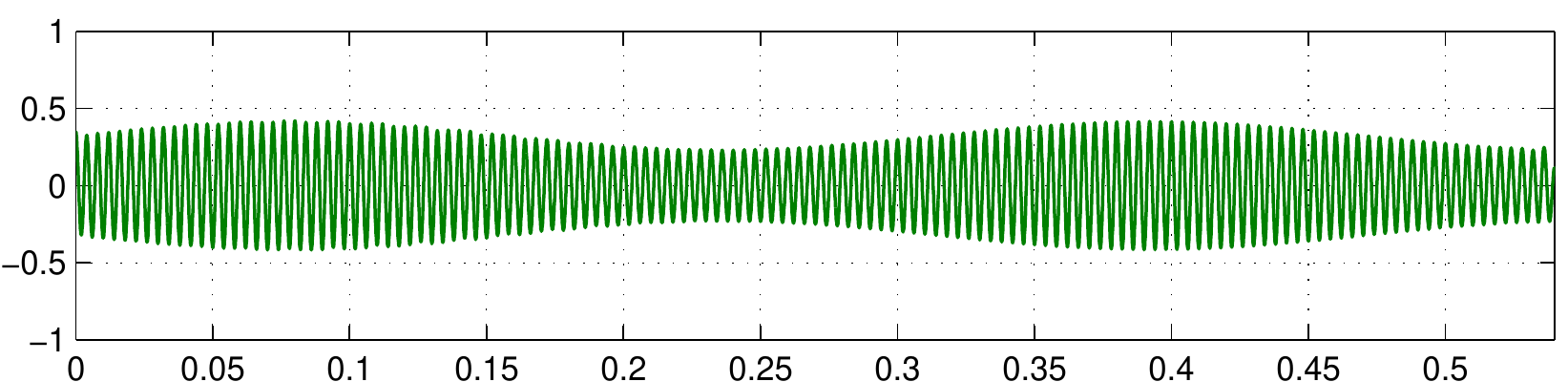}
  		\label{fig:SiftingF}
  		}
    \caption{This sequence of plots illustrates the steps of a first iteration of the sifting algorithm. In (a) the example signal composed of the components, $\varphi_0(t)$ (\textcolor{MyGreen}{\solidSrule[3.5mm]}) and $\varphi_1(t)$ (\textcolor{MyPurple}{\solidLrule[3.5mm]}); (b) the superposition of the components $r(t)$ ({\solidSrule[3.5mm]}) and the input to the sifting algorithm; (c)-(d) the upper envelope $u(t)$ ({\protect\lineDotBlueS}) and lower envelope $l(t)$  ({\protect\lineSqrRedS}) of $r(t)$; (e) average of the upper and lower envelopes $e(t) \approx \varphi_1(t)$; and (f) IMF estimate at first iteration $r(t)-e(t) \approx \varphi_0(t)$.}
  		\label{fig:Sifting}
\end{figure}

We illustrate a single iteration of the sifting algorithm in Fig.~\ref{fig:Sifting}. In Fig.~\ref{fig:SiftingA}, we show two unknown components, $\varphi_0(t)$ (\textcolor{MyGreen}{\solidMrule[3.5mm]}) and $\varphi_1(t)$ (\textcolor{MyPurple}{\solidLrule[3.5mm]}) and in Fig.~\ref{fig:SiftingB}, we show the signal under analysis $r(t)=\varphi_0(t)+\varphi_1(t)$ ({\solidSrule[3.5mm]}) which is the input to the sifting algorithm. Figs.~\ref{fig:SiftingC} and (d) show the location and interpolation of the extrema as given in Steps \ref{step:findMaxima} - \ref{step:interpMinima}. Interpolations lead to estimates of the upper envelope $u(t)$ (\lineDotBlue) and lower envelope $l(t)$ (\lineSqrRed) of $r(t)$. The average of the upper and lower envelopes $e(t)\approx\varphi_1(t)$ is shown in Fig.~\ref{fig:SiftingE}. At the end of the first iteration of sifting,  $r(t)-e(t) \approx \varphi_0(t)$ and is shown in Fig.~\ref{fig:SiftingF}.

\subsection{Improving the EMD Algorithm} \label{ssec:IEMD}
\begin{figure}[ht]
\centering
  	\begin{minipage}[b]{0.99\linewidth}
  		\centering
  		\subfigure[]{
  		\includegraphics[width = 0.75\linewidth]{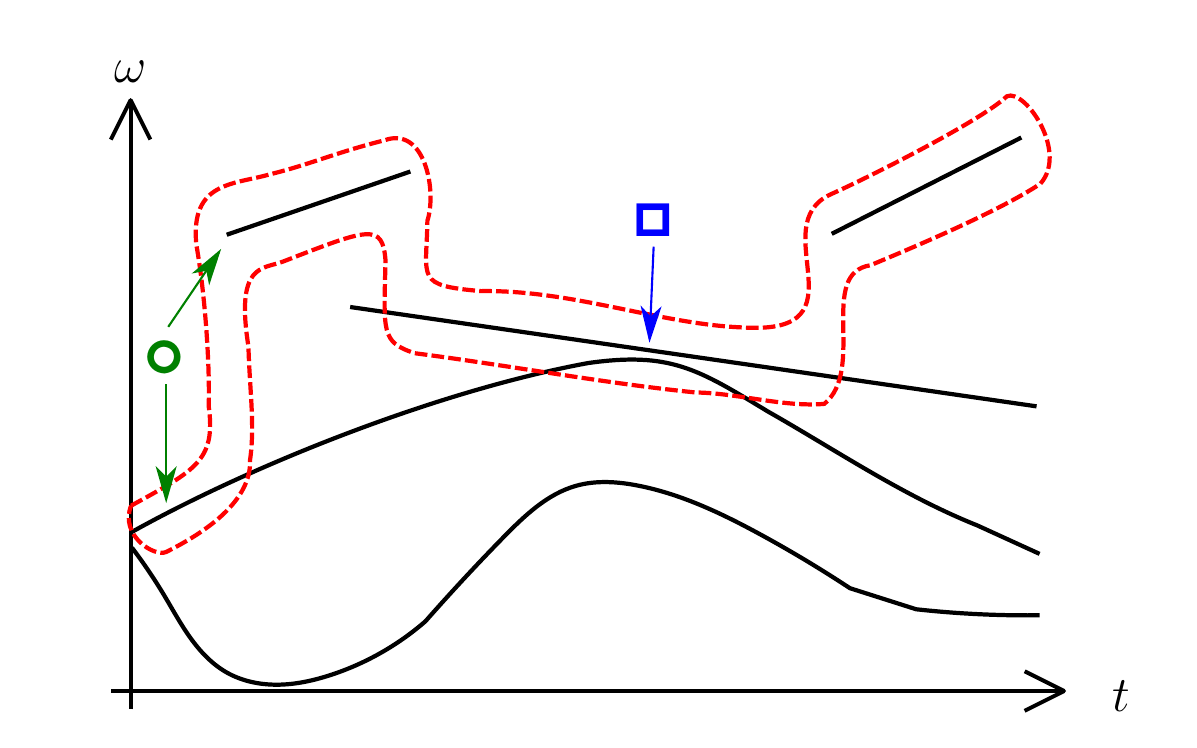}
  		\label{fig:EMD_IllustrationA}
  	}
  	\end{minipage}	
  	\begin{minipage}[b]{0.99\linewidth}
  		\centering	
  		\subfigure[]{
  		\includegraphics[width = 0.75\linewidth]{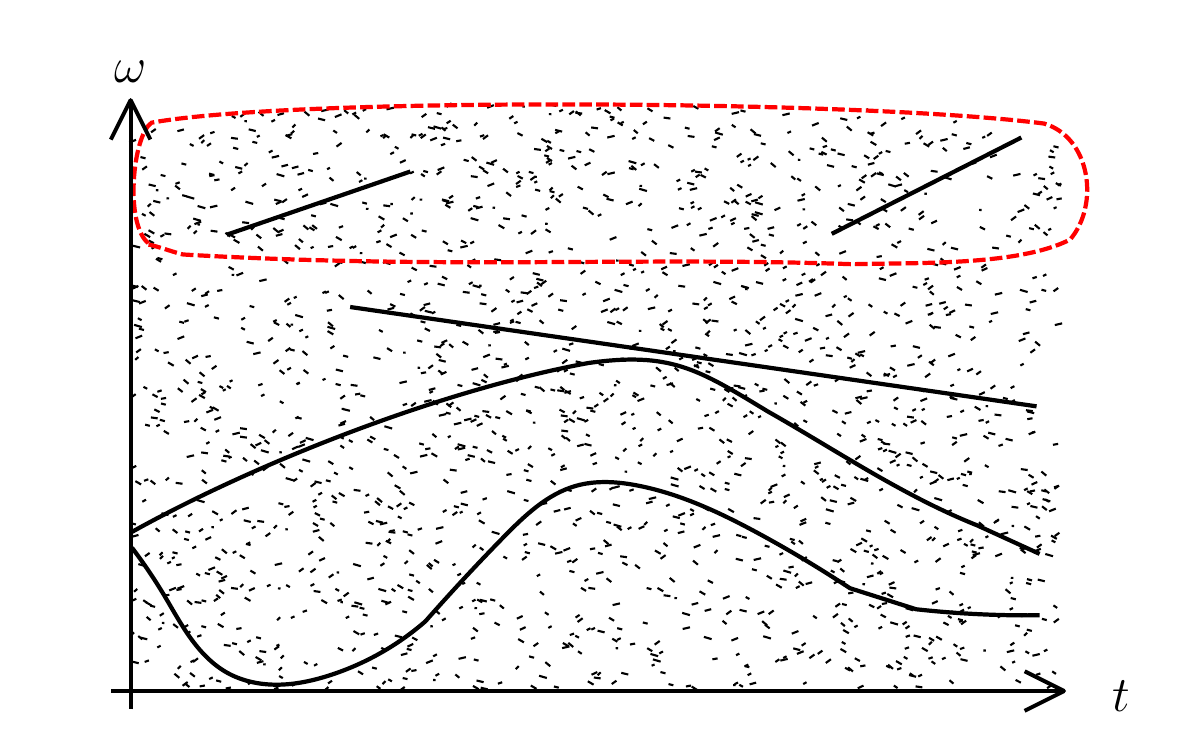}
  		\label{fig:EMD_IllustrationB}
  	}
  	\end{minipage}
  	\caption{In (a) and (b), the assumed components are indicated with \solidMrule[3.5mm] and the first component or high frequency IMF, identified with the sifting algorithm, is indicated within the \textcolor{MyRed} {\protect\dashedrule} frame. In (a), the mode mixing problem is apparent where we see components of disparate scales being in the same IMF (indicated by \textcolor{MyGreen}{\CircPipe}) and components of similar scale in the same IMF (indicated by \textcolor{MyBlue}{\SquarePipe}). In (b), adding noise and ensemble averaging may assist in resolving mode mixing.}
  	  \label{fig:EMD_Illustration}	
\end{figure}

The major problem in the EMD algorithm is mode mixing, which is defined as a single IMF either consisting of components of disparate scales or components of similar scale residing in the same IMFs \cite{wu2009ensemble}.  Mode mixing is a consequence of signal intermittency, or more specifically relative component intermittency. As a result, the particular component(s) tracked by the sifting algorithm in a particular IMF at any instant may change as intermittent components begin or end \cite{DeeringMasking}. This is illustrated in Fig.~\ref{fig:EMD_IllustrationA}, where in the left part of the figure, we show components of disparate scales being in the same IMF denoted by \textcolor{MyGreen}{\CircPipe}, while in the center part of the figure we show two components of similar scale in the same IMF denoted by \textcolor{MyBlue}{\SquarePipe}. The ability of EMD to resolve two components considering both the relative IAs and IFs of components, was examined and quantified by Rilling \cite{rilling2008one}. However as was noted, resolving closely-spaced components may not be the ultimate goal, provided that the decomposition is suitably matched to some meaningful interpretation \cite{rilling2008one}.

Two commonly used methods of mitigating mode mixing are EEMD \cite{wu2009ensemble} and tone masking \cite{DeeringMasking}. EEMD (given in Algorithm \ref{alg:EEMD}, where $\mathbb{E}[\cdot]$ denotes the expectation) utilizes zero-mean white noise, $w^{(i)}(t)$ to perturb the signal so a component may be tracked properly over an ensemble average. As the illustration in Fig.~\ref{fig:EMD_IllustrationB} shows, noise can be used to assist the sifting algorithm. Inserting noise with high enough power gives the sifting algorithm something to track when the highest frequency component is intermittent, then vanishes in the ensemble average. Although injecting noise can help to track components properly, a carefully designed masking signal can result in better performance. A situation where a carefully designed masking signal may be beneficial is illustrated in Fig.~\ref{fig:EMD_IllustrationA} denoted by \textcolor{MyBlue}{\SquarePipe}.

\begin{algorithm}[H]
\caption{Ensemble Empirical Mode Decomposition}
  \begin{algorithmic}[1] 
    \Procedure{$\lbrace \bar\varphi_k(t)\rbrace={\mathrm{EEMD}}$} { $x(t)$ }
      \State initialize: $I$ is the number of trials
      \State $\varphi^{(i)}_k(t) = \mathrm{EMD}\left(~x(t) + w^{(i)}(t)~\right),~i=1,\ldots,I$
      \State $\bar\varphi_k(t) = \mathbb{E}[ \varphi^{(i)}_k(t)]$ 
    \EndProcedure
  \end{algorithmic}
  \label{alg:EEMD}
\end{algorithm}

EEMD is not without its disadvantages as it is more computationally complex, loses the perfect reconstruction property, propagates IMF estimation error, results in inconsistent numbers of IMFs across the trials, and the resulting set of averaged IMFs $\{\bar\varphi_k(t)\}$ are not necessarily IMFs \cite{torres2011complete}. Torres proposed the ``complete EEMD'' given in Algorithm \ref{alg:CEEMD} to address some of these issues \cite{torres2011complete}. Complete EEMD defines a procedure $\mathrm{EMD}_k(\cdot)$ which returns the $k$th IMF obtained using EMD \cite{torres2011complete}. This method of EEMD requires fewer sifting iterations, a smaller ensemble size, and recovers the completeness property of the original EMD algorithm to within the numerical precision of the computer \cite{torres2011complete}. 

\begin{algorithm}[H]
\caption{Complete EEMD}
  \begin{algorithmic}[1] 
    \Procedure{$\lbrace \bar\varphi_k(t)\rbrace={\mathrm{CEEMD}}$} { $x(t)$ }
	\State initialize: $k=1$, $\beta_k$ is a SNR factor,\par and $I$ is the number of trials
	\State $\bar\varphi_0(t) = \displaystyle{\frac{1}{I}\sum_{i=1}^{I}} \mathrm{SIFT}(~x(t)+\beta_0 w^{(i)}(t)~)$
     \State $x_{0}(t) = x(t)-\bar\varphi_{0}(t)$
      \While{$x_{k-1}(t)\neq0$ and $x_{k-1}(t)$ is not monotonic}
		\State \hspace{-3mm} $\bar\varphi_{k}(t) = \displaystyle{\frac{1}{I}\sum_{i=1}^{I}}\mathrm{SIFT}(~x_{k-1}(t)+\beta_k \mathrm{EMD}_k(~w^{(i)}(t)~)$
		\State $x_{k}(t) = x_{k-1}(t)-\bar\varphi_{k}(t)$
		\State $k\leftarrow k+1$
      \EndWhile
            \State $\bar\varphi_{k}(t)=x_{k-1}(t)$
    \EndProcedure
  \end{algorithmic}
  \label{alg:CEEMD}
\end{algorithm}

Rather than using a noise signal, a deterministic signal $v(t)$ can also be used as a perturbation and then removed after IMF estimation. The ``tone masking'' technique is given in Algorithm \ref{alg:TM} and can be used, for example, in place of Step \ref{step:Sift} in Algorithm \ref{alg:EMD} \cite{DeeringMasking}. Advanced forms of tone masking have been proposed and are termed Signal-Assisted EMD (SA-EMD) \cite{FlandrinHuang2005, senroy2007two, guanlei2009time, li2011sinusoidal}. These methods have additional advantages in their ability to track closely-spaced components, but require careful selection of the masking signal. 

\begin{algorithm}[H]
\caption{Tone Masking}
  \begin{algorithmic}[1]
    \Procedure{$\bar\varphi(t)={\mathrm{TM}}$} { $x(t),v(t)$ }  
	  	\State $x^{(+)}(t) = x(t) + v(t)$ and $x^{(-)}(t) = x(t) - v(t)$
	  	\State $\varphi^{(+)}(t) = \mathrm{SIFT}(~x^{(+)}(t)~)$ and $\varphi^{(-)}(t) = \mathrm{~SIFT}(~x^{(-)}(t)~)$\label{step:twoSifts} 
	  	\State $\bar\varphi(t) = \displaystyle{ \frac{\varphi^{(+)}(t)+\varphi^{(-)}(t)}{2} }$
     \EndProcedure 
  \end{algorithmic}
  \label{alg:TM}
\end{algorithm}

\subsection{IMF Demodulation} \label{ssec:demod}

In order to obtain the IA/IF parameters in the AM--FM model, we are required to demodulate the IMFs returned by EMD.  One approach, used by Huang, is to apply the HT to each IMF in order to obtain estimates of IA and IF. This analysis is often referred to as the Hilbert-Huang Transform (HHT) \cite{huang1998empirical}. Using the HT for IMF demodulation is \emph{inappropriate} because as argued, the HT assumes SHCs and HC---which is \emph{not} true since an IMF is in general an AM--FM component without HC. A similar argument can be made of many of the other demodulation methods that have been considered \cite{huang2009instantaneous}. A second approach, as advocated in Part II, is to let the assumed form of the AM--FM component imply a complex extension. When viewed in the context of  LSA, the definition of the IMF given in Section \ref{sec:EMD} forces a unique complex extension to the IMF---justifying our view of the IMF as a latent component specified by a real signal.

Inherent in \textbf{C2} of the IMF definition, are the following assumptions:
\begin{itemize}[labelwidth = \widthof{~~~~~A1:},leftmargin=!]
    \item[\textbf{A1:}] $a(t_p)=|s(t_p)|$, where $|s(t_p)|$ are the extrema of $s(t)$.
    \item[\textbf{A2:}] $a(t),t  \notin\{ t_p \}$ is inferred by cubic spline interpolation.
\end{itemize}
The first assumption, \textbf{A1} can be viewed as quadrature forced to zero at $\{t_p\}$. To see this, note that from (\ref{eq:AMFMcomp}), we have
\begin{equation}
    |a(t)| = \sqrt{s^2(t)+\sigma^2(t)}
\end{equation}
and thus, $\sigma(t_p)=0$ and $a(t_p) = |s(t_p)|$. \textbf{A1} also implies non-negativity of $a(t_p)$, however, non-negativity of $a(t)$ is not guaranteed for all $t$ due to the interpolation.

The second assumption, \textbf{A2} can be viewed as a relative bandlimiting condition on $a(t)$ controlled by two factors: the choice of interpolator and the density of extrema. To see this, note that $a(t)$ between extrema of $s(t)$ is defined by an interpolator.  This implies that $a(t)$ will only have as much variation between extrema as the interpolator allows. However, with dense extrema (high IF) much less restrictive constraints have been imposed on $a(t)$ than if extrema are sparse (low IF)---thus our view as a relative bandlimiting condition on $a(t)$. Also note, using an interpolator other than the cubic spline is likely to change the IA, effectively changing the resulting IMFs \cite{rato2008hht,hawley2010some}.

In the IMF definition, \textbf{C1} requires that the number of extrema and the number of zero crossings must be either equal or differ at most by one. A general AM--FM component may not satisfy this condition, for example as a result of sign reversal of $\omega(t)$ or in cases where \textbf{A1} and \textbf{A2} are not satisfied. Assumption of positive IF on the AM--FM component in Part II, eliminates the possibility of sign reversal of the IF but may not be true for all signal classes.
 
Rato proposed an AM demodulation approach given in Algorithm \ref{alg:AMest}, that is consistent with the decomposition assumptions in Algorithm \ref{alg:SIFT} \cite{rato2008hht}. Starting with an IMF estimate $\hat{\varphi}(t)$, we obtain an estimate for IA $\hat{a}(t)$ which can then be used to estimate the IF via the FM signal
\begin{equation}
    \hat{s}_{\mathrm{FM}}(t) = \frac{\hat{\varphi}(t)}{\hat{a}(t)}.
    \label{eq:sFM}
\end{equation}
However, the estimate in (\ref{eq:sFM}) may result in $|\hat{s}_\mathrm{FM}(t)|> 1$. Thus, Huang proposed an iterative normalization procedure, given in Algorithm \ref{alg:iterIMFdemod}, which removes the AM from the signal to obtain a more accurate $\hat{s}_\mathrm{FM}(t)$ \cite{huang2009instantaneous}. Although the iterative procedure improves demodulation accuracy, we noticed that it can be susceptible to oscillating artifacts introduced by overfitting of the cubic spline interpolator. As a result, we have found that these artifacts can be minimized by limiting the number of iterations to three. 
 
\begin{algorithm}[H]
\caption{IMF IA Estimation}
  \begin{algorithmic}[1]
    \Procedure{$\hat{a}(t)=\mathrm{IAest}$} { $\hat{\varphi}(t)$ }  
            \State $r(t) = |{\hat{\varphi}}(t)|$
            \State find all local maxima: $u_p=r(t_p),\ p = 1,2,\ldots$ 
            \State interpolate: $\hat{a}(t)=\text{CublicSpline}(\{t_p,u_p\})$
     \EndProcedure 
  \end{algorithmic}
  \label{alg:AMest}
\end{algorithm}

\begin{algorithm}[H]
\caption{Obtaining a real FM signal from an IMF}
  \begin{algorithmic}[1]
    \Procedure{$\hat{s}_\mathrm{FM}(t)=\mathrm{iterAMremoval}$} { $\hat{\varphi}(t)$ }  
        \State initialize: $g(t) = \hat{\varphi}(t)$, ${b}(t)\neq1$, and $n=0$
           \While{${b}(t)\neq1$ and $n<3$} 
                \State ${b}(t)=\text{IAest}(~g(t)~)$
                \State $g(t) \gets g(t)/b(t)$ 
                \State $n\gets n+1$
           \EndWhile
           \State $\hat{s}_\mathrm{FM}(t) = g(t)$
     \EndProcedure 
  \end{algorithmic}
  \label{alg:iterIMFdemod}
\end{algorithm}

We can directly obtain two estimates of the IF, $\pm\hat\omega(t)$ by substituting $a(t)=1$ in (\ref{eq:AMFMcompB}), $s(t)=\hat{s}_\mathrm{FM}(t)$ in (\ref{eq:AMFMcompC}), and then computing (\ref{eq:NonNegInstFreq}). Direct FM demodulation is not straightforward because different approaches exist for obtaining the IF from $\hat{s}_\mathrm{FM}(t)$, that although are mathematically equivalent, may differ in numerical stability \cite{huang2009instantaneous}.

\begin{figure}[H]
    \centering
    \ifAnimations
\else  
		\includegraphics[width=8 cm]{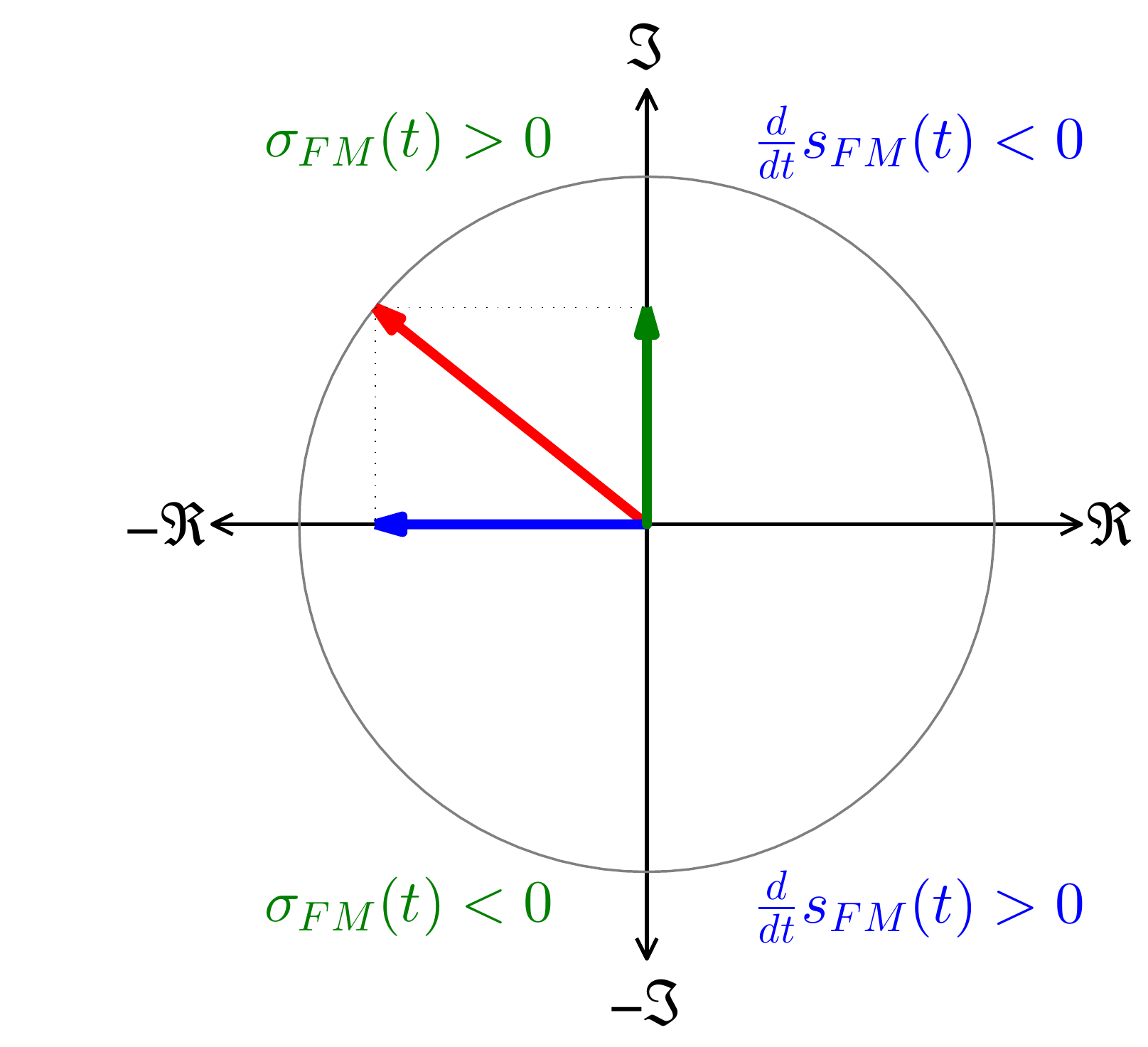}
	\fi	
	\caption{In this figure (animated in web version), $\hat{s}_\mathrm{FM}(t)$, $\hat{\sigma}_\mathrm{FM}(t)$ are represented by the blue and green vectors, respectively.  The amplitude-normalized $\hat{\psi}_{FM}(t)$, is represented by the red vector.  The magnitude of $\hat{\sigma}_\mathrm{FM}(t)$ is easily calculated. The sign of $\hat{\sigma}_\mathrm{FM}(t)$ is obtained as follows. We note that when $\hat{s}_\mathrm{FM}(t)$ is decreasing (blue vector moves left), $\hat{\sigma}_\mathrm{FM}(t)$ is always positive (green vector is in the upper half plane) and when $\hat{s}_\mathrm{FM}(t)$ is increasing (blue vector moves right), $\hat{\sigma}_\mathrm{FM}(t)$ is always decreasing. Thus by reversing the sign of the derivative of $\hat{s}_\mathrm{FM}(t)$ we obtain the sign of $\hat{\sigma}_\mathrm{FM}(t)$.}
    \label{fig:IMFdemod}
\end{figure}

We propose a FM demodulation method to address the numerical stability issues.  We begin by estimating the quadrature in (\ref{eq:AMFMcompC}) as
\begin{equation}
	\hat{\sigma}_\mathrm{FM}(t) = -\text{sgn}\left[\frac{d}{dt}\hat{s}_\mathrm{FM}(t)\right] \sqrt{1^2-\hat{s}_\mathrm{FM}^2(t)}
	\label{eq:SigmaEST}
\end{equation}
where $-\text{sgn}\left[\frac{d}{dt}\hat{s}_\mathrm{FM}(t)\right]$ is required to obtain an appropriate four quadrant estimate with assumed positive IF. In Fig.~\ref{fig:IMFdemod} , we have provided an illustration (animated in the web version) to assist the reader with understanding of how we arrive at (\ref{eq:SigmaEST}). The computationally unstable points, $\{t_0\}$ occur where $\hat\sigma_\mathrm{FM}(t_0)= 0$, thus we can replace a small range around these points $(t_0-\epsilon, t_0+\epsilon)$ with interpolated values. Then
\begin{equation}
	\hat{\theta}(t) = \text{arg}\left[\hat{s}_\mathrm{FM}(t) + j \hat{\sigma}_\mathrm{FM}(t)\right]
	\label{eq:ThetaEST}
\end{equation}
and the IF is obtained with (\ref{eq:NonNegInstFreq}). We can optionally smooth the resulting IF estimate. Our complete IMF demodulation algorithm is listed in Algorithm \ref{alg:IMFdemod}.

We note to the reader that the proposed IMF demodulation procedure, which does not involve the HT or FT, can in some cases return the same results but in general this is not true.  For example, demodulation of $\cos(\omega_0 t)$ with either the IMF or HT demodulation returns a SHC with HC. On the other hand, HT demodulation of the triangle waveform returns a single AM-FM component with HC as given in (\ref{eq:trianglesingleamfmifP1}), while IMF demodulation returns IF given by (\ref{eq:trianglesinglefmifP1}) and constant IA. 

\subsection{Proposed Algorithm for HSA Assuming IMFs} \label{ssec:EMDproposed}
The various recommendations and improvements listed this section appear to have never been combined and investigated together.  Therefore, we propose to incorporate the most desirable features of the complete EEMD and tone masking for addressing the mode-mixing problem, Rato's improvements to the sifting algorithm, and our proposed demodulation method, into a single algorithm given in Algorithm \ref{alg:ProposedEMD}. 

\begin{algorithm}[H]
\caption{IMF demodulation}
  \begin{algorithmic}[1]
    \Procedure{$[ \hat{a}(t),\hat{\omega}(t) ]=\mathrm{IMFdemod}$} { $\hat{\varphi}(t)$ }  
        \State $\hat{a}(t) = \mathrm{IAest}(~\hat\varphi(t)~)$
        \State $\hat{s}_\mathrm{FM}(t)=\mathrm{iterAMremoval}(~\hat{\varphi}(t)~)$
        \State $\hat{\sigma}_\mathrm{FM}(t) = \displaystyle{ -\text{sgn}\left[\dfrac{d}{dt}\hat{s}_\mathrm{FM}(t)\right] \sqrt{1^2-\hat{s}_\mathrm{FM}^2(t)} }$
        \State Find $\{t_{0}\}$ such that $\hat{\sigma}_\mathrm{FM}(t_0)=0$
        \State For each $t_0$, replace $\left(\hat{\sigma}_\mathrm{FM}(t_0-\epsilon),\hat{\sigma}_\mathrm{FM}(t_0+\epsilon)\right)$\par with interpolation
        \State $\hat{\omega}(t) = \displaystyle{ \dfrac{d}{dt}\text{arg}\left[\hat{s}_\mathrm{FM}(t) + j \hat{\sigma}_\mathrm{FM}(t)\right] }$
     \vspace{0.7mm}\EndProcedure 
  \end{algorithmic}
  \label{alg:IMFdemod}
\end{algorithm}

\begin{algorithm}[H]
\caption{HSA--IMF Algorithm}
\label{alg:Proposed}
  \begin{algorithmic}[1]
    \Procedure{$\lbrace  \bar\varphi_k(t), \hat{a}_k(t), \hat{\omega}_k(t)  \rbrace=\mathrm{HSA--IMF}$} { $x(t)$ }  
            \State initialize: $x_{-1}(t) = x(t)$, $k=0$, $\beta_k$ is a SNR factor, $\varepsilon$ is an energy threshold, and $I$ is the number of trials
           \While{$\int |x_{k-1}(t)|^2 dt > \varepsilon$\par and $x_{k-1}(t)$ is not monotonic}\label{step:outermost}
     		\State $\bar\varphi_{k}(t) = \dfrac{1}{I}\displaystyle{\sum_{i=1}^{I} \mathrm{TM}(~{x}_{k-1}(t),~\beta_{k} v^{(i,k)}(t)~)}$\label{step:getIMF}
            \State $[ \hat{a}_{k}(t), \hat{\omega}_{k}(t) ] = \mathrm{IMFdemod}(~\bar\varphi_{k}(t)~)$ \label{step:demod}
     		\State ${x}_{k}(t)={x}_{k-1}(t)-\bar\varphi_{k}(t)$
     		\State $k \gets k+1$ 
           \EndWhile
           \State $\bar\varphi_{k}(t)=x_{k-1}(t)$
     \EndProcedure 
  \end{algorithmic}
  \label{alg:ProposedEMD}
\end{algorithm}

Although the masking signal $v^{(i,k)}(t)$ is usually a carefully chosen real-valued AM--FM signal \cite{senroy2007two, guanlei2009time}, we choose $v^{(i,k)}(t)$ as a white, Gaussian, lowpass-filtered noise with cutoff frequency below the amplitude-weighted IF \cite{strom1977amplitude} of the previously found IMF.  This replaces the need to use sifted white noise in complete EEMD, however, we do allow for more sophisticated approaches for choosing the masking signal. Note that this introduces a \emph{feedback loop} between decomposition and demodulation. This changes EMD from simply a decomposition method, to a full HSA method because the IA and IF estimates of the $k$th component are inherently computed and then used to design the masking signal for the $(k+1)$th component.

\section{Examples using the HSA--IMF Algorithm} \label{sec:Examples}
In this section, we demonstrate the HSA--IMF algorithm for signal analysis and compare to conventional STFT using both synthetic and real-world signals. In these examples, we use the proposed visualization method from Part II orthographically projected onto the time-frequency plane, for plotting the instantaneous parameters resulting from HSA.  We also plot the STFT magnitude (STFTM) using the same colormap to facilitate comparisons \cite{NiccoliPMK}. The STFT is computed using a 4096-point hamming window with a 1 sample frame advance. We have provided {\sc matlab} functions for these examples \cite{HSAlink}.

\subsection{Synthetic Signals} \label{ssec:ExamplesSynth}
We provide two examples using synthetic signals with known underlying signal models. The synthetic signal examples demonstrate the proposed algorithm on a single component in a noiseless environment.  In this case, mode mixing is not a problem and we initialize $I=1$, $\alpha=0.95$, and $\beta_k=0$ in Algorithm \ref{alg:ProposedEMD}.

\begin{figure*}[t]
\centering
  	\begin{minipage}[b]{0.49\linewidth}
  		\centering
  		\subfigure[]{
  		\includegraphics[width = 0.99\linewidth]{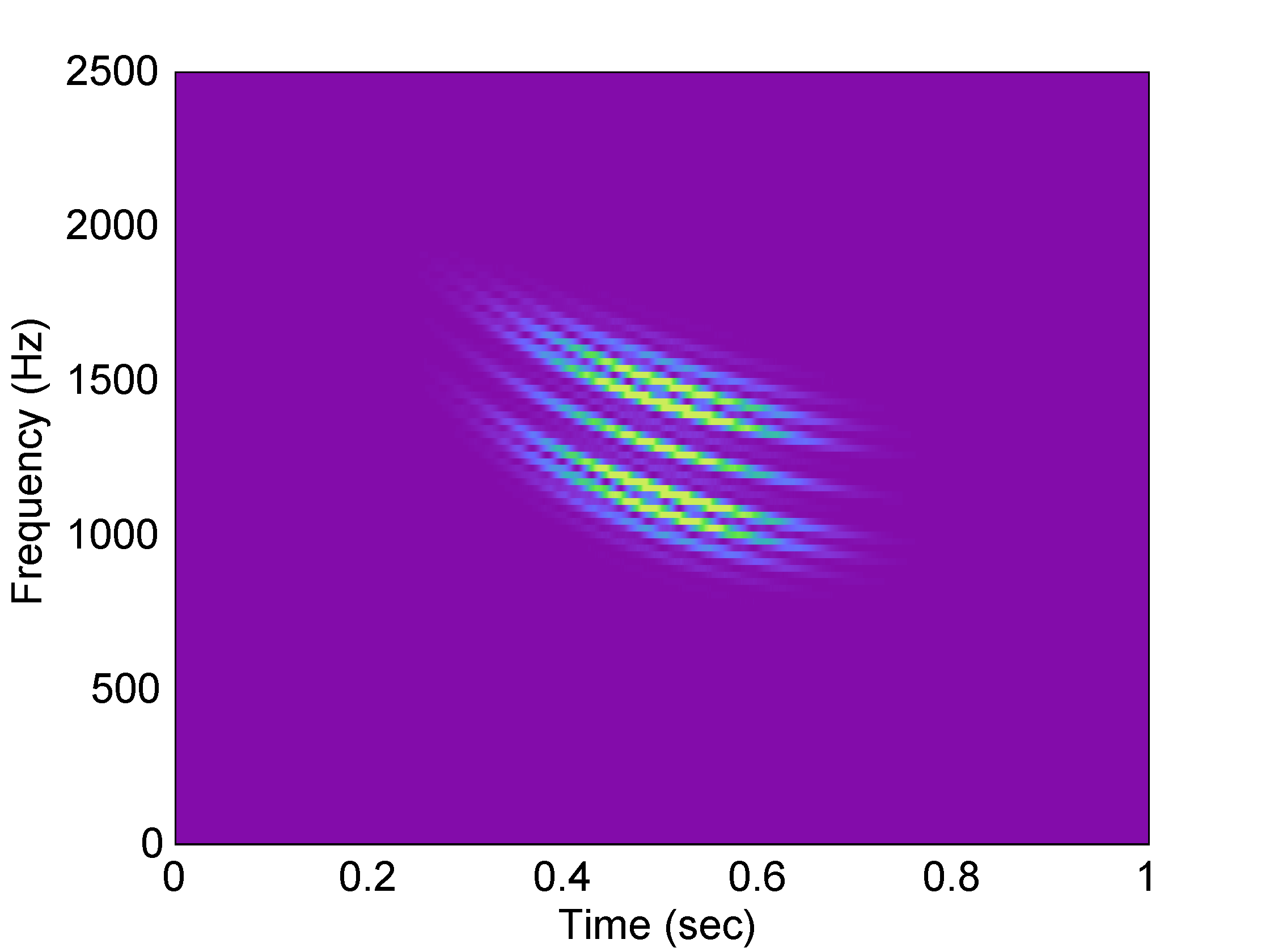}
  		\label{fig:SyntEx1a}
  	}
  	\end{minipage}	
  	\begin{minipage}[b]{0.49\linewidth}
  		\centering	
  		\subfigure[]{
  		\includegraphics[width = 0.99\linewidth]{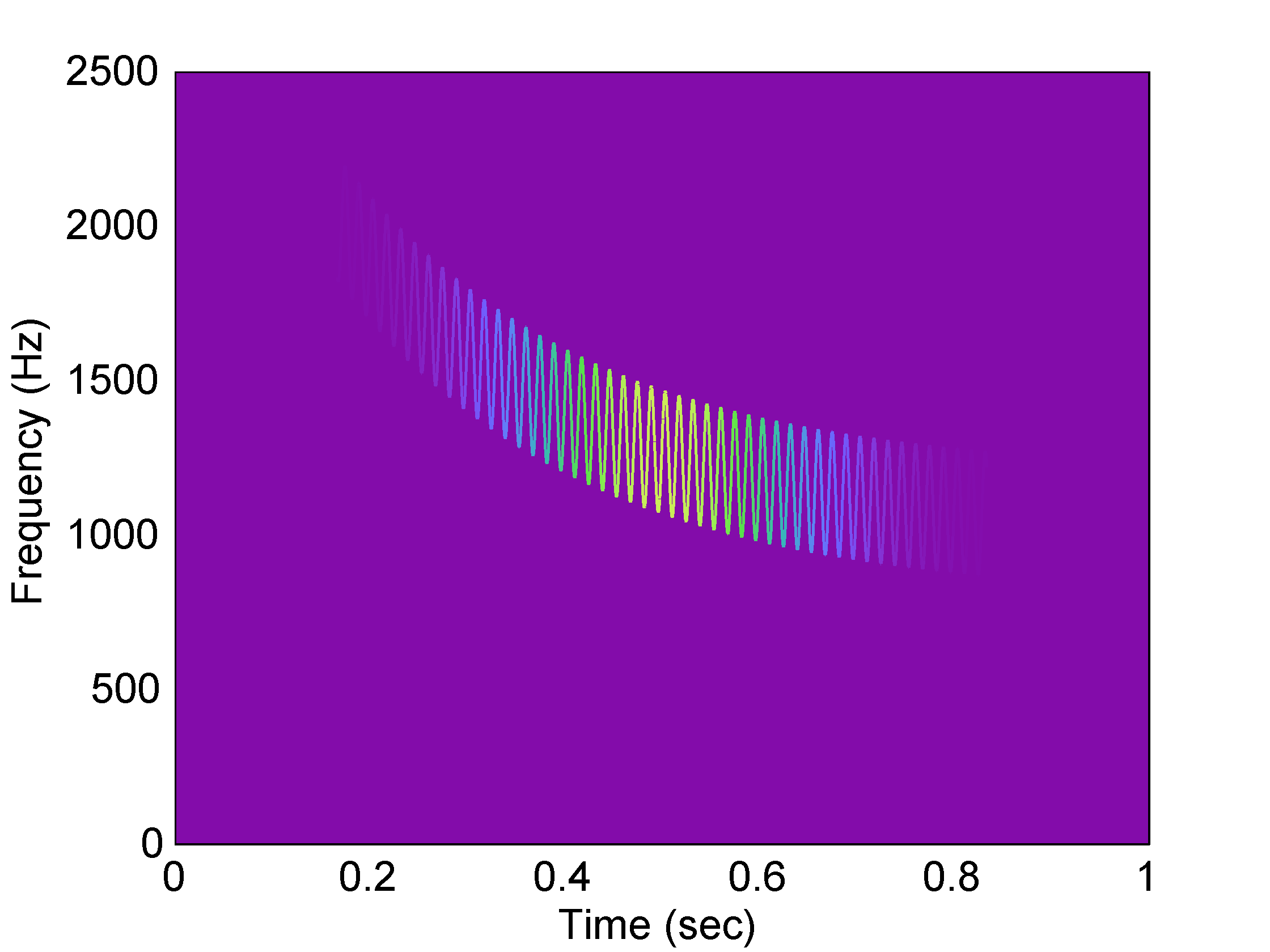}
  		\label{fig:SyntEx1b}
  	}
  	\end{minipage}
  \caption{(a) STFTM and (b) Hilbert spectrum for the fast-varying FM and slow-varying AM synthetic signal given in (\ref{eq:ex1})-(\ref{eq:ex1FMmess}) in Example 1 in Subsection \ref{ssec:ExamplesSynth}. The wideband FM message results in harmonic structure under Fourier analysis and a fast-frequency-varying component under HSA.}
  	\label{fig:SyntEx1}
\end{figure*}

In the first example, we analyze a signal with a slow-varying AM and fast-varying FM given by
\begin{equation}
    x(t) =\Re\left\{ a(t) \exp\left[j\left(6000\pi t + \int\limits_{-\infty}^{t} m(\tau) d\tau \right)\right] \right\}
    \label{eq:ex1}
\end{equation}
where the IA is
\begin{equation}
    a(t) = e^{-(t-0.5)^2/25}
    \label{eq:ex1AMmess}
\end{equation}
and the FM message is
\begin{equation}
    m(t) =  250\sin(140\pi t)+2000[\exp(-4t)-1].
    \label{eq:ex1FMmess}
\end{equation}
Fig.~\ref{fig:SyntEx1a} shows the STFTM where we see classic harmonic structure resulting from the inherent assumption of SHCs despite (\ref{eq:ex1}) consisting of only a single component.  Fig.~\ref{fig:SyntEx1b} shows a single component in the Hilbert spectrum with a fast FM oscillation consistent with the underlying signal model in (\ref{eq:ex1}). The IA in (\ref{eq:ex1AMmess}), consisting of a Gaussian envelope centered at $t=0.5$, is visible as color variation in Fig.~\ref{fig:SyntEx1a}. The FM message in (\ref{eq:ex1FMmess}) consists of a 70 Hz oscillation superimposed on a decaying exponential; this FM message is offset by a 3000 Hz carrier. The FM message is reflected by the vertical behavior in Fig.~\ref{fig:SyntEx1b} where we see the IF sweeping from 3000 Hz to 1000 Hz with a 70 Hz oscillation. Due to more appropriate assumptions of an underlying component, this example shows that the HSA--IMF parametrization can more closely match the underlying signal model than traditional methods.

\begin{figure*}[t]
\centering
  	\begin{minipage}[b]{0.49\linewidth}
  		\centering
  		\subfigure[]{
  		\includegraphics[width = 0.99\linewidth]{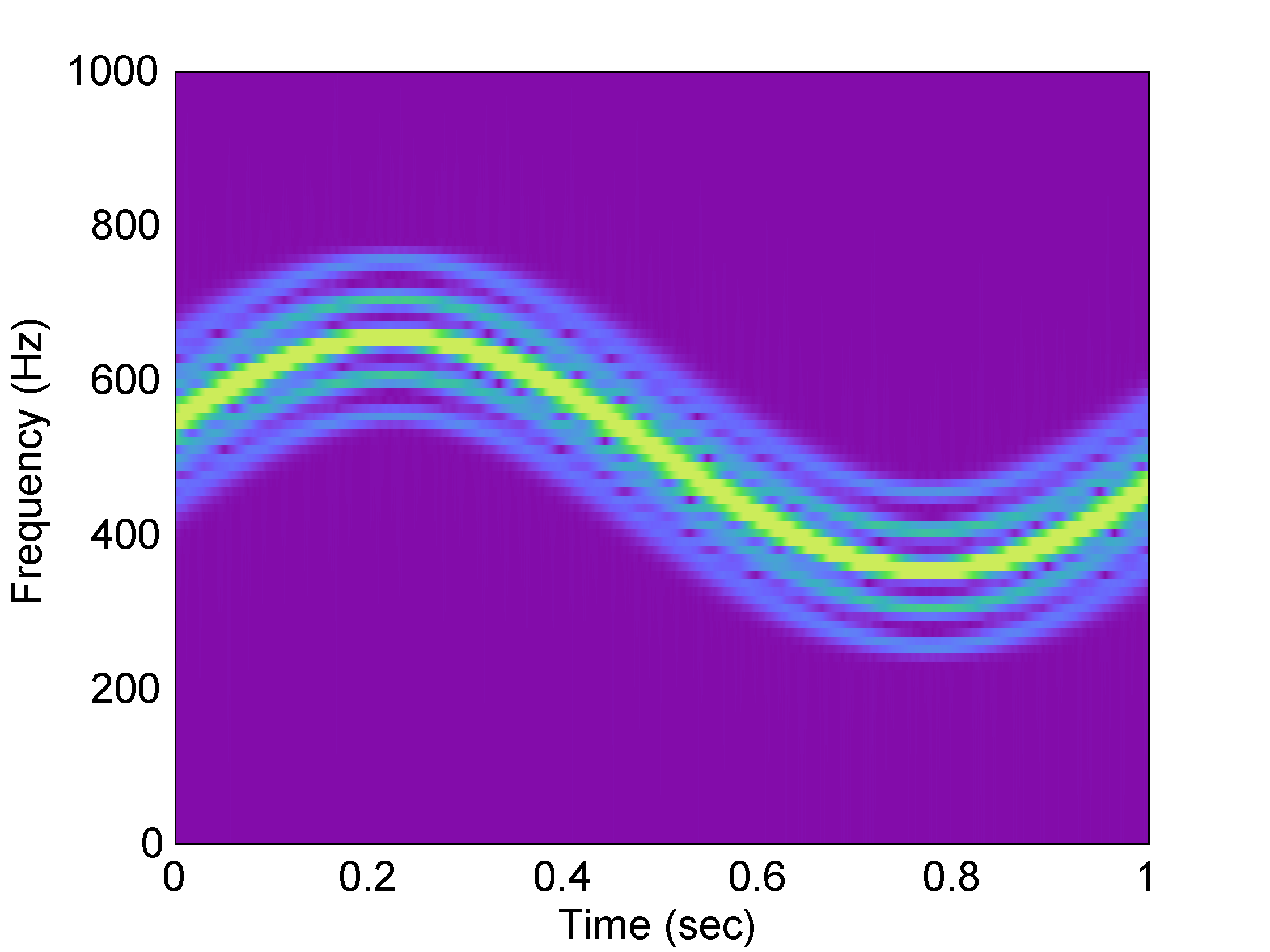}
  		\label{fig:SyntEx2a}
  	}
  	\end{minipage}	
  	\begin{minipage}[b]{0.49\linewidth}
  		\centering	
  		\subfigure[]{
  		\includegraphics[width = 0.99\linewidth]{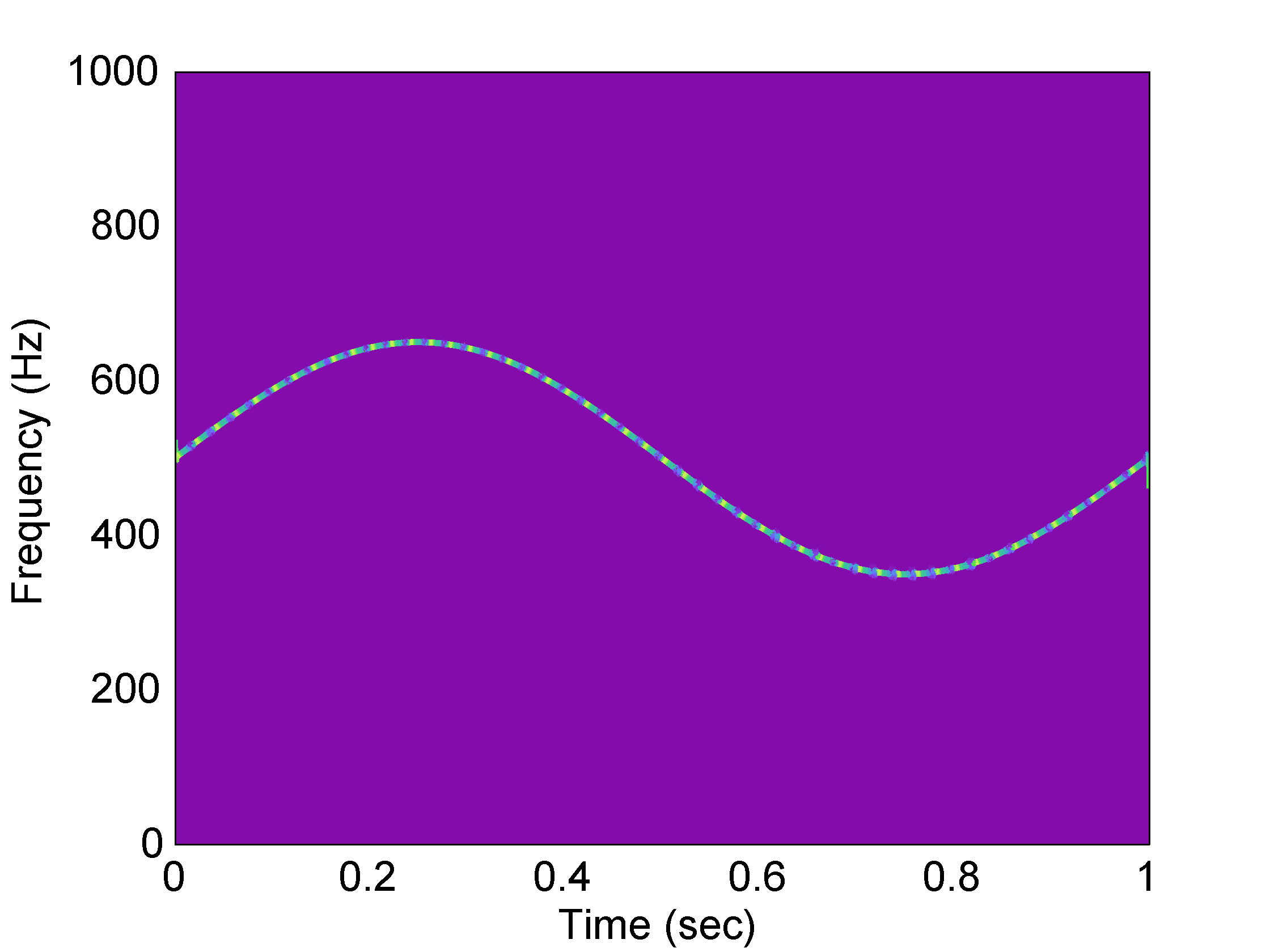}
  		\label{fig:SyntEx2b}
  	}
  	\end{minipage}
  \caption{(a) STFTM and (b) Hilbert spectrum for the fast-varying AM and slow-varying FM synthetic signal given in (\ref{eq:ex2})-(\ref{eq:ex2FMmess}) in Example 2 in Subsection \ref{ssec:ExamplesSynth}. The wideband IA results in harmonic structure under Fourier analysis and a fast-amplitude-varying component under HSA.}
  	\label{fig:SyntEx2}
\end{figure*}

In the second example, we analyze a signal with a fast-varying AM and slow-varying FM given by
\begin{equation}
    x(t) = \Re\left\{  a(t) \exp\left[j\left( 1000\pi t   + \int\limits_{-\infty}^{t}  m(\tau)  d\tau\right)\right] \right\}
    \label{eq:ex2}
\end{equation}
where the IA is
\begin{equation}
    a(t) =  \dfrac{1}{2} + \dfrac{1}{3}\sin(100\pi t) + \dfrac{1}{5}\sin(200\pi t)
    \label{eq:ex2AMmess}
\end{equation}
and the FM message is
\begin{equation}
    m(t) = 150\sin(2\pi t).
    \label{eq:ex2FMmess}
\end{equation}
Fig.~\ref{fig:SyntEx2a} shows the STFTM where we again see classic harmonic structure resulting from the inherent assumption of SHCs despite (\ref{eq:ex2}) consisting of only a single component. Fig.~\ref{fig:SyntEx2b} shows a single component in the Hilbert spectrum with a fast AM oscillation, captured by color variation in the plot line, consistent with the underlying signal model in (\ref{eq:ex2}). 

As described in Part II, signals composed of narrowband components can have similar Fourier and Hilbert spectra, however, signals composed of wideband components can have spectra which are quite different. Fourier analysis of wideband signals may result in multiple ridges in terms of spectral frequencies, which may be further decomposed into narrowband components. However, for many real-world signals, a wideband component, such as the AM--FM component in (\ref{eq:AMFMcomp}), is the more appropriate choice rather than multiple, narrowband components. The appearance of structure in the Fourier spectrum may indicate the presence of a wideband component(s) in the signal as demonstrated in Example 1, Fig.~\ref{fig:SyntEx1a} (fast-varying IF) and in Example 2, Fig.~\ref{fig:SyntEx2a} (fast-varying IA).

Both decompositions, i.e.~STFT and HSA--IMF are equally valid models for the real signal $x(t)$---the resulting decompositions simply correspond to the different assumptions of the underlying components, i.e.~constant-frequency components and IMF components, respectively. The complex extensions assumed in STFT and implied in HSA--IMF are fundamentally different and because of this, the quadrature signal $y(t)$ can be different for the two analysis methods. This ultimately results in a different $z(t)$ for each analysis and consequently different instantaneous parameterizations despite both models produce the same $x(t)$ in (\ref{eq:RealObservation}).

\subsection{Real-World Signals} \label{ssec:ExamplesRealWorld}
We provide two examples using real-world audio signals. In Algorithm \ref{alg:ProposedEMD}, we initialize $I=200$, $\alpha=0.95$, and through experimentation, $\beta_k=4$ for all $k$. In addition, we have applied a 1 ms moving-average filter to smooth the IF estimate. 

\begin{figure*}[t]
\centering
  	\begin{minipage}[b]{0.49\linewidth}
  		\centering
  		\subfigure[]{
  		\includegraphics[width = 0.99\linewidth]{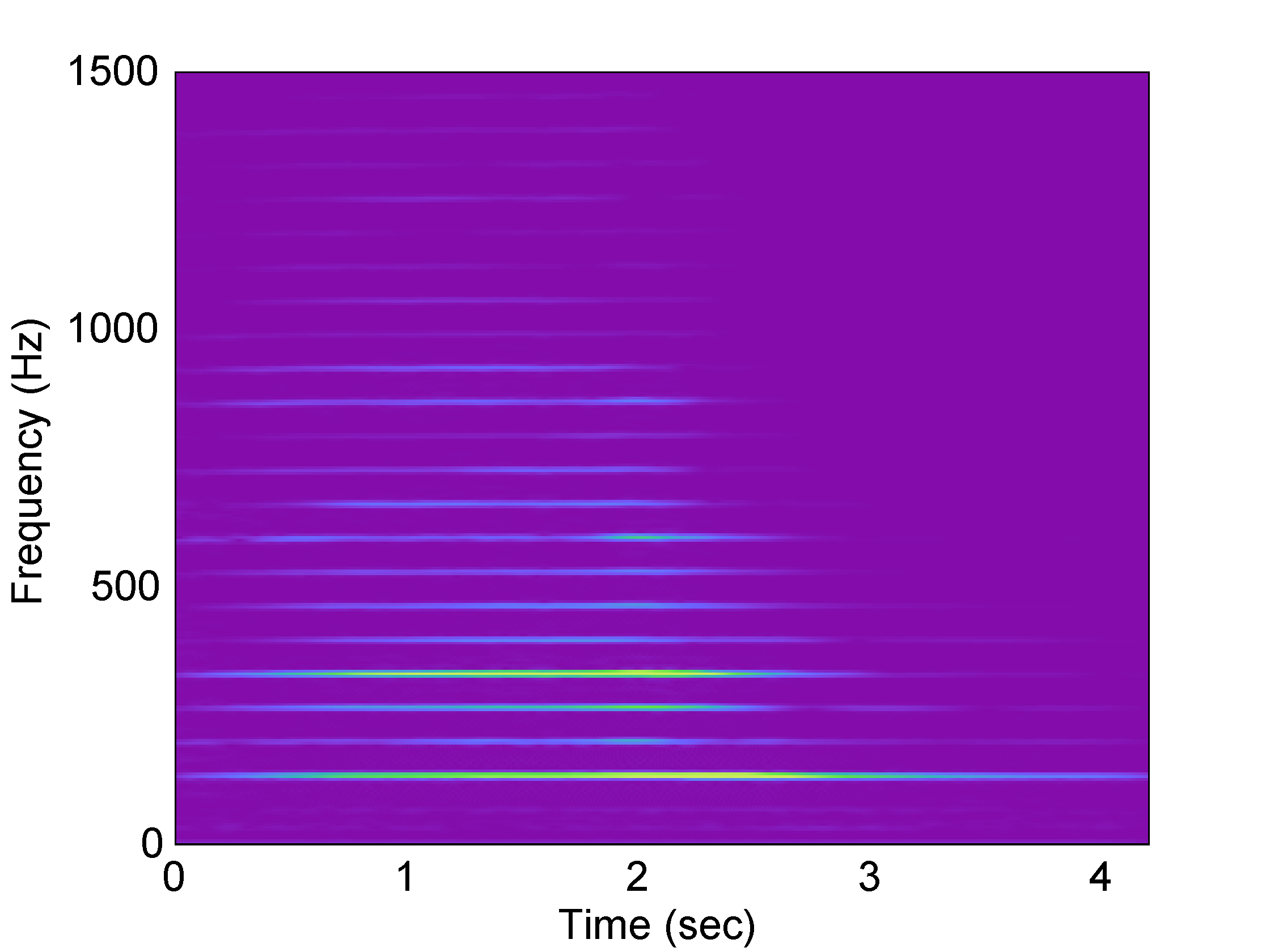}
  		\label{fig:EX3a}
  	}
  	\end{minipage}	
  	\begin{minipage}[b]{0.49\linewidth}
  		\centering	
  		\subfigure[]{
  		\includegraphics[width = 0.99\linewidth]{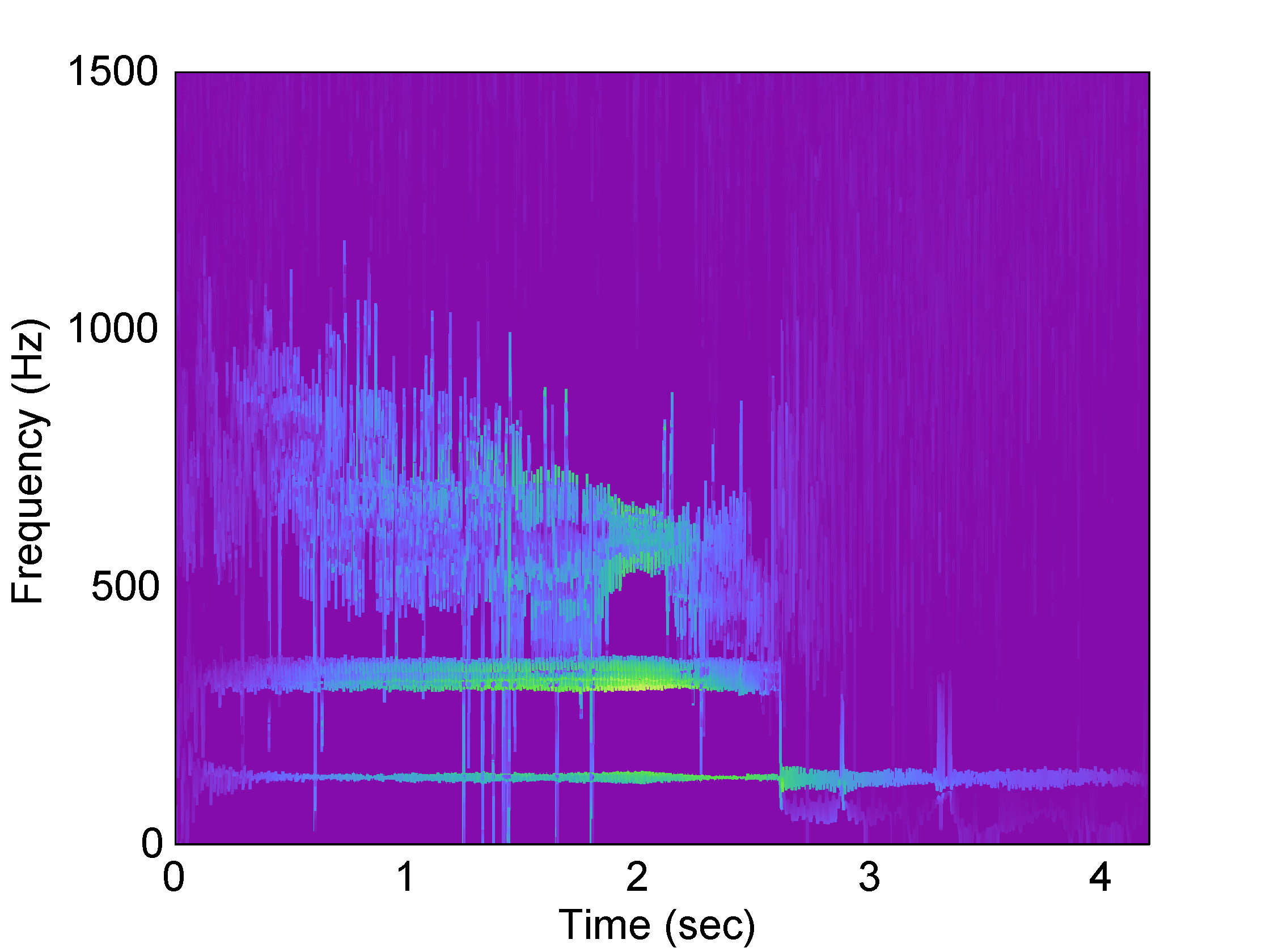}
  		\label{fig:EX3b}
  	}
  	\end{minipage}
  \caption{(a) STFTM and (b) Hilbert spectrum for the cello recording in Example 1 in Subsection \ref{ssec:ExamplesRealWorld}. The lower two components in (b), range in IF from 120-140 Hz and from 300-360 Hz corresponding to the dominant spectral lines in the Fourier spectrum at 133 Hz and 333 Hz. The harmonics above 500 Hz in (a) and the upper component in (b) with IF ranging from 500-1000 Hz partially accounts for the spectral richness of this instrument's note.}
  	\label{fig:cello}
\end{figure*}

In the first example, we analyze a recording ($f_s = 22.050$ kHz) of a single note played on a cello. Fig.~\ref{fig:EX3a} shows the STFTM where we again see classic harmonic structure resulting from the inherent assumption of SHCs. We note the fundamental frequency is approximately 67 Hz (15 spectral lines evenly-spaced over 1000 Hz) and two dominant spectral lines at the second harmonic ($\sim$133 Hz) and at the fifth harmonic ($\sim$333 Hz).  In addition, there is a brief dominant spectral line at $t=2$ s corresponding to the ninth harmonic at $\sim$600 Hz.

Fig.~\ref{fig:EX3b} plots the five components returned by the HSA--IMF algorithm, where we see three dominant components.  The lower two components, range in IF from $\sim$120-140 Hz and from $\sim$300-360 Hz corresponding to the dominant spectral lines in the Fourier spectrum.  The upper component also exhibits significant energy at $t=2$ at $\sim$550-750 Hz corresponding to the brief dominant spectral line, i.e.~ninth harmonic in the Fourier spectrum.

\begin{figure*}[t]
\centering
  	\begin{minipage}[b]{0.49\linewidth}
  		\centering
  		\subfigure[]{
  		\includegraphics[width = 0.99\linewidth]{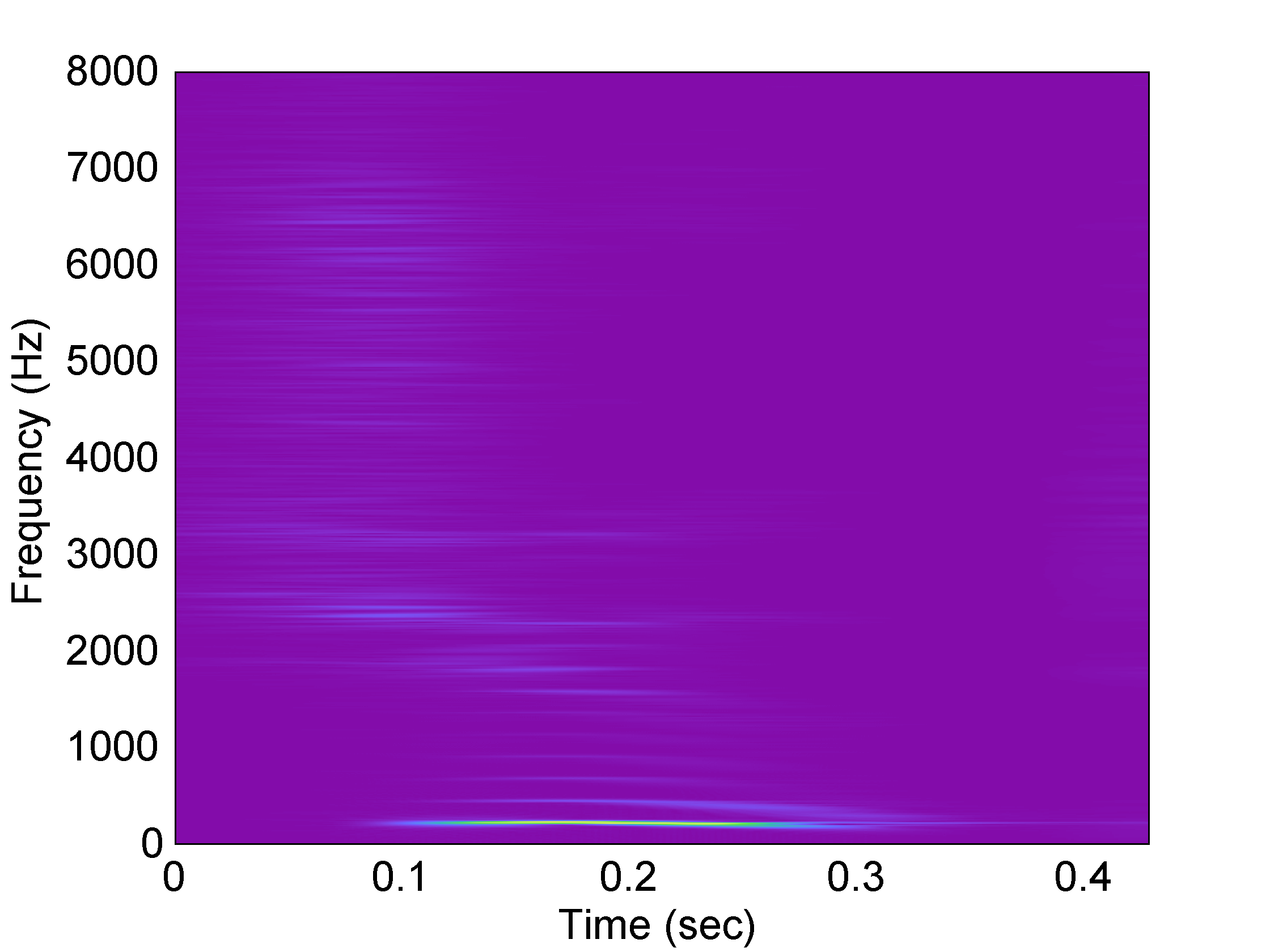}
  		\label{fig:EX4a}
  	}
  	\end{minipage}	
  	\begin{minipage}[b]{0.49\linewidth}
  		\centering	
  		\subfigure[]{
  		\includegraphics[width = 0.99\linewidth]{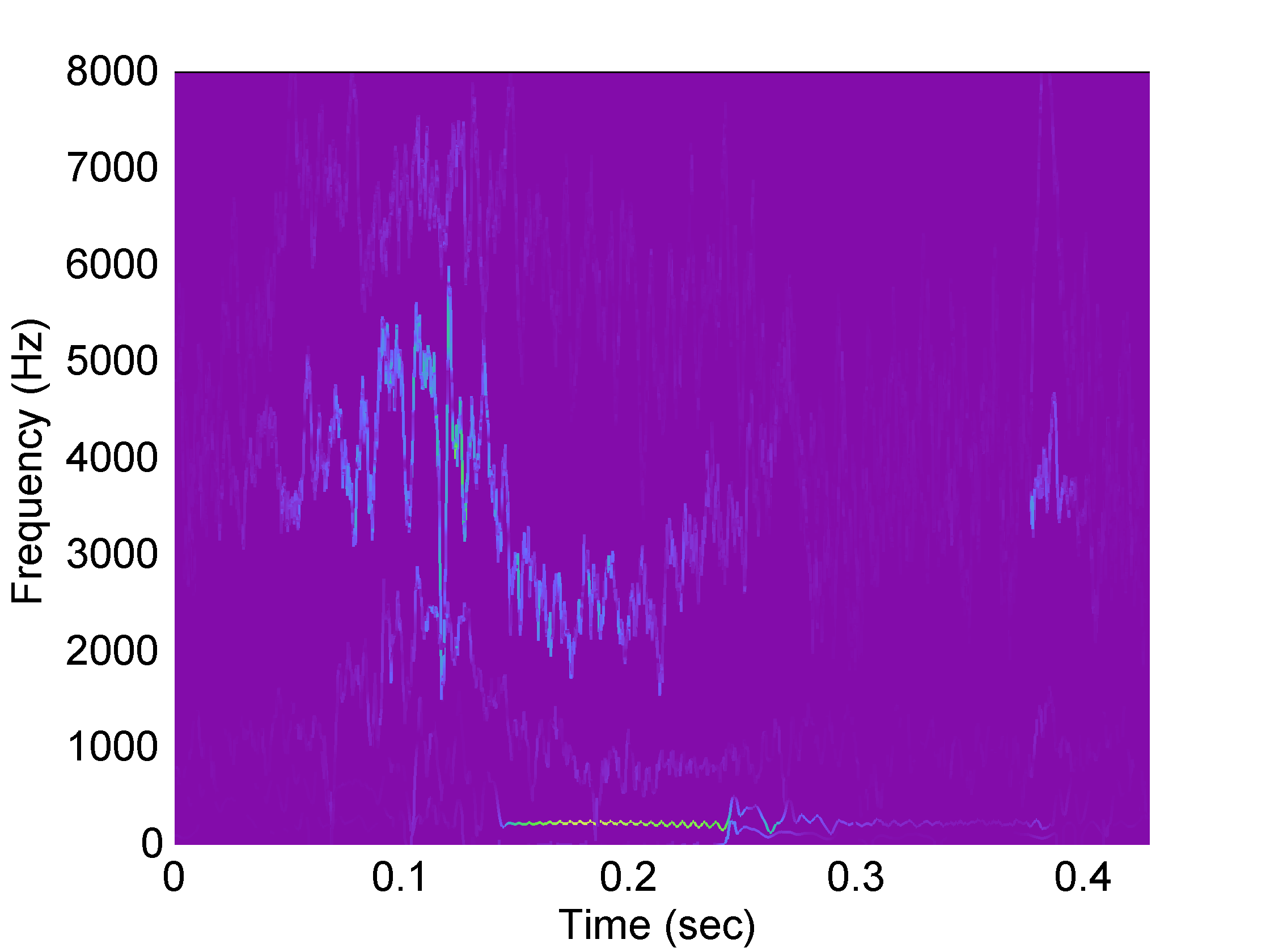}
  		\label{fig:EX4b}
  	}
  	\end{minipage}
  \caption{(a) STFTM and (b) Hilbert spectrum for the speech recording ``shoot'' in Example 2 in Subsection \ref{ssec:ExamplesRealWorld}. The ``SH'' fricative appears in three components with IF ranges of 6000-7000 Hz, 2500-5000 Hz, and 1000-2500 Hz. The vowel ``UW'' is clearly captured in a single component near 230 Hz. Unlike in the Fourier spectrum, the stop ``T'' is clearly captured in the Hilbert spectrum by the first component near $t=0.4$.}
  	\label{fig:Speech}
\end{figure*}

In the second example, we analyze a recording ($f_s = 44.1$ kHz) of the word ``shoot.''  Fig.~\ref{fig:EX4a} shows the STFTM where we see the spectral energy of the fricative ``SH'' over $0 \leq t \leq 0.15$ s, scattered over the range 0 to 8 kHz.  The spectral energy for the vowel ``UW'' over $0.15 \leq t \leq 0.25$, is concentrated at a fundamental of $\sim$230 Hz. The spectral energy for the stop ``T'' over $0.37 \leq t \leq 0.4$ s, is very weak and spread across the band and hence not visible in the plot.

Fig.~\ref{fig:EX4b} plots the five components returned by the HSA--IMF algorithm. The ``SH'' fricative appears in three components with IF ranges of $\sim$6000-7000 Hz (zeroth component), $\sim$2500-5000 Hz (first component), and $\sim$1000-2500 Hz (second component) but is mostly captured in the second component with quickly varying AM and FM. The vowel ``UW'' is clearly captured in a single component near 230 Hz exhibiting some FM variation conjectured to be natural jitter. Unlike in the Fourier spectrum, the stop ``T'' is clearly captured in the Hilbert spectrum by the first component near $t=0.4$.


\section{Discussion} \label{sec:discussion}

\subsection{Resolving Closely-Spaced Components} \label{ssec:ResolvingClosely}   
EMD has been criticized for its inability to resolve closely-spaced components and there have been numerous studies and analyzes on the resolving ability \cite{DeeringMasking, rilling2008one, feldman2009analytical}.  Rilling has investigated EMD's (Algorithm \ref{alg:EMD}) ability to resolve two tones as a function of a relative amplitude parameter and a relative frequency spacing parameter.  This analysis describes regions in this parameter space where EMD returns one or two components \cite{rilling2008one}. As Rilling noted, the goal in EMD is not resolving closely-spaced components but rather resolving components that are suitably matched to an underlying signal model or compatible with assumptions made on the signal model \cite{rilling2008one}.

As an example, consider two infinitely-long tones, $\cos\left(\omega_a t\right)$ and $\cos\left(\omega_b t\right)$, with $\omega_b>\omega_a$.  We can express the sum of these tones as
\begin{subequations}
	\label{eq:beat}
	\begin{align}
	  	\label{eq:beatA}
	    x(t) &= \Re\left\lbrace \exp\left(j\omega_a t\right)+\exp\left(j\omega_b t\right) \right\rbrace \\[-0.1em]
	  	\label{eq:beatB}
	  	     &= \Re\left\lbrace 2\cos\left[\left(\omega_b-\omega_a\right)t/2\right]\exp\left[j\left(\omega_b+\omega_a\right)t/2\right]\right\rbrace.
	\end{align}
\end{subequations}
If $\omega_a$ and $\omega_b$ are sufficiently far apart, both Fourier analysis and EMD will resolve two SHCs as in (\ref{eq:beatA}).  On the other hand, if $\omega_a$ and $\omega_b$ are not sufficiently far apart, EMD will resolve a single IMF as in (\ref{eq:beatB}). As is well known, when $\omega_a$ and $\omega_b$ are closely-spaced, the signal exhibits a beat effect.  In the human auditory system, these closely-spaced tones are not perceived as two distinct tones but rather a single AM tone \cite{oshaughnessy}. As Deering points out, EMD may correspond to the psychoacoustics of human hearing \cite{DeeringMasking}.  As Rilling points out, a decomposition into SHCs may not be an appropriate solution ``...if the aim is to get a representation matched to physics (and/or perception) rather than to mathematics'' \cite{rilling2008one}.

A generalized example of this beat effect was given in Example 2 in Subsection \ref{ssec:ExamplesSynth}.  To clarify the connection to auditory perception of beating, ignore the slow-varying FM and hence the plot lines in Fig.~\ref{fig:SyntEx1} would be horizontal and not sinusoidal. The example, when analyzed with the STFT shows five closely-spaced tones as in Fig.~\ref{fig:SyntEx1a}. However, if these tones are closely spaced then as noted, they may be perceived as a single tone with AM variation. This is demonstrated with an analysis using HSA--IMF and shown in Fig.~\ref{fig:SyntEx1b}.

A similar example regarding an FM signal, was given in Example 1 in Subsection \ref{ssec:ExamplesSynth}, and is essentially the same waveform used in FM synthesis pioneered by Chowning \cite{chowning1977synthesis}.  In Chowning's work, he expressed the FM signal as a superposition of SHCs weighted by Bessel functions and showed rich spectra associated with this signal.  Fig.~\ref{fig:SyntEx2a} illustrates this rich spectra via the presence of multiple harmonics.  In the AM--FM model, such rich spectra may be encapsulated in a single component as illustrated in Fig.~\ref{fig:SyntEx2b} and as has been suggested, a decomposition into SHCs may not be an appropriate solution to describe the underlying signal model.

\subsection{Remarks on Computation} \label{ssec:Computation}
The HSA--IMF algorithm contains a triple-nested loop that incurs significant computation depending on the signal and parameter choices.  Clearly, if the signal has many underlying components, EMD will require more computation due to the extrema searches and interpolations. Unfortunately, there is no way to predict the number of components ahead of time. 

The outermost loop in Algorithm \ref{alg:ProposedEMD} Step \ref{step:outermost} iteratively removes the IMF estimated by tone masking until termination conditions are reached. Hence, there is no way to predict the number of iterations required for termination of the loop. The middle loop, ensemble averages the IMFs returned from tone masking in Algorithm \ref{alg:ProposedEMD} Step \ref{step:getIMF}.  This average is controlled by a fixed number of trials, $I$ (the choice of $I$ is discussed in the next subsection). The inner loop results from the tone masking procedure in Algorithm \ref{alg:ProposedEMD} Step \ref{step:getIMF} calling the sifting algorithm twice, in Algorithm \ref{alg:TM} Step \ref{step:twoSifts}. This in turn calls Algorithm \ref{alg:SIFT}, in which Step \ref{step:Stop} iteratively estimates an IMF. Within this inner loop, the step-size introduced in (\ref{eq:introduceAlpha}) and used in Step \ref{step:Alpha}, also controls the speed at which termination conditions are reached. Of these loops, the innermost loop requires the most computation due to the search for extrema and interpolation.  Taken together, the iterative nature of the outer and inner loops, coupled with computationally complex inner loop can require significant computation depending on the signal length and sample rate.  As has been pointed out, signal oversampling is required for robust estimates of the IMFs further increasing computation.  Finally, IMF demodulation in Step \ref{step:demod} occurs in the outermost loop and does not add significant computational burden.  Much of this computation can occur in parallel \cite{waskito2010parallelizing, chen2010gpgpu, chang2011parallel, waskito2011evaluation}, for example, the ensemble averaging in Algorithm \ref{alg:ProposedEMD} Step \ref{step:getIMF} and the search for extrema in Algorithm \ref{alg:SIFT} Steps \ref{step:findMaxima} and \ref{step:findMinima}.

We have implemented Algorithm \ref{alg:ProposedEMD} in {\sc matlab} where the trials are computed in parallel using a \texttt{parfor} loop and have timed the computation for the synthetic and real-world signal examples in this paper.  Our PC consists of an eight-core AMD FX at 4.01 GHz with 32 GB RAM. The results are given in Table \ref{tbl:benchmarks} where we see for relatively short audio signals, decomposition may require relatively large $\beta$ and $I$ leading to long computation times.

\begin{table}[ht]
\caption{Benchmarks}
\begin{center}
\begin{tabular}{lccccc}
\toprule
\head{Example}             & \head{\boldsymbol{$f_s$} (kHz)} & \head{Duration (s)} &\head{\boldsymbol{$\beta$}} & \head{\boldsymbol{$I$ }}   & \head{Computation Time (s)}\\
\midrule
Synthetic 1         & 44.1        &  1           & 0      & 1      & 1.3\\
Synthetic 2         & 44.1        &  1           & 0      & 1      & 0.9\\
Cello               & 22.05       &  4.21        & 4      & 200    & 320.1\\
Speech              & 44.1        &  0.43        & 4      & 200    & 52.4\\
\bottomrule
\end{tabular}  
\end{center}
\label{tbl:benchmarks}
\end{table}%

\subsection{HSA--IMF Robustness} \label{ssec:Robustness}

The goal in HSA is to obtain an AM--FM decomposition with meaningful interpretation.  Thus, proper identification of the underlying components is required from Algorithm \ref{alg:ProposedEMD}. Two parameters must therefore be carefully chosen: the SNR factor $\beta_k$ and the number of trials $I$.  As a reminder, $\beta_k$ weights the additive masking signal used to mitigate mode mixing and $I$ minimizes, through ensemble averaging, the influence of the additive masking signal on the resulting IMF; both of these parameters appear in Step \ref{step:getIMF}.

In our experience with real-world signals where the underlying signal model is unknown, we begin by selecting $\beta_k=0$ and $I=1$ and visualizing the resulting IMFs by plotting the time-real plane or time-frequency plane. Mode mixing will be evident in the time-real plane when the frequency of the waveform changes abruptly. Mode mixing will also be evident in the time-frequency plane when the IMF is similar to the illustrated IMF within the \textcolor{MyRed}{\protect\dashedrule} frame in Fig.~\ref{fig:EMD_IllustrationA}. If mode mixing is present, we increase $\beta_k$ and $I$.  This process is repeated until reasonable IMFs are obtained keeping in mind the associated computational load for large $I$.  Although this process for decomposition is somewhat heuristic, such refinement is typically present in all time-frequency analyzes \cite{cohen1995time}, e.g.~choice of window length and type for Fourier analysis and mother wavelet selection for wavelet analysis.

The step-size parameter $\alpha$ introduced in (\ref{eq:introduceAlpha}) and used in Algorithm \ref{alg:SIFT} Step \ref{step:Alpha}, is of secondary importance and merely scales the trend which is removed from the signal being sifted.  This scaling is used to minimize the impact of possible overshoot in the trend.  In our experience, selecting $\alpha=0.95$ gives satisfactory performance, noting that lower values lead to additional iterations and more computation.

The termination condition in the sifting algorithm may be too restrictively chosen hence preventing convergence.  In our implementation we include a maximum number of iterations, typically 50, to guarantee an exit.  As a final point, IMFs with excessively large values are omitted from the ensemble.

\section{Future Research in HSA--IMF} \label{sec:further}
In the course of this research, several avenues for further algorithm improvement are apparent.  These include: alternative masking signals, alternate interpolators, and improvements to IMF demodulation. In the HSA--IMF algorithm, we use lowpass filtered noise as the masking signal.  However, for certain signal analysis problems, more sophisticated masking signals have been proposed  \cite{senroy2007improved, senroy2007two, guanlei2009time}.  Ideally, the masking signals would have properties similar to the underlying components which would likely require prior knowledge of the signal model. At the heart of EMD is the iterative estimation of the IMFs which are completely determined by interpolation. As noted earlier, the cubic spline interpolator may be susceptible to overfitting which can lead to poor IMF estimates. Alternatives to cubic spline interpolation have been investigated including B-spline and Akima interpolators \cite{riemenschneider2005b, qin2006envelope, chen2006b, rato2008hht, kopsinis2008improved, wang2010intrinsic, bouchikhi2012multicomponent}, however, these do not appear to offer significant improvements.  Nevertheless, improvements in interpolation may lead to more robust decomposition and demodulation. Finally, IMF demodulation requires estimation of the IF which uses Huang's iterative normalization procedure.  Unfortunately, the required interpolation in the normalization procedure can result in an overfitting of the cubic spline leading to incorrect instantaneous estimates. As noted earlier, changing the cubic spline interpolator in Algorithm \ref{alg:IMFdemod} changes the IMF.  Thus a change of the interpolator in the IA estimator requires the same change in the sifting algorithm.

Although IMFs are latent AM--FM components, there are other classes of AM--FM components that are not IMFs opening up alternate possibilities. The IMF is only important due to its computability via sifting. Hence, it may be possible to define other useful AM--FM components (not defined by \textbf{C1} and \textbf{C2}) corresponding to different assumptions on the form of the component (not \textbf{A1} and \textbf{A2}), that ultimately lead to a replacement of the sifting algorithm as the component tracker, and thus new AM--FM decomposition methods.

\section{Summary} \label{sec:SummaryP3}
In the final part of this paper, we have reported an end-to-end,  solution for the estimation of instantaneous parameters of the AM--FM model assuming IMF components. By leveraging the theory developed in Part II of these papers to interpret and demodulate the results returned by EMD, we have provided a complete numerical method for HSA. We provided examples of HSA--IMF on synthetic signals and argued that the resulting decompositions were more representative of the underlying signal models as compared to conventional Fourier analysis.  Examples of HSA--IMF on real-world signals were shown to allow for alternative and possibly more useful interpretation of the underlying signal model. Finally, we discussed computational aspects of the proposed algorithm. 

\FloatBarrier


\hypertarget{ConcluMark}{}
\bookmark[level=part,dest=ConcluMark]{Conclusions}
\begin{center}\textbf{Conclusions}\end{center}
\setcounter{section}{0}
In this paper, we began by reframing the classic complex extension problem as a Latent Signal Analysis (LSA) problem where the objective is to determine the complex-valued latent signal, $z(t)$ from the real-valued observation, $x(t)=\Re\{z(t)\}$. We used this framework to argue against the assumption of Harmonic Correspondence (HC), and hence the use of Gabor's quadrature method and the Hilbert Transform (HT) for complex extension.  By relaxing the HC assumption, there are many choices for the quadrature and furthermore, many of these complex extensions can be useful in modeling physical phenomena.  

Next, we presented a theory for the Hilbert spectrum framed as generalized LSA problem in which we seek a representation of $z(t)$ consisting of a superposition of latent AM--FM components parameterized by a set of Instantaneous Amplitude (IA)/Instantaneous Frequency (IF) pairs. The use of latent AM--FM components admits many possible forms of the component therefore for a given $x(t)$, there is considerable freedom in the signal model. We presented the analogue of the LSA problem in the frequency domain where we showed that a latent spectrum cannot be uniquely tied to a to the spectrum of the real observation because of the structure imposed by the real operator. We showed that without the HC assumption, there are many choices for the latent spectrum and these spectra will not have in Hermitian symmetry. We proposed a novel 3D visualization of the Hilbert spectrum which plots $\omega(t)$ vs.~$s(t)$ vs.~$t$ and coloring with respect to $|a(t)|$ and allows for the visualization of the instantaneous parameters. We have recast time-frequency analysis and Gabor's analogies to quantum mechanics to an analysis method where uncertainty is in the quadrature signal and not in frequency.  Further, by moving away from Simple Harmonic Components (SHCs) with HC, we allow for a new and powerful way to analyze nonstationary signals. 

We recognized that an Intrinsic Mode Function (IMF) is a latent AM--FM component and leveraged HSA theory to interpret both the IMF and Empirical Mode Decomposition (EMD).  With this interpretation we show that the definition of an IMF unambiguously forces a unique complex extension. Furthermore, we also recognize fundamental problems with the Hilbert-Huang Transform (HHT), i.e.~that the HT is inappropriate for IMF demodulation and proposed an IMF demodulation method that is compatible the the IMF definition. Finally, we utilized the EMD algorithm with our modifications and IMF demodulation to calculate the IA/IF parameters of $x(t)$, thus providing a numerical method for Hilbert Spectral Analysis (HSA).


\hypertarget{AcknowMark}{}
\bookmark[level=part,dest=AcknowMark]{Acknowledgments}
\section*{Acknowledgments}

The authors wish to thank Prof.~Joe Lakey of the Dept.~of Mathematical Sciences at New Mexico State University and Prof.~Antonia Papandreou-Suppappola of the School of Electrical, Computer and Energy Engineering at Arizona State University for reviewing an early draft of this manuscript. 
\newpage

\hypertarget{BibMark}{}
\bookmark[level=part,dest=BibMark]{References}
\bibliographystyle{elsarticle-num}
\bibliography{main.bib}
\end{document}